\documentclass[11pt]{article}
\usepackage[utf8]{inputenc}
\usepackage{amsmath,amssymb,epsfig,bbm}
\usepackage{stmaryrd,mathabx}
\usepackage{comment}
\usepackage{color}
\usepackage[T1]{fontenc}

\usepackage{interval}

\usepackage[textsize=tiny,textwidth=20mm]{todonotes}
\usepackage{enumitem}
\usepackage{varwidth}
\setlist{nolistsep}
\usepackage{hyperref}
\usepackage{mathtools}

\usepackage{amsthm}

\newtheorem{defi}{Definition}
\newtheorem{prop}[defi]{Proposition}
\newtheorem{theo}[defi]{Theorem}
\newtheorem{theofr}[defi]{Théorème}
\newtheorem{conj}[defi]{Conjecture}
\newtheorem{lemm}[defi]{Lemma}
\newtheorem{lemmfr}[defi]{Lemme}
\newtheorem{coro}[defi]{Corollary}

\theoremstyle{definition}
\newtheorem{rema}[defi]{Remark}
\newtheorem{exem}[defi]{Example}
\newtheorem{exems}[defi]{Examples}

\newcommand{\bdefi}{\begin{defi}}
\newcommand{\edefi}{\end{defi}}
\newcommand{\bprop}{\begin{prop}}
\newcommand{\eprop}{\end{prop}}
\newcommand{\btheo}{\begin{theo}}
\newcommand{\etheo}{\end{theo}}
\newcommand{\btheofr}{\begin{theofr}}
\newcommand{\etheofr}{\end{theofr}}
\newcommand{\blemm}{\begin{lemm}}
\newcommand{\elemm}{\end{lemm}}
\newcommand{\blemmfr}{\begin{lemmfr}}
\newcommand{\elemmfr}{\end{lemmfr}}
\newcommand{\brema}{\begin{rema}}
\newcommand{\erema}{\end{rema}}
\newcommand{\bexer}{\begin{exem}}
\newcommand{\eexer}{\end{exem}}
\newcommand{\bexems}{\begin{exems}}
\newcommand{\eexems}{\end{exems}}
\newcommand{\bconj}{\begin{conj}}
\newcommand{\econj}{\end{conj}}
\newcommand{\bcoro}{\begin{coro}}
\newcommand{\ecoro}{\end{coro}}
\newcommand{\dem}{\noindent{\bf Proof. }}

\usepackage{mathrsfs}
\renewcommand\mathcal{\mathscr}

\renewcommand{\H}{{\cal H}}
\newcommand{\I}{{\cal I}}

\newcommand{\M}{{\cal M}}
\newcommand{\N}{{\cal N}}
\newcommand{\OOO}{{\cal O}}

\newcommand{\V}{{\cal V}}

\newcommand{\maths}[1]{{\mathbb #1}}  

\newcommand{\BB}{\maths{B}}
\newcommand{\CC}{\maths{C}}

\newcommand{\HH}{\maths{H}}

\newcommand{\KK}{\maths{K}}

\newcommand{\NN}{\maths{N}}
\newcommand{\OO}{\maths{O}}
\newcommand{\PP}{\maths{P}}
\newcommand{\QQ}{\maths{Q}}
\newcommand{\RR}{\maths{R}}
\newcommand{\SSS}{\maths{S}}

\newcommand{\ZZ}{\maths{Z}}

\newcommand{\aaa}{{\mathfrak a}}
\newcommand{\bbb}{{\mathfrak b}}

\newcommand{\mmm}{{\mathfrak m}}

\newcommand{\ra}{\rightarrow}
\newcommand{\bs}{\backslash}
\newcommand{\sm}{\smallsetminus}

\newcommand{\wt}[1]{{\widetilde{#1}}}

\newcommand{\ga}{\gamma}
\newcommand{\Ga}{\Gamma}

\newcommand{\cqfd}{\hfill$\Box$}

\newcommand{\arcosh}{\operatorname{arcosh}}

\newcommand{\Ax}{\operatorname{Ax}}

\newcommand{\bigO}{\operatorname{O}}

\newcommand{\card}{\operatorname{Card}}

\newcommand{\Fix}{\operatorname{Fix}}

\newcommand{\Heis}{\operatorname{Heis}}

\newcommand{\id}{\operatorname{id}}
\renewcommand{\Im}{\operatorname{Im}}

\newcommand{\nod}{\mathbf d}
\newcommand{\Motonod}{\wh{\mathbf d}}
\newcommand{\Naivenod}{\mathbf d}
\newcommand{\Nr}{\operatorname{{\tt N}}}

\newcommand\Perp{\operatorname{Perp}}

\renewcommand{\Re}{\operatorname{Re}}

\newcommand{\smallo}{\operatorname{o}}
\newcommand{\ssm}{\!\smallsetminus\!}

\newcommand{\Vol}{\operatorname{Vol}}
\newcommand{\vol}{\operatorname{vol}}

\newcommand{\hdr}{{\HH}^2_\RR}
\newcommand{\htr}{{\HH}^3_\RR}

\newcommand{\hnr}{{\HH}^n_\RR}
\newcommand{\hnc}{{\HH}^n_\CC}
\newcommand{\hnk}{{\HH}^n_\KK}

\newcommand{\hnh}{{\HH}^n_\HH}

\newcommand{\PSL}{\operatorname{PSL}}
\newcommand{\SL}{\operatorname{SL}}
\newcommand{\GL}{\operatorname{GL}}

\newcommand{\SLOK}{\operatorname{SL}_{2}(\OOO_K)}
\newcommand{\PSLOK}{\operatorname{PSL}_{2}(\OOO_K)}

\newcommand{\PSLC}{\operatorname{PSL}_{2}(\CC)}

\newcommand{\PSLR}{\operatorname{PSL}_{2}(\RR)}
\newcommand{\SLZ}{\operatorname{SL}_{2}(\ZZ)}
\newcommand{\PSLZ}{\operatorname{PSL}_{2}(\ZZ)}

\newcommand{\PGL}{\operatorname{PGL}}

\newcommand{\flow}[1]{{g^{#1}}}  

\newcounter{const}
\newcounter{fig}

\usepackage[percent]{overpic}
\usepackage{pict2e}

\title{Divergent geodesics, ambiguous closed geodesics \\
  and  the binary additive divisor problem}
\author{Jouni Parkkonen \and Fr\'ed\'eric Paulin}
\date{\today}

\begin{document}
\bibliographystyle{alphanum}
\maketitle
\begin{abstract} 
We give an asymptotic formula as $t\to+\infty$ for the number of
common perpendiculars of length at most $t$ between two divergent
geodesics or a divergent geodesic and a compact locally convex subset
in negatively curved locally symmetric spaces with exponentially
mixing geodesic flow, presenting a surprising non-purely exponential
growth.  We apply this result to count ambiguous geodesics in the
modular orbifold recovering results of Sarnak, and to confirm and
extend a conjecture of Moto\-ha\-shi on the binary additive divisor
problem in imaginary quadratic number fields.
\end{abstract}

\section{Introduction}
\label{sec:intro}
Let $M$ be a noncompact finite volume complete connected negatively
curved locally symmetric good orbifold.

\noindent
\begin{minipage}[t]{4cm}
\vspace{0cm}
\begin{center}
\!\!\includegraphics[width=3.55cm]{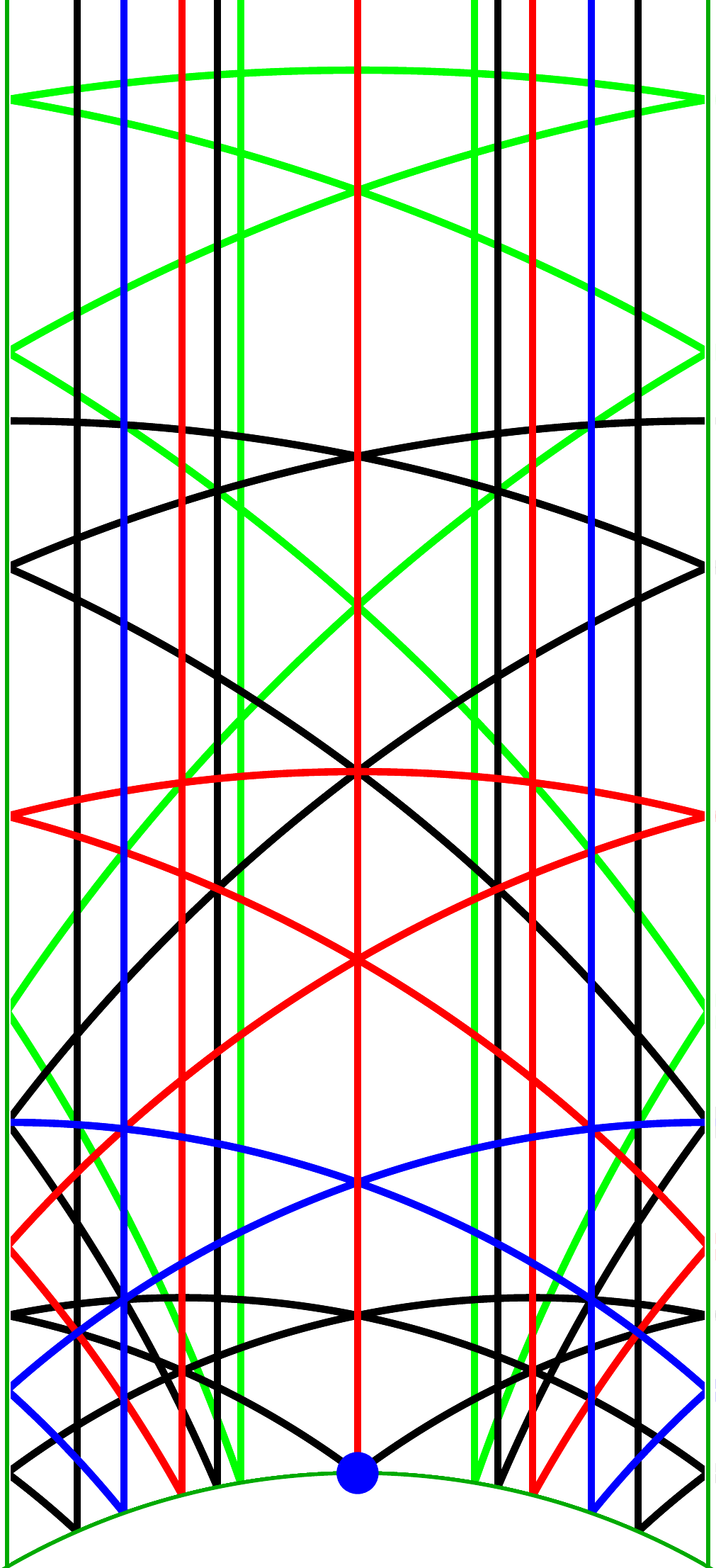}
\end{center}
\vspace*{-0.2cm}
\addtocounter{fig}{1}
{\small Figure \arabic{fig}~:~\hypertarget{divgeod}{Divergent
geo-\\desics in $\Ga_\ZZ\bs\hdr$.}}
\end{minipage}
\hspace{.15cm}
\begin{minipage}[t]{10.75cm} 
\setlength\parindent{15pt}
\noindent ~~~ 
A locally geodesic line $\ell: \RR\to M$ that is a
proper mapping is a {\em divergent geodesic} in $M$. The distribution
of divergent geodesics has been very actively studied in recent
years. We refer for instance to \cite{DavSha18,ParPauSay24} for
equidistribution results of divergent orbits, in the space of lattices
of $\RR^2$ for the first one, in finite volume complete connected
negatively curved good Riemannian orbifolds for the second one. See
for instance \cite{ShaZhe19, DavSha20, SolTam23, DanPauSay23} for
higher rank results.
  
With $\hdr$ the upper halfspace model of the real hyperbolic plane,
Figure \hyperlink{divgeod}{\arabic{fig}} shows some divergent
geodesics in the modular orbifold $\PSLZ\bs\hdr$ (lifted to the usual
fundamental domain of $\Ga_\ZZ=\PSLZ$ with its boundary
identifications).  See Section \ref{sec:divergent} for explanations.

Let $D^-$ and $D^+$ be two properly immersed closed locally convex
subsets of $M$. For instance, $D^-$ and $D^+$ can be the images of two
divergent geodesics in $M$. A {\em common perpendicular}\footnote{See
\cite[\S 2.3]{ParPau17ETDS} for definitions when the boundary of $D^-$
or $D^+$ is not smooth.}  from $D^-$ to $D^+$ is a locally geodesic
path in $M$ starting perpendicularly from $D^-$ and arriving
perpendicularly to $D^+$.  In this paper,
\end{minipage}

\noindent\begin{minipage}{4cm}
\begin{center}
\includegraphics[width=2.5cm]{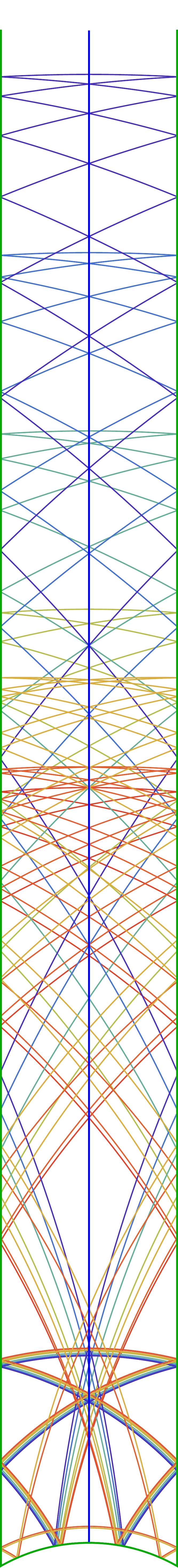}
\end{center}\vspace*{-0.2cm}\addtocounter{fig}{1}
{\small Figure \arabic{fig}~: \hypertarget{commperp}{Common
per\-pen\-dic\-u\-lars in  $\Ga_\ZZ\bs\hdr$\,.
}
}
\end{minipage}
\label{page:tallfigure}
\hspace{.1cm}
\begin{minipage}{10.9cm}
we
prove an effective asymptotic counting result on the set of the common
perpendiculars between the images of two divergent geodesics in $M$.

~~~ For instance, the image $\ell$ in the modular orbifold
$\PSLZ\bs\hdr$ of the imaginary axis in $\hdr$ is a divergent geodesic
in $\PSLZ\bs\hdr$. Figure \hyperlink{commperp}{\arabic{fig}} shows
several common perpendiculars (lifted to the usual fundamental domain
of $\Ga_\ZZ=\PSLZ$ with its boundary identifications) between $\ell$ and
itself (with an extra symmetry, see Remark \ref{rem:ambi2} for
explanations).

~~~ For all $s> 0$, we denote by $\N_{D^-,\,D^+}(s)$ the cardinality
  of the set of common perpendiculars from $D^-$ to $D^+$ with length
  at most $s$, considered with multiplicities (see Equations
  \eqref{eq:defmultipli} and \eqref{eq:defcountingfunct} for
  precisions).  The counting function $\N_{D^-,\,D^+}$ has been
  studied for particular triples $(M,D^-,D^+)$ at least since the
  1940's for example in \cite{Delsarte42, Huber59, Herrmann62,
    Margulis69, EskMcMul93, KonOh11, OhSha12} and \cite{Kim15}.  
    In all these cases,   we have 
    $$
    \lim_{s\to+\infty}\frac{\N_{D^-,\,D^+}(s)}{ c\, e^{\alpha s}}=1
    $$ for 
    constants $c$ and $\alpha$ that depend on $D^-$, $D^+$ and $M$. 
  See \cite{ParPau16LMS} for a more detailed review and an explanation
  of how the above works are related to counting common
  perpendiculars.  As we shall explain in Section \ref{sec:divergent},
  this general purely exponential asymptotic behaviour of
  $\N_{D^-,\,D^+}(s)$ as $s\ra+\infty$ proven in
  \cite[Thm.~1]{ParPau17ETDS} does not apply when $D^-$ or $D^+$ is a
  divergent geodesic.  In this paper, we prove that when $D^-$ or
  $D^+$ is a divergent geodesic, the number of common perpendiculars
  actually no longer has a purely exponential growth in terms of an
  upper bound on their lengths.

\btheo\label{theo:divergentperpintro} Let $M$ be a noncompact finite
volume complete connected real hyperbolic good orbifold of dimension
$n\geq 2$.  Let $D^-$ and $D^+$ be the images of two divergent
geodesics in $M$.  Then there exists a constant $C_{D^-,\,D^+}>0$ such
that as $s\to+\infty$, we have
\[
\N_{D^-,D^+}(s)=\frac{C_{D^-,\,D^+}}{\Vol(M)}\,s^2\, e^{(n-1) \,s}+
\bigO(s\,e^{(n-1) s})\,.
\]
\etheo

~~~ The constant $C_{D^-,\,D^+}$ is made explicit in Theorem
\ref{theo:divergentperp2}. See Theorem \ref{theo:1geoddivreal} for a
version of this theorem when $D^-$ is instead assumed for instance to
be compact, already providing a non-purely exponential growth. The
size of the error term in Theorem \ref{theo:divergentperpintro} is
optimal, as explained below.  See Theorem
\ref{theo:1+2geoddivcompquat} and its following comment for the
version of Theorems \ref{theo:1geoddivreal} and
\ref{theo:divergentperp2} valid for the other locally symmetric
spaces. In Sections \ref{sec:realhypgeom} and
\ref{sec:complexhypgeom}, we give fine results on the lengths of
common perpendiculars that are ending high in Margulis neighborhoods
of the cusps of $M$. These results will be 
\end{minipage}

\bigskip
\noindent 
crucial for the proofs of our geometric main results, Theorems
\ref{theo:1geoddivreal}, \ref{theo:divergentperp2} and
\ref{theo:1+2geoddivcompquat}, that are given in Sections
\ref{sec:divergent} and \ref{sec:complexhypgeom}, and that introduce a
new counting disintegration process, that will explain the non-purely
exponential behavior.

Earlier geometric counting results with growth that is not purely
exponential in negatively curved spaces include the case of closed
geodesics with an upper bound on their length starting with Bowen and
Margulis, see for instance \cite{ParPol83, EskMir11} and
\cite[Coro.~9.15]{PauPolSha15}, as well as the results of
\cite{Vidotto19} and \cite{PeiTapVid20} when the manifold $M$ has
infinite Bowen-Margulis measure or the covering group is of
convergence type. See also \cite{Sarnak07}.

As a first arithmetic application of our geometric counting results,
we recover in Section \ref{sec:ambi} counting results of Sarnak
\cite{Sarnak07} on ambiguous and reciprocal ambiguous\footnote{See
Section \ref{sec:ambi} for definitions} conjugacy classes of primitive
hyperbolic elements in $\PSL_2(\ZZ)$, related to the ambiguous
integral binary quadratic forms of Gauss. For instance, the following
result is due to \cite[Eq.~(12)]{Sarnak07}, with a dif\-ferent
multiplicative constant.

\btheo\label{theo:ambi-intro} The number $\N_{A}(s)$ of ambiguous
primitive closed geodesics in $\PSLZ\bs\hdr$ of length at most $s$
satisfies, as $s\ra+\infty$, that
\[
\N_{A}(s)=\frac{3}{4\,\pi^2}\;s^2\,e^{\frac s2}+\bigO(s\;e^{\frac s2})\,.
\]
\etheo

In the very special case when $M=\PSLZ\bs\hdr$ is the modular
orbifold, we prove in Section \ref{sec:ingham} that Theorem
\ref{theo:divergentperpintro} follows from the asymptotics on the
binary additive divisor problem (see for instance
\cite{Ingham27,Estermann31,HeathBrown79,Motohashi94}). These
arithmetic results produce an asymptotic formula with an additional
term of the form $b_1\, se^s$ with $b_1\ne 0$ and an error term of the
form $\bigO(e^s)$, confirming that the general error term obtained in
Theorem \ref{theo:divergentperpintro} is of the correct order.

Let $K$ be an imaginary quadratic number field, with discriminant
$D_K$, ring of integers $\OOO_K$ and Dedekind zeta function $\zeta_K$.
We denote by $\Naivenod_K:\OOO_K\ssm\{0\}\to\NN$ the (naive) {\em
  number of divisors function} of $\OOO_K$, with
$\Naivenod_K(x)=\card\{d\in \OOO_K\ssm\{0\}:d\mid x\}$ for every $x\in
\OOO_K\ssm\{0\}$. In Section \ref{sec:motohashi}, we use Theorem
\ref{theo:divergentperp2} to prove the following new arithmetic
application.

\btheo\label{theo:mainarithintro} As $X\ra+\infty$, we
have
\[
\sum_{x\in\OOO_K\,\ssm\,\{0,-1\}\,:\; |x|^2\le X}
\Naivenod_K(x)\,\Naivenod_K(x+1) =
\frac{8\,\pi^3}{|D_K|^{3/2}\,\zeta_K(2)}\; X(\ln X)^2+\bigO(X \ln X)\,.
\]
\etheo

In Remark \ref{rem:Moto} \hyperlink{Moto2}{(2)}, we show that this
result confirms a particular case of a conjecture of Motohashi
\cite[p.~277]{Motohashi01} when $K=\QQ(i)$, and gives a generalization
to any imaginary quadratic number field.  It might be possible to
improve our error term given by Theorem \ref{theo:mainarithintro}
using arithmetic methods.  See also \cite{SavVar03} that solves the
special case $K=\QQ(i)$ of Theorem \ref{theo:mainarithintro} with a
less explicit constant.
 
The proof of Theorem \ref{theo:mainarithintro} given in Section
\ref{sec:motohashi} uses arithmetic hyperbolic $3$-manifolds.  Let
$\htr$ be the upper halfspace model of the $3$-dimension real
hyperbolic space, let $M$ be the Bianchi orbifold $\PSL_2(\OOO_K)\bs
\htr$, and let $\ell$ be the image in $M$ of the vertical axis of
$\htr$. Then $\ell$ is a divergent geodesic, and the key idea is to
link (this is new and not immediate) the counting of the binary additive
divisor problem with the counting of common perpendiculars between
$\ell$ and itsef. Then we apply the asymptotic of Theorem
\ref{theo:divergentperpintro}.

We apply Theorems \ref{theo:1geoddivreal} and
\ref{theo:1+2geoddivcompquat} \hyperlink{theocompquat1}{(1)} in
\cite{ParPauHD} to count pairs of Farey neighbours in the rational
field $\QQ$, in quadratic imaginary number fields and in definite
quaternion algebras over $\QQ$.

\medskip
\noindent{\small {\it Acknowledgements:} This research was supported
  by the French-Finnish CNRS IEA PaCap.}

\section{Geometric and measure-theoretic background}
\label{sec:geomback}

Let $\wt M$ be a negatively curved Riemannian symmetric space with
dimension at least $2$ and sectional curvature normalized to have
maximum $-1$.  Then $\wt M$ is isometric to the hyperbolic space
$\HH^n_\KK$ with dimension $n$ over $\KK=\RR,\CC,\HH, \OO$ (with $n=2$
in this last case), with the above normalization of its Riemannian
metric.  Let $\Ga$ be a discrete subgroup of the group
$\operatorname{Isom}(\wt M)$ of isometries of $\wt M$, and let $M=
\Ga\bs \wt M$ be the quotient (complete, connected) locally symmetric
good orbifold. We assume throughout this paper that $M$ is noncompact
and has finite volume. We refer for instance to \cite[\S 2.1]
{BroParPau19} for background on CAT$(-1)$ spaces.

Let $\partial_\infty\wt M$ be the boundary at infinity of $\wt M$, let
$T^1\wt M$ be the unit tangent bundle of $\wt M$, and let $T^1M$ be
the unit tangent bundle of $M$, which identifies as an orbifold with
$\Ga\bs T^1\wt M$. We denote the footpoint maps by $\wt p_\bullet :
T^1\wt M\ra \wt M$ and $p_\bullet:T^1M\ra M$, so that the following
diagram, whose vertical maps are the canonical projections modulo
$\Ga$, is commutative
\begin{equation}\label{eq:ppbullet}
  \begin{array}{ccc}
T^1\wt M & \stackrel{\wt p_\bullet}{\longrightarrow} & \wt M\medskip\\
   ^{\wt p}\downarrow\;\; & & \;\;\downarrow{}^p\smallskip\\
 T^1M & \stackrel{p_\bullet}{\longrightarrow} &  \;M\;.
 \end{array}
\end{equation}
Let $\delta$ be the {\it critical exponent} of $\Ga$, which equals the
topological entropy of the geodesic flow on $T^1M$, see for instance
\cite[Theo.~6.1]{PauPolSha15}.  By for instance \cite[Theo.~4.4
  (i)]{Corlette90}, if $\wt M=\HH^n_\KK$ is the hyperbolic $n$-space
over $\KK=\RR,\CC,\HH, \OO$ (with $n=2$ in this last case), then since
$\Ga$ is a lattice in $\operatorname{Isom}(\wt M)$, we have
\begin{equation}\label{eq:valexpcrit}
  \delta=(\dim_\RR\KK)(n+1)-2\,
\end{equation}

\smallskip
Let $\H$ be a (closed) horoball in $\wt M$, and let $\xi$ be its point
at infinity. For every $x\in\partial \H$, let $t\mapsto x_t$ be the
geodesic line in $\wt M$ starting from the point at infinity $\xi$
such that $x_0=x$, and let $x_{+\infty}\in\partial_\infty \wt
M\ssm\{\xi\}$ be its terminal point at infinity. As defined for
instance in \cite[{Appendix}]{HerPau97} on $\partial_\infty \wt
M\ssm\{\xi\}$ and using the homeomorphism $x\mapsto x_{+\infty}$ from
$\partial \H$ to $\partial_\infty \wt M\ssm\{\xi\}$, the {\it
  Hamenstädt distance} $d_\H$ on $\partial \H$ is defined by
\begin{equation}\label{eq:defihamendist}
\forall \;x,y\in\partial \H,\quad
d_\H(x,y)=\lim_{t\ra+\infty}e^{\frac{1}{2}\,d(x_t,y_t)-t}\,.
\end{equation}
As introduced in \cite[\S 2.1]{HerPau02a}, the {\it cuspidal distance}
$d'_\H$ on $\partial \H$ is defined, for all $x,y\in\partial \H$ by
setting $d'_\H(x,y)$ to be the greatest lower bound of all $r>0$ such
that the horosphere centered at $y_{+\infty}$, at signed distance
$-\ln(2 \,r)$ from $\partial \H$ along the geodesic line $t\mapsto
y_t$, meets the geodesic line $t\mapsto x_t$. When $\wt M=\hnr$, we
have $d'_\H=d_\H$ by loc. cit.. The cuspidal distance is indeed a
distance by loc. cit. since $\wt M$ is a symmetric space, and it is
equivalent to the Hamenstädt distance by \cite[Rem.~2.6]{HerPau02a}.

For every isometry $\ga$ of $\wt M$, for all $x,y\in\partial \H$, we
have $d_{\ga\H}(\ga x,\ga y)= d_\H(x,y)$ and similarly $d'_{\ga\H}(\ga
x,\ga y)= d'_\H(x,y)$.

\blemm\label{lem:majodistham}
For all $x,y\in\partial \H$, if $D$ is the image of the
map $t\mapsto y_t$, then $d_\H(x,y)\leq e^{d(x,D)}$.
\elemm

\dem Let $\xi$ be the point at infinity of $\H$. Let $q$ be the
closest point to $x$ on the geodesic line $D$ and let $\H'$ be the
horoball centered at $\xi$ with $q\in\partial \H'$ (see Figure
\addtocounter{fig}{1}\hyperlink{majdistham}{\arabic{fig}}\addtocounter{fig}{-1}).
Let $p$ be the intersection point with the image of the geodesic line
$t\mapsto x_t$ of the horosphere $\partial \H'$. Recall that two
horospheres centered at the same point at infinity are
equidistant. Since the points $p$ and $q$ are the closest points on
$\partial \H'$ to $x$ and $y$ respectively, we hence have

\begin{center}
  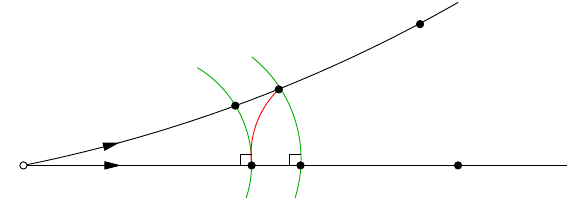
\end{center}
\vspace*{-0.8cm}\addtocounter{fig}{1}
\begin{center}
  {\small Figure \arabic{fig}~: \hypertarget{majdistham}{Majoration
      of the Hamenstädt distance}.}
\end{center}

\noindent $d(y,q)= d(x,p) \leq d(x,q)=d(x,D)$.  By the triangle
inequality, for every $t\geq 0$, we have
\[
d(x_t,y_t)\leq d(x_t,x)+d(x,q)+d(q,y)+d(y,y_t)\leq 2\,d(x,D)+2\,t\,.
\]
The result then follows by the definition \eqref{eq:defihamendist} of
the Hamenstädt distance.  \cqfd

\medskip
We denote by $\|\mmm\|$ the total mass of any finite measure $\mmm$.
Since $\wt M$ is a negatively curved symmetric space and $M$ has
finite volume, there exists up to multiplication by a positive scalar
a unique (measurable) family $(\mu_x)_{x\in\wt M}$ of
Patterson-Sullivan measures on $\partial_\infty\wt M$ for $\Ga$, with
full support, which is actually equivariant under
$\operatorname{Isom}(\wt M)$ (see for instance the book \cite[\S 4]
{BroParPau19} for definitions, since we will be using its
notation). When $\wt M=\hnr$, we normalize these measures so that
$\|\mu_x\|= \Vol(\SSS^{n-1})$ for every $x\in\wt M$. When $\wt M=\hnc$
(respectively $\wt M=\hnh$), we normalize these measures as in
\cite[\S 4]{ParPau17MA} just before Lemma 12 (respectively as in
\cite[\S 7]{ParPau22MPCPS} just before Lemma 7$\cdot$2). The details
of these normalizations are not needed: We will directly use the
computations of the above references that use them.

We denote by $m_{\rm BM}$ the {\it Bowen-Margulis measure} of $M$
associated with this choice of Patterson-Sullivan measures.  Under our
assumption on $M$, it is finite nonzero and mixing, its lift to $\wt
M$ is invariant under $\operatorname{Isom}(\wt M)$ and $m_{\rm BM}$
coincides, up to a positive scalar, with the Liouville measure on
$T^1M$ as well as, when normalized to be a probability measure, with
the measure of maximal entropy for the geodesic flow on $T^1M$ (see
for instance \cite[\S 6, \S 7]{PauPolSha15}). Note that since $M$ is a
finite volume locally symmetric good orbifold with negative curvature,
the geodesic flow $(\flow{t})_{t\in\RR}$ on $T^1M$ is {\it
  exponentially mixing}, that is, there exist $\ell\in\NN$ and
$\kappa,c>0$ such that, denoting by $\|\cdot\|_{\ell,2}$ the Sobolev
$W^{\ell,2}$-norm, for all $\phi,\psi\in C^\ell_c(T^1M)$ and
$t\in\RR$, we have
\[
\Big|\;\int \phi\circ \flow{t}\;\psi \;dm_{\rm BM}-\frac{1}{\|m_{\rm BM}\|}
\;\int \phi \;dm_{\rm BM}\int \psi \;dm_{\rm BM}\;\Big|\leq
c\;e^{-\kappa\,|t|}\;\|\phi\|_{\ell,2}\;\|\psi\|_{\ell,2}\,.
\]
This is well-known by specialists, see for instance \cite{EdwOh21}
when $\wt M=\HH^n_\RR$.  It also follows from the more general
\cite[Theo.~1.1]{Khalil26} with the help of the Appendix B in
\cite{ChaColSch05} in order to obtain the error term in the above
form.

\medskip
By definition, a {\it properly immersed closed locally convex subset}
$D$ of $M$ is the image by the orbifold covering map $\wt M\ra M$ of a
proper nonempty closed convex subset $\wt D$, thereafter called {\em a
  lift} of $D$ in $\wt M$, with stabilizer $\Ga_{\wt D}$ in $\Ga$ such
that the family $(\ga \wt D)_{\ga\in\Ga/\Ga_{\wt D}}$ of subsets of
$\wt M$ is locally finite. We denote by $m(D)$ the order of the
pointwise stabiliser of $\wt D$. We denote by $\partial^1_+ D$ and
$\partial^1_- D$ the {\em outer and inner unit normal bundles} of
$\partial D$ respectively. See \cite{ParPau14ETDS}, generalising
\cite{OhSha12,OhSha13}, for definition when $\partial D$ is not
smooth, or \cite[\S 2.4]{BroParPau19}.

Let us now recall (see \cite[Eq.~(11)]{ParPau14ETDS}) the formula for
the outer/inner skinning measures $\sigma^\pm_D$ of $D$ associated
with the above choice of Patterson-Sullivan measures $(\mu_x)_{x\in\wt
  M}$. Let $p_{\wt D}$ be the closest point projection from $(\wt
M\cup\partial_\infty \wt M) \ssm\partial_\infty\wt D$ to $\wt D$, let
$\beta: (\partial_\infty \wt M\times \wt M\times \wt M)\ra\RR$ be the
Busemann function of $\wt M$, and let $x_*$ be any basepoint in $\wt
M$. Let $\wt\sigma^\pm_{\wt D}$ be the measure on $T^1\wt M$ (with
support contained in $\partial^1_\pm\wt D$) defined as follows: For
every unit normal vector $w\in\partial^1_\pm\wt D$, with $w_\pm$ the
point at $\pm\infty$ of the geodesic line it defines, we have
\begin{equation}\label{eq:defskinning}
  d\wt\sigma^\pm_{\wt D}(w) =
  e^{-\delta\;\beta_{w_\pm}(x_*,\,p_{\wt D}(w_\pm))}\;d\mu_{x_*}(w_\pm)\;.
\end{equation}
Then $\sigma^\pm_D$ is the measure induced by $\wt\sigma^\pm_D$ on
$T^1M$, with support contained in $\partial^1_\pm D$, by the locally
finite $\Ga$-invariant measure $\sum_{\ga\in \Ga/\Ga_{\wt D}}
\ga_*\wt\sigma^\pm_{\wt D}$ on $T^1\wt M$, using the orbifold covering
$\wt p:T^1\wt M\ra T^1M=\Ga\bs T^1\wt M$, see for instance \cite[\S
  2.6]{PauPolSha15}. When $\wt D$ is a geodesic line or a horoball in
$\wt M$, the measure $\wt\sigma^\pm_{\wt D}$ is invariant under the
group of isometries of $\wt M$ preserving $\wt D$, and $\sigma^\pm_D$
is, up to a multiplicative constant, the induced Riemannian measure on
$\partial^1_\pm D$ (see \cite{OhSha12,OhSha13,ParPau17ETDS,ParPau17MA,
  ParPau22MPCPS}, we will give more precision when needed later on).

\medskip
For every horoball $\H$ in $\wt M$, let $\wt\sigma^-_\H$ be the inner
skinning measure of $\H$ (associated with the above choice of
Patterson-Sullivan measures).  Since $\wt M$ is a negatively curved
symmetric space, the group of isometries of $\wt M$ acts transitively
on the set of horoballs of $\wt M$. Furthermore, the group of
isometries of $\wt M$ preserving $\H$ acts transitively on $\partial
\H$ and leaves $\wt \sigma^-_\H$ invariant since $M$ has finite
volume. Let $B_{d'_\H}(x,r)$ be the ball of radius $r>0$ and center
$x\in\partial\H$ for the cuspidal distance $d'_\H$ on $\partial
\H$. Let us define
\begin{equation}\label{eq:defiXi}
\Xi_{\wt M}=\wt \sigma^-_\H\big(\wt p_\bullet^{\;-1}(B_{d'_\H}(x,1))\big)\,,
\end{equation}
which depends neither on the horoball $\H$ in $\wt M$ nor on the point
$x\in\partial\H$. A computation of this constant will be given in
Equation \eqref{eq:Xirealhyp} and in Lemma \ref{lem:computXinonreel}.

\section{A lemma in real hyperbolic geometry}
\label{sec:realhypgeom}

Let
\[
\hnr=\big(\;\{(x,y)\in\RR^{n-1}\times\RR:y>0\}, \quad ds^2_{\hnr}=
\frac{1}{y^2} (dx^2+dy^2)\;\big)
\]
be the upper halfspace model of the real hyperbolic space of dimension
$n$ (with constant sectional curvature $-1$). Recall that
$\partial_\infty \hnr=(\RR^{n-1}\times\{0\})\cup\{\infty\}$, that
\begin{equation}\label{eq:defiHinftyreal}
\H_\infty=\{(x,y)\in\hnr:y\geq 1\}
\end{equation}
is a horoball in $\hnr$ centered at $\infty$, and that the geodesic
line in $\hnr$ from $(x,0)\in \partial_\infty \hnr$ to $\infty$,
through the horosphere $\partial \H_\infty= \{(x,y)\in\hnr : y=1\}$ at
time $t=0$, is the map $t\mapsto (x,e^{t})$. Furthermore, the map
$x\mapsto (x,1)$, from $\RR^{n-1}$ endowed with the standard Euclidean
distance to $\partial\H_\infty$ endowed with the Hamenstädt distance
$d_{\H_\infty}$, is an isometry.

\blemm\label{lem:distancetodivline} Let $\H$ be a horoball in $\hnr$
and let $\wt\ell$ be a geodesic line in $\hnr$ that enters $\H$
perpendicularly at $\wt\ell(0) \in\partial\H$.

\smallskip\noindent \hypertarget{lemhypreeli}{(i)} Let $\wt\ell'$ be a
geodesic line in $\hnr$ that exits $\H$ perpendicularly at
$\wt\ell'(0) \in\partial\H$.  For every $t\ge 0$, we have
\[
d(\wt\ell'(t),\wt\ell\,)=t+\ln d_\H(\wt\ell'(0),\wt\ell(0))+\ln 2
+\bigO\big(d_\H(\wt\ell'(0),\wt\ell(0))^{-2}e^{-2t}\big)\,.
\]

\smallskip\noindent\hypertarget{lemhypreelii}{(ii)} Let $D$ be a
closed convex subset of $\hnr$ disjoint from $\H$ and let
$x_0\in\partial\H$ be the closest point to $D$ in $\H$. If
$d_\H(x_0,\wt\ell(0))\geq 1$, then
\[
d(D,\wt\ell\,)=d(D,\H)+\ln d_\H(x_0,\wt\ell(0))+\ln 2+
\bigO\big(d_\H(x_0,\wt\ell(0))^{-2}e^{-2\,d(D,\,\H)}\big)\,.
\]
Furthermore, if $d_\H(x_0,\wt\ell(0))\geq 2$, then the closest point
to $D$ on  $\wt\ell$ belongs to $\H$.
\elemm

\begin{center}
  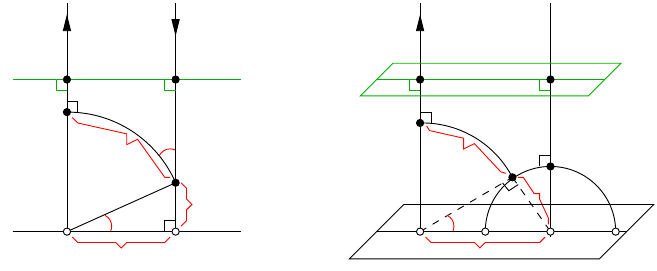
\end{center}
\vspace*{-0.8cm}\addtocounter{fig}{1}
\begin{center}
{\small Figure \arabic{fig}~: \hypertarget{lemrealhyp}{Real hyperbolic distance formulas to geodesic lines}.}
\end{center}
\dem \hyperlink{lemhypreeli}{(i)} Since two geodesic lines meeting
perpendicularly a horosphere have one common point at infinity, the
geodesic lines $\wt\ell$ and $\wt\ell'$ are coplanar. By symmetry, we
may assume that $n=2$, that $\H=\H_\infty$ and that there exists $a>0$
such that for every $t\in\RR$, we have $\wt\ell(t)=(0,e^t)$ and
$\wt\ell'(t)= (a,e^{-t})$. Note that then
$a=d_\H(\wt\ell(0),\wt\ell'(0))$. Recall that in a right-angled
hyperbolic triangle with one vertex at infinity, with finite opposite
side length $u$ and acute angle $\alpha$, the {\it angle of
  parallelism formula} gives $\cosh u=
\frac{1}{\sin\alpha}$.\footnote{See for instance \cite[Thm.~7.9.1
  (ii)]{Beardon83}.}  Hence (see Figure
\hyperlink{lemrealhyp}{\arabic{fig}} on the left), we have as wanted
\begin{align*}
d(\wt\ell'(t),\wt\ell\,)&=\arcosh\frac{\sqrt{a^2+e^{-2t}}}{e^{-t}}
= \ln(\sqrt{a^2e^{2t}+1}+a\,e^{t})\\
&=\ln(a\,e^{t}(1+\sqrt{1+a^{-2}e^{-2t}}\,))=
t+ \ln a+\ln 2+\bigO(a^{-2}e^{-2t})\,.
\end{align*}

\smallskip\noindent\hyperlink{lemhypreelii}{(ii)} Let $x_{\wt\ell}\in
D$ and $p_{x_{\wt\ell}}$ be the endpoints of the common perpendicular
between $D$ and $\wt\ell$, and let $x_{\H}\in D$ be the closest point
in $D$ to $\H$. Let $\wt\ell'$ be the geodesic line through $x_{\H}$
exiting $\H$ at time $0$ perpendicularly at $x_0$. Let $t=d(D,\H)$ and
$a=d_\H(x_0,\wt\ell(0))$. Since
\[
d(x_{\wt\ell},\wt\ell\,)=d(D,\wt\ell\,)\leq d(x_\H,\wt \ell\,)
=t+ \ln a+\ln 2+\bigO(a^{-2}e^{-2t})
\]
by Assertion \hyperlink{lemhypreeli}{(i)} applied with the above
$\wt\ell'$, we only need to prove a similar lower bound on
$d(x_{\wt\ell},\wt\ell)$.

The union of the geodesic lines perpendicular to $\wt\ell'$ at $x_\H$
is a (totally geodesic) hypersurface $S$ that separates $D$ and $\H$
(and is a supporting hyperplane of $D$). Replacing $x_{\wt\ell}$ by
the intersection point of the geodesic segment $[p_{x_{\wt\ell}},
  x_{\wt\ell}]$ with $S$ does not increase $d(x_{\wt\ell}, \wt\ell)$.
Let $D'$ be the geodesic line in $S$ through $x_\H$ and $x_{\wt\ell}$.
Up to rotating $D'$ around $\wt\ell'$ until it lies in the copy of the
hyperbolic plane containing $\wt\ell$ and $\wt\ell'$, which does not
increase $d(x_{\wt\ell},\wt\ell\,)$, we may assume that $\wt\ell$,
$\wt\ell'$ and $D'$ are coplanar. We may assume that $n=2$,
$\H=\H_\infty$, $\wt\ell(t)= (0,e^t)$, $\wt\ell'(t)= (0,e^{-t})$ as in
Assertion \hyperlink{lemhypreeli}{(i)}, so that $D'$ is the geodesic
line with points at infinity $(a-e^{-t},0)$ and $(a+e^{-t},0)$. Since
by assumption $a\geq 1>e^{-t}$, the point $x_{\wt\ell}$ is the vertex
at its right angle of the right-angled Euclidean triangle with other
vertices $(0,0)$ and $(a,0)$. When furthermore $a\geq 2$, then
$a-e^{-t}\geq 1$, hence $p_{x_{\wt\ell}}$ belongs to $\H_\infty$. The
Euclidean length of the side between $x_{\wt\ell}$ and $(a,0)$ is
equal to $e^{-t}$, since $D$ is the open Euclidean halfcircle centered
at $(a,0)$ with Euclidean radius $e^{-t}$.  Again applying the angle
of parallelism formula (see Figure
\hyperlink{lemrealhyp}{\arabic{fig}} on the right), we have
as wanted
\begin{align*}
  d(x_{\wt\ell},\wt\ell\,)&=\arcosh\big(\frac{a}{e^{-t}}\big)=
  \ln(a\,e^{t}+\sqrt{a^2e^{2t}-1})\\
&=\ln(a\,e^{t}(1+\sqrt{1-a^{-2}e^{-2t}}\,))=
t+ \ln a+\ln 2+\bigO(a^{-2}e^{-2t})\,.\quad\Box
\end{align*}

\section{Common perpendiculars of divergent geodesics}
\label{sec:divergent}

Before giving precisely the counting function whose asymptotic we
will study in this paper, let us recall the structural properties of
the noncompact finite volume complete connected negatively curved good
orbifold $M=\Ga\bs\wt M$. Let $\operatorname{Par}_\Ga$ be the subset
of $\partial_\infty\wt M$ consisting of the fixed points of the
parabolic elements of $\Ga$.  The {\it set of cusp} of $M$ is the
finite set $\Ga\bs \operatorname{Par}_\Ga$ of $\Ga$-orbits of
parabolic fixed points of $\Ga$, whose elements are denoted by $e_1,
\dots, e_m$.

Let $(H_\xi)_{\xi\in\operatorname{Par}_\Ga}$ be a $\Ga$-equivariant
family of (closed) horoballs with pairwise disjoint interiors, with
$H_\xi$ centered at $\xi$ for every $\xi$ in $\operatorname{Par}_\Ga$,
which is {\it precisely invariant}: If the intersection
$\ga\stackrel{\circ}{\H_\xi} \cap \stackrel{\circ}{\H_{\xi'}}$ is
nonempty, then $\ga \xi=\xi'$.  For every $i\in\llbracket
1,m\rrbracket$, the closure $\V_{e_i}$ of $\Ga\bs\big(\bigcup_{\xi\in
  e_i}\stackrel{\circ}{H_\xi}\big)$ is called the {\it Margulis
  neighbourhood} of the cusp $e_i$.  The closure of $M\ssm\bigcup_{1\leq
  i\leq m} \V_{e_i}=\Ga\bs\big(\wt M\ssm
\bigcup_{\xi\in\operatorname{Par}_\Ga}H_\xi\big)$ is a compact subset
of $M= \Ga\bs\wt M$.

\medskip
We understand the locally geodesic lines $\ell$ in $M$ in the orbifold
sense. They are possibly not injective maps from $\RR$ to $M$ with
multiplicities. Since the fixed point set of an isometry of $\wt M$ is
totally geodesic, the orbifold stabilizer of a positive length
subsegment of $\ell$ is equal to the orbifold stabilizer of the whole
$\ell$. Specific examples when $\ell$ is not injective are the
following ones.

We say that a locally geodesic line $\ell$ in $M$ (or its image) is
{\em reversible} if it has a lift $\wt\ell:\RR\ra \wt M$ such that an
element of $\Ga$ interchanges the two endpoints at infinity of the
geodesic line $\wt\ell$.  Let $\iota_{\rm rev} (\ell(\RR))=1$ if
$\ell$ is reversible, and $\iota_{\rm rev} (\ell(\RR))=2$ otherwise.
We say that $\ell$ (or its image) is {\em strongly reversible} if
there is such an element of order $2$.

When $\wt M=\hdr$ and $\Ga= \PSL_2(\ZZ)$, a locally geodesic line in
$M=\Ga\bs \wt M$ is reversible if and only if it is strongly
reversible if and only if it passes through $\Ga \cdot i$.  In this
context, (strongly) reversible closed geodesics are called {\em
  reciprocal geodesics} because of their connection with reciprocal
integral binary quadratic forms, see \cite{Sarnak07} for more details.
See \cite{Sarnak07}, \cite{ErlSou24} and \cite{ErlParPau} for counting
and equidistribution results of reciprocal and strongly reversible
closed geodesics in negatively curved spaces, and Remark
\ref{rem:ambi2} for an example.

Recall that a locally geodesic line $\ell: \RR\to M$ that is a proper
mapping is a {\em divergent geodesic} in $M$.  By the above
description of $M$, a locally geodesic line $\ell: \RR\to M$ is a
divergent geodesic in $M$ if and only if there are times
$t_-,t_+\in\RR$ with $t_-\leq t_+$ at which $\ell$ meets at a right
angle the boundary of two Margulis neighbourhoods $\V_-$ and $\V_+$ of
cusps of $M$, that we refer to as the {\em initial} and {\em terminal
  Margulis neighbourhoods} of $\ell$. They are possibly equal, for instance  when
$M$ has only one cusp or when $\ell$ is reversible.  The images
of the subrays $\ell|_{\interval[open left] {-\infty}{\,t_-}}$ and
$\ell|_{\interval[open right]{t_+} {\,+\infty}}$ of $\ell$ are
contained in $\V_-$ and $\V_+$.  The image $D_\ell$ of $\ell$ is a
properly immersed closed locally convex subset in $M$, and we have
$\partial^1_+ D_\ell=\partial^1_- D_\ell$ and $\iota_*
\sigma^-_{D_\ell} = \sigma^+_{D_\ell}$, where $\iota :v\mapsto -v$ is
the time reversal map of $T^1M$.

\medskip
\noindent{\bf Examples: } Since the set of parabolic fixed points of
$\PSLZ$ is $\QQ\cup\{\infty\}$, the divergent geodesics in the modular
orbifold $\PSLZ\bs\hdr$ are the images of the vertical geodesic lines
starting from $\infty$ and ending at a rational point $\frac{p}{q}$
with $p,q$ coprime. Figure \addtocounter{fig}{-3}\hyperlink{commperp}{\arabic{fig}}\addtocounter{fig}{3} in the introduction gives all
the divergent geodesics ending at a rational point with positive
denominators at most $6$. Here are three further pictures, with
divergent geodesics defined by the rational numbers $3/8$, $31/80$ and
$3/10$ from the left to the right.  The last one, passing through $i$,
is strongly reversible. Following the path of each geodesic in the quotient
orbifold requires to use the boundary identifications $z\mapsto z+1$
and $z\mapsto -\frac 1z$ of the usual fundamental domain
$\{z\in\CC:-\frac{1}{2}\le\Re z\le \frac{1}{2}, |z|\ge 1\}$ of $\PSLZ$.
\begin{center}
 \includegraphics[width=4.2cm]{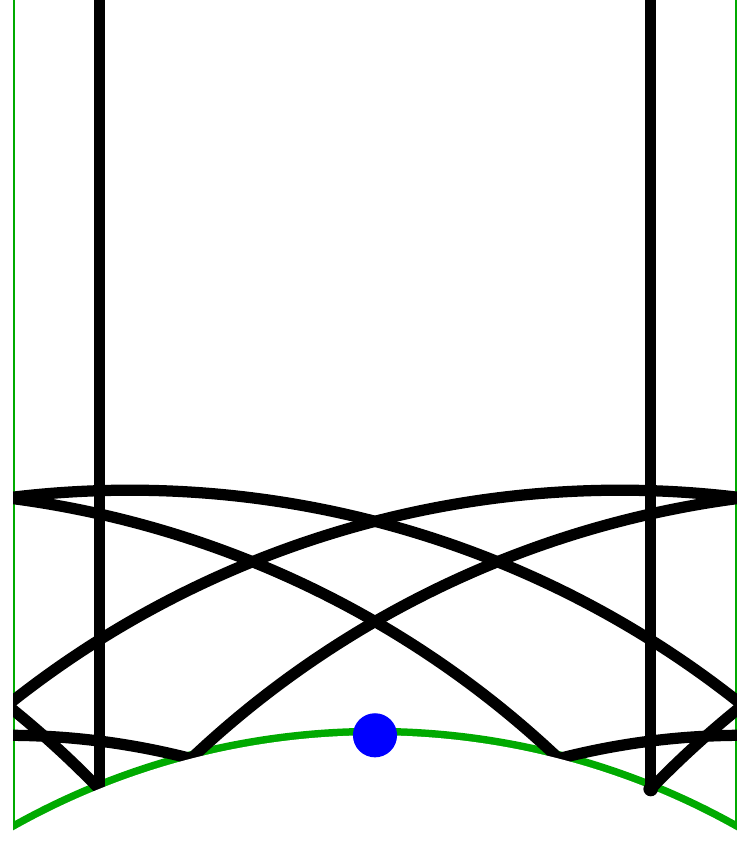}\hspace*{0.5cm}
  \includegraphics[width=4.2cm]{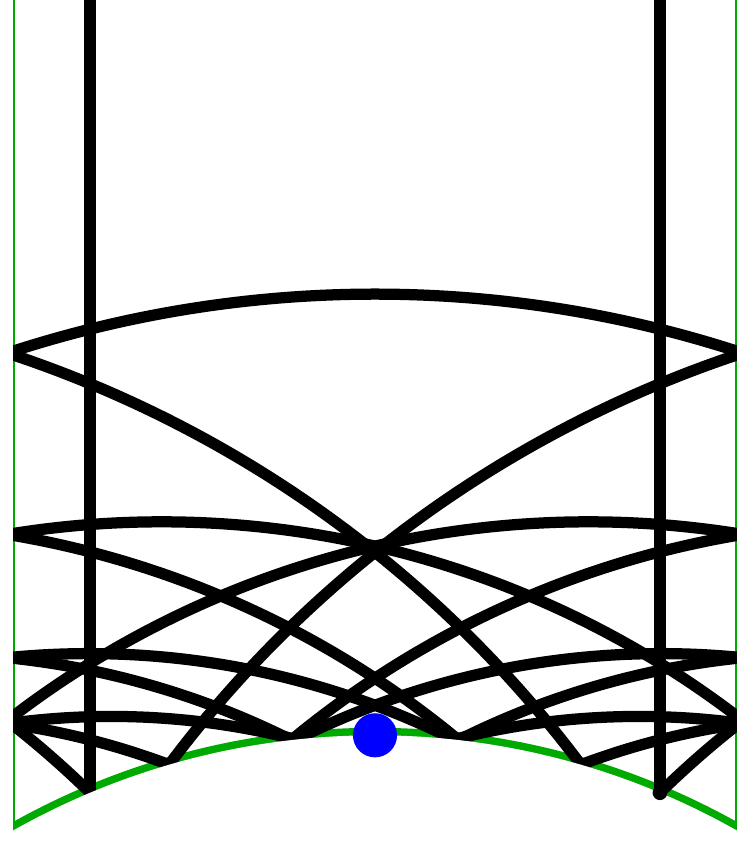}\hspace*{0.5cm}
  \includegraphics[width=4.2cm]{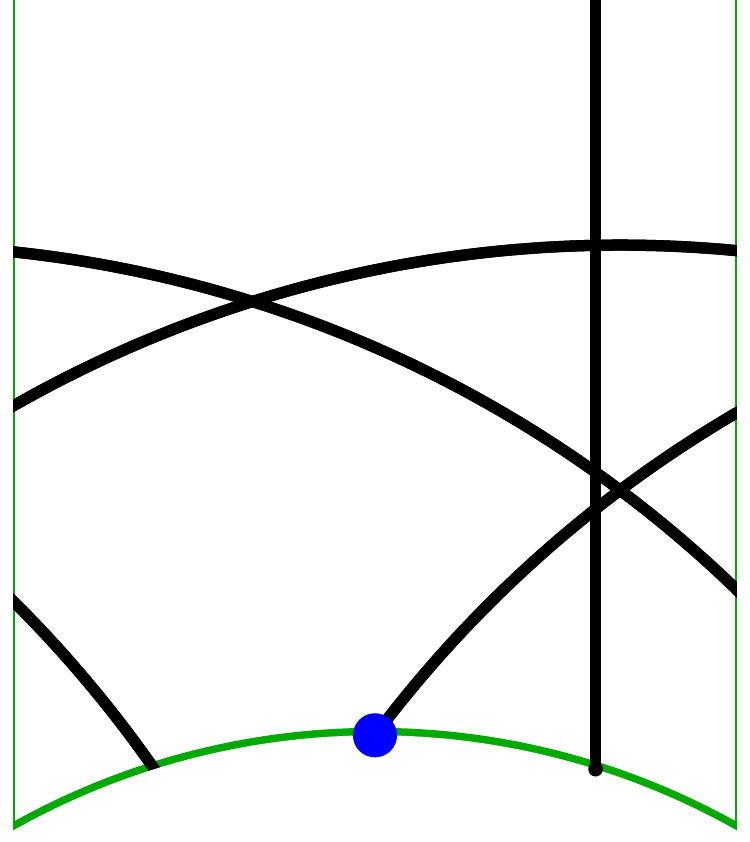}
\end{center}
\vspace*{-0.6cm}\addtocounter{fig}{1}
\begin{center}
{\small Figure \arabic{fig}~: \hypertarget{moredivgeod}{Divergent geodesics on $\PSLZ\bs\hdr$  with vertical lift ending at $3/8$, $31/80$ and
    $3/10$}.}
\end{center}
\medskip
Let $D^-$ and $D^+$ be two properly immersed closed locally convex
subsets of $M$, with lifts $\wt D^-$ and $\wt D^+$. For all $\ga,\ga'$
in $\Ga$, the convex sets $\ga \wt D^-$ and $\ga' \wt D^+$ have a
common perpendicular (as defined in the introduction) if and only if
their closures $\overline{\ga \wt D^-}$ and $\overline{\ga' \wt D^+}$
in $\wt M\cup\partial_\infty \wt M$ do not intersect. This common
perpendicular $\alpha_{\ga \wt D^-,\,\ga' \wt D^+}$ starts from $\ga
\wt D^-$ at time $t=0$ with unit tangent vector $\dot{\alpha}_{\ga\wt
  D^-,\,\ga'\wt D^+}(0)\,$ and ends in $\ga' \wt D^+$ at time $t=d(\ga
\wt D^-,\ga' \wt D^+)$ with unit tangent vector $\dot{\alpha}_{\ga\wt
  D^-,\,\ga'\wt D^+}(d(\ga \wt D^-,\ga' \wt D^+))\,$.  Its {\em
  multiplicity}\footnote{See \cite[\S 3.3]{ParPau17ETDS} and \cite[\S
  12.1]{BroParPau19} for precisions.}  is
\begin{equation}\label{eq:defmultipli}
m_{\ga\wt D^-,\ga'\wt D^+}=
\frac 1{\card(\ga\Ga_{\wt D^-}\ga^{-1}\cap\ga'\Ga_{\wt D^+}{\ga'}^{-1})}\,,
\end{equation}
which equals $1$ when $\Ga$ acts freely on $T^1\wt M$ (for instance
when $\Ga$ is torsion-free).  Note that for every $\ga''\in\Ga$, we
have
\[
\ga''\alpha_{\ga \wt D^-,\,\ga' \wt D^+}=
\alpha_{\ga''\ga \wt  D^-,\,\ga''\ga' \wt D^+}\quad\text{and}\quad
m_{\ga''\ga\wt D^-,\ga''\ga'\wt D^+}=m_{\ga\wt D^-,\ga'\wt D^+}\,.
\]

Recall that $\wt p:T^1\wt M\ra T^1M=\Ga\bs T^1\wt M$ is the canonical
projection. Let $\Omega^-$ and $\Omega^+$ be measurable subsets of
$\partial^1_+ D^-$ and $\partial^1_- D^+$, and let $\wt\Omega^-=
\partial^1_+ \wt D^-\cap\wt p^{\,-1}(\Omega^-)$ and $\wt\Omega^+=
\partial^1_-\wt D^+\cap\wt p^{\,-1}(\Omega^+)$ be the subsets of all
elements of $\partial^1_+ \wt D^-$ and $\partial^1_- \wt D^+$ mapping
to $\Omega^-$ and $\Omega^+$ by $\wt p$. Note that $\wt\Omega^-$ and
$\wt\Omega^+$ are invariant under $\Ga_{\wt D^-}$ and $\Ga_{\wt D^+}$
respectively. The counting function $\N_{\Omega^-,\,\Omega^+}$ is
defined\footnote{See \cite[page 86]{ParPau17MA} for precisions.} by
\begin{equation}\label{eq:defcountingfunct}
\N_{\Omega^-,\,\Omega^+}:t\mapsto
\sum_{\substack{\Ga(\ga\Ga_{\wt D^-},\,\ga'\Ga_{\wt D^+})\in
  \Ga\bs((\Ga/\Ga_{\wt D^-})\times (\Ga/\Ga_{\wt D^+}))\\
  \overline{\ga \wt D^-}\,\cap \,\overline{\ga' \wt D^+}\,=
    \emptyset,\; d(\ga\wt D^-,\ga'\wt D^+)\leq t \\
    \dot{\alpha}_{\ga\wt D^-,\,\ga'\wt D^+}(0)\,\in\, \ga\wt\Omega^-,\;
    \dot{\alpha}_{\ga\wt D^-,\,\ga'\wt D^+}(d(\ga\wt D^-,\ga'\wt D^+))\,
    \in\, \ga'\wt\Omega^+}} m_{\ga\wt D^-,\ga'\wt D^+}\,.
\end{equation}
where $\Ga$ acts diagonally on $\Ga/\Ga_{\wt D^-}\times\Ga/\Ga_{\wt
  D^+}$. In order to simplify the notation, let $\N_{D^-,\,D^+}
=\N_{\partial^1_+ D^-,\,\partial^1_- D^+}$, $\N_{\Omega^-,\,D^+}
=\N_{\Omega^-,\partial^1_- D^+}$ and $\N_{D^-,\,\Omega^+}
=\N_{\partial^1_+ D^-,\,\Omega^+}$.

\medskip
As mentionned in the introduction, the general purely exponential
asymptotic theorem on the counting function $\N_{D^-,\,D^+}(s)$ as
$s\ra+\infty$ proven in \cite[Thm.~1]{ParPau17ETDS}, which requires
the finiteness of the skinning measures of $D^-$ and $D^+$, does not
apply when $D^-$ or $D^+$ is the image of a divergent
geodesic. Indeed, if $\ell$ is a divergent geodesic, then the skinning
measure of its image has infinite total mass.  Theorems
\ref{theo:1geoddivreal} and \ref{theo:divergentperp2} below show that,
when $D^-$ or $D^+$ is a divergent geodesic, the growth of the
counting function $\N_{D^-,\,D^+}$ is no longer purely exponential.

\btheo\label{theo:1geoddivreal}
Let $M$ be a noncompact finite volume complete connected real
hyperbolic good orbifold.  Let $D^+$ and $D^-$ be nonempty properly
immersed closed locally convex subsets of $M$.  Assume that $D^-$ has
nonzero finite (outer) skinning measure and that $D^+$ is the image of
a divergent geodesic in $M$.  Then as $s\ra+\infty$, we have
\[
\N_{D^-,D^+}(s)=
\frac{\Ga(\frac{n}{2})\,\iota_{\rm rev}(D^+)\,\|\sigma^+_{D^-}\|}
  {2^n\,\sqrt{\pi}\;\Ga(\frac{n+1}{2})\,m(D^+)\,\Vol M}\;
  s\, e^{(n-1) \,s} + \bigO(e^{(n-1) s})\,.
\]
\etheo

\dem This proof gives more detail than might seem necessary, and uses
as much as possible the general notation of the beginning of Sections
\ref{sec:geomback} and \ref{sec:divergent}, in order to serve proving
Theorem \ref{theo:divergentperp2} (also when $M$ is real hyperbolic)
and Theorem \ref{theo:1+2geoddivcompquat} (when $\wt M\neq \hnr$). We
believe that this process will be easier for the reader.

Let $\ell$ be a divergent geodesic in $M$ whose image is $D^+$.  Let
$\V_-$ and $\V_+$ be the initial and terminal Margulis neighbourhoods
of $\ell$. Let $t_-$ be the first exit time of $\ell$ from $\V_-$ and
let $t_+$ be the last entry time of $\ell$ into $\V_+$. Let $\wt\ell$
be a lift of $\ell$ in $\wt M$, and let $\wt D^+$ be the image of
$\wt\ell$. For simplicity, let $m^+=m(D^+)$ and $\iota_{\rm rev}^+=
\iota_{\rm rev}(D^+)$.  Let
\begin{align*}
  \Omega_- & =\big\{v\in\partial^1_-D^+:
  p_\bullet(v)\in\ell(\,\interval[open]{-\infty}{t_-}\,)\big\}\,,\\
  \Omega_0\, & =\big\{v\in\partial^1_-D^+:
  p_\bullet(v)\in\ell(\interval{t_-}{t_+})\big\}\quad\textrm{and}\\
  \Omega_+ & =\big\{v\in\partial^1_-D^+:
  p_\bullet(v)\in\ell(\,\interval[open]{t_+}{+\infty}\,)\big\}\;.
\end{align*}
We denote by $\H=H_{\wt\ell(+\infty)}$ the horoball of the equivariant
family $(H_\xi)_{\xi\in\operatorname{Par}_\Ga}$ with point at infinity
$\wt\ell(+\infty)$, that is a lift of $\V_+$. Let $\wt D^-$ be a lift
of $D^-$.

\medskip\noindent{\bf Case \hypertarget{Case1}{1}. } Let us first
assume that $\ell$ is not reversible.  Let us prove that the
subsets $\Omega_-$ and $\Omega_+$ are disjoint.  Assume for a
contradiction that there exists an element $\ga$ of $\Ga$ mapping an
element in $\partial^1_-\wt D^+$ with footpoint in
$\wt\ell(\,\interval[open] {-\infty}{t_-}\,)$ to another element with
footpoint in $\wt\ell (\,\interval[open]{t_+}{+\infty} \,)$.  Then
$\ga$ would map a point $\wt\ell(s_-)$ with $s_-\in \interval[open]
{-\infty}{t_-}$ to a point $\wt\ell(s_+)$ with $s_+\in\interval[open]
{t_+}{+\infty}\,$. Since the equivariant family
$(H_\xi)_{\xi\in\operatorname{Par}_\Ga}$ is precisely invariant, the
element $\ga$ would also map $\wt\ell(-\infty)$ to
$\wt\ell(+\infty)$. Therefore the restriction of $\ga$ to $\wt D^+$
would be the central symmetry with respect to the unique fixed point
of $\ga$ on $\wt D^+$, thereby exchanging its two points at
infinity. Thus $\ell$ would be reversible, a contradiction.

In particular, since $\partial^1_-D^+=\Omega_-\cup\Omega_0 \cup
\Omega_+$, for every $s\geq 0$, we have
\begin{equation}\label{eq:NDDOmegapm}
  \N_{D^-,\Omega_-}(s)+\N_{D^-,\Omega_+}(s)\leq
\N_{D^-,D^+}(s)\leq\N_{D^-,\Omega_-}(s)+\N_{D^-,\Omega_0}(s)+\N_{D^-,\Omega_+}(s)\;.
\end{equation}

We shall prove that as $s\ra+\infty$, we have
\begin{equation}\label{eq:NOmegaplus}
  \N_{D^-,\Omega_+}(s)=
  \frac{\|\sigma^+_{D^-}\|\;\Xi_{\wt M}}{2^\delta\,m^+\,\|m_{\rm BM}\|}\;
  s\,e^{\delta\,s}+\bigO(e^{\delta\,s})\,.
\end{equation}
By the independence property on the horoball in the definition of
$\Xi_{\wt M}$ in Equation \eqref{eq:defiXi}, the same proof replacing
the horoball $\H= H_{\wt\ell(+\infty)}$ by the horoball
$H_{\wt\ell(-\infty)}$ (that is a lift of $\V_-$) gives that as
$s\ra+\infty$, we will have
\begin{equation}\label{eq:NOmegamoins}
  \N_{D^-,\Omega_-}(s)=
  \frac{\|\sigma^+_{D^-}\|\;\Xi_{\wt M}}{2^\delta\,m^+\,\|m_{\rm BM}\|}\;
  s\,e^{\delta\,s}+\bigO(e^{\delta\,s})\,.
\end{equation}
Let $D_0=\ell([t_-,t_+])$, which is a compact nonempty properly
immersed locally convex subset of $M$, hence has a nonzero finite
inner skinning measure. By \cite[Theo.~1]{ParPau17ETDS}, as
$s\ra+\infty$, we therefore have
\begin{equation}\label{eq:NOmegazero}
  \N_{D^-,\Omega_0}(s)\leq \N_{D^-,D_0}(s)= \bigO(e^{\delta\,s})\,.
\end{equation}
Thus by Equations \eqref{eq:NDDOmegapm}, \eqref{eq:NOmegaplus},
\eqref{eq:NOmegamoins} and \eqref{eq:NOmegazero}, and since
$\iota_{\rm rev}^+=2$ when $\ell$ is not reversible, as
$s\ra+\infty$, we will have
\begin{equation}\label{eq:NDDabstrait}
  \N_{D^-,D^+}(s)=
  \frac{\iota_{\rm rev}^+\,\|\sigma^+_{D^-}\|\;\Xi_{\wt M}}
       {2^\delta\,m^+\,\|m_{\rm BM}\|}\;
  s\,e^{\delta\,s}+\bigO(e^{\delta\,s})\,.
\end{equation}

By \cite[Prop.~20 (2)]{ParPau17ETDS}, since $M$ is real hyperbolic in
the assumptions of Theorem \ref{theo:1geoddivreal}, since the
Patterson-Sullivan measures have been normalized in Section
\ref{sec:geomback} to have total mass $\Vol(\SSS^{n-1})$, the metric
measured space $(\partial \H, d_\H, (\wt p_\bullet)_*\wt\sigma^-_\H)$
is isomorphic to $\RR^{n-1}$ endowed with its usual Euclidean distance
and with $2^{n-1}$ times its Lebesgue measure. Furthermore, the
Hamenstädt distance $d_\H$ is equal to the cuspidal distance $d'_\H$
since $M$ is real hyperbolic. Hence by the definition of $\Xi_{\wt M}$
in Equation \eqref{eq:defiXi}, we have
\begin{equation}\label{eq:Xirealhyp}
  \Xi_{\wt M}= 2^{n-1}\,\Vol(\BB_{n-1})= 
  \frac{2^{n-1}\,\pi^{\frac{n-1}{2}}}{\Ga(\frac{n+1}{2})}\,.
\end{equation}
By \cite[Prop.~20 (1)]{ParPau17ETDS}, we have
\begin{equation}\label{eq:PP17ETDS20_1}
  \|m_{\rm BM}\|=2^{n-1}\,\Vol(\SSS^{n-1})\,\Vol(M)=
  2^{n-1}\,\frac{2\,\pi^{\frac{n}{2}}}{\Ga(\frac{n}{2})}\,\Vol(M)\,.
\end{equation}
We also have $\delta=n-1$. Hence Theorem \ref{theo:1geoddivreal} will
follow from Equation \eqref{eq:NDDabstrait}, once we have proven
Equation \eqref{eq:NOmegaplus}.
Up to changing the parametrisation of the geodesic line $\wt \ell$ by
a translation, we may assume that $t_+=0$ to simplify the
notation. We now start the work on the sum defining $\N_{D^-,\Omega_+}
(s)$ given by Equation \eqref{eq:defcountingfunct} in order to prove
Equation \eqref{eq:NOmegaplus}.

\medskip\noindent{\bf Step \hypertarget{Step1}{1}. } The first step is
to simplify the set of indices in this sum.

The map from $\Ga\bs((\Ga/\Ga_{\wt D^-})\times (\Ga/\Ga_{\wt D^+}))$
to $\Ga_{\wt D^-}\bs\Ga/\Ga_{\wt D^+}$ defined by
\[
\Ga(\ga\Ga_{\wt D^-},\ga'\Ga_{\wt D^+})\mapsto
\Ga_{\wt D^-}\ga^{-1}\ga'\Ga_{\wt D^+}
\]
is a bijection, whose inverse is $[\ga]=\Ga_{\wt D^-}\ga\Ga_{\wt D^+}
\mapsto\Ga(\ga^{-1}\Ga_{\wt D^-},\Ga_{\wt D^+})$, by an immediate
checking. In order to simplify the notation, let 
\[
z_{\ga\wt D^-,\,\wt  D^+} =
{\alpha}_{\ga\wt D^-,\,\wt D^+}(d(\ga\wt D^-,\wt D^+))\,.
\]
Since $\wt D^+$ is not reversible, we have $\partial^1_-{\wt
  D^+} \cap \wt p^{-1}(\Omega_+)=\wt p_\bullet^{\;-1}(\wt \ell(\,
\interval[open] {0}{+\infty}\,))$. Thus by Equation
\eqref{eq:defcountingfunct} and by a change of variable $\ga\mapsto
\ga^{-1}$, for every $s\geq 0$, we have
\begin{align}
  \N_{D^-,\,\Omega_+}(s)&=
  \sum_{\substack{[\ga]\in\Ga_{\wt D^-}\bs\Ga/\Ga_{\wt D^+}\,:\;
    0<d(\ga^{-1} \wt D^-,\wt D^+)\leq s \\ 
    \dot{\alpha}_{\ga^{-1}\wt D^-,\,\wt D^+}(d(\ga^{-1} \wt D^-,\wt D^+))\, \in\,
    \wt p_\bullet^{\;-1}(\wt\ell(\,\interval[open]{0}{+\infty}\,))}}
  m_{\ga^{-1}\wt D^-,\wt D^+} \nonumber\\ &=
  \sum_{[\ga]\in\Ga_{\wt D^+}\bs\Ga/\Ga_{\wt D^-} \,:\;
    0<d(\ga \wt D^-,\wt D^+)\leq s,\; z_{\ga\wt D^-,\,\wt D^+}\,
    \in\, \wt\ell(\,\interval[open]{0}{+\infty}\,)} m_{\ga\wt D^-,\wt D^+}\,.
  \label{eq:defcountfunct1}
\end{align}

\medskip\noindent{\bf Step \hypertarget{Step2}{2}. } The second step
in the proof of Equation \eqref{eq:NOmegaplus} is to prove that the
contribution to the above sum defining $\N_{\wt D^-,\Omega_+}(s)$ of
the indices with multiplicities different from $1$ is negligible.

By Equation \eqref{eq:valexpcrit}, the critical exponent of a positive
codimension totally geodesic subspace of $\wt M$ is at most
$\delta-1$. Note that the stabilizer $\Ga_{\wt D^+}$ of $\wt D^+$ in
$\Ga$ is finite since $\wt\ell(+\infty)$ is a parabolic fixed point of
$\Ga$, hence no loxodromic element in $\Ga$ preserves $\wt D^+$. The
fixed point set of a nontrivial isometry with finite order of $\wt M$
is a totally geodesic subspace with positive codimension. Hence as
$[\ga]$ ranges over $\Ga_{\wt D^+}\bs\Ga/\Ga_{\wt D^-}$ with $d(\ga\wt
D^-,\wt D^+)>0$, the common perpendiculars between $\ga\wt D^-$ and
$\wt D^+$ with multiplicity $m_{\ga\wt D^-,\wt D^+}\neq 1$ are
contained in finitely many positive codimension totally geodesic
subspaces of $\wt M$. By Equation \eqref{eq:defmultipli}, the
multiplicities $m_{\ga\wt D^-,\wt D^+}$ are at most $1$. Hence as
$s\ra+\infty$, by the same proof as the one that follows, we have
\begin{align}
&\sum_{[\ga]\in\Ga_{\wt D^+}\bs\Ga/\Ga_{\wt D^-}:\;0<d(\ga\wt
    D^-,\wt D^+)\leq s,\; m_{\ga\wt D^-,\wt D^+}\neq 1}
  m_{\ga\wt D^-,\wt D^+} \nonumber\\\leq \;&
  \card\{[\ga]\in\Ga_{\wt D^+}\bs\Ga/\Ga_{\wt D^-}:0<d(\ga\wt
  D^-,\wt D^+)\leq s,\;\; m_{\ga\wt D^-,\wt D^+}\neq 1\}
  \nonumber\\ =\;&\bigO(s\,e^{(\delta-1)s})=\bigO(e^{\delta\,s})\,.
\label{eq:decaymultnon1}
\end{align}
Hence by Equation \eqref{eq:defcountfunct1}, we have
\begin{equation}
  \N_{D^-,\,\Omega_+}(s)=\card\big\{\ga\in\Ga_{\wt D^+}\bs\Ga/\Ga_{\wt D^-}:
   \begin{array}{l}0<d(\ga \wt D^-,\wt D^+)\leq s\\z_{\ga\wt D^-,\,\wt D^+}
    \in \wt\ell(\,\interval[open]{0}{+\infty}\,)\end{array}\big\}+
    \bigO(e^{\delta\,s}) \,.
  \label{eq:defcountfunct2}
\end{equation}
The map from $(\Ga_\H\bs\Ga/\Ga_{\wt D^-})\times \Ga_\H$ to $\Ga_{\wt
  D^+} \bs\Ga/\Ga_{\wt D^-}$ defined by $([\ga],\beta)\mapsto \Ga_{\wt
  D^+} \beta \ga\Ga_{\wt D^-}$ is onto, and will be used in Step
\hyperlink{Step6}{6}, Equation \eqref{eq:defcountfunct3}, in order to
disintegrate the counting in Equation \eqref{eq:defcountfunct2} above
a counting in $\Ga_\H\bs\Ga/\Ga_{\wt D^-}$, the point being that the
images of $\H$ and $\wt D^-$ in $M$ both have finite skinning
measures. But this requires some preliminary work.

\medskip\noindent{\bf Step \hypertarget{Step3}{3}. } The third step in
the proof of Equation \eqref{eq:NOmegaplus} is to define two
exceptional finite subsets $F$ in $\Ga_{\H}\bs\Ga/ \Ga_{\wt D^-}$ and
$F'$ in $\Ga_\H$ whose contribution to the counting that will occur in
Step \hyperlink{Step6}{6}, Equation \eqref{eq:defcountfunct3}, will be
proven to be negligible in Step \hyperlink{Step5}{5}.

The quotient space $\Ga_\H\bs \partial \H$ is compact, since the point
at infinity of the horoball $\H$ is the parabolic fixed point
$\wt\ell(+\infty)$ of $\Ga$. Since the family $(\ga'\wt
D^-)_{\ga'\in\Ga/\Ga_{\wt D^-}}$ is locally finite, the subset $F$ of
elements $[\ga]\in\Ga_{\H}\bs\Ga/ \Ga_{\wt D^-}$ such that the closed
convex subsets $\H$ and $\ga \wt D^-$ are not disjoint is finite. For
every $\ga\in\Ga$ such that $\Ga_{\H} \ga \Ga_{\wt D^-}\notin F$, we
denote by $p_\ga\in\partial \H$ and $p'_\ga\in\partial (\ga\wt D^-)$
the two endpoints of the common perpendicular between $\H$ and $\ga
\wt D^-$ (see Figure
\addtocounter{fig}{1}\hyperlink{figstep3}{\arabic{fig}}\addtocounter{fig}{-1}).
Note that $p_\ga$ and $p'_\ga$ are unchanged by the multiplication on
the right of $\ga$ by an element of $\Ga_{\wt D^-}$.

Again since $\Ga_\H\bs \partial \H$ is compact, there exists a
constant \addtocounter{const}{1}$c_{\arabic{const}}>0$ (for instance
the radius of a closed ball in $(\partial \H,d'_\H)$ with center
$\wt\ell(0)$ which maps onto $\Ga_\H\bs \partial \H\,$) such that we
may fix from now on a representative $\ga$ of every $[\ga]\in
\Ga_{\H}\bs \Ga/\Ga_{\wt D^-}\ssm F$ (by multiplying it on the left by
an element of $\Ga_\H$) so that $d'_\H(p_\ga,\wt\ell(0)) \leq
c_{\arabic{const}}$.

As $M$ is assumed to be real hyperbolic in Theorem
\ref{theo:1geoddivreal}, we define $\KK=\RR$ and $c_\KK=2$.
Since the isometric action of $\Ga_\H$ on $(\partial \H,d'_\H)$ is
discrete, there exists a finite subset $F'$ of $\Ga_\H$ such that for
every $\beta\in\Ga_\H\ssm F'$, we have the inequality $d'_\H(\wt
\ell(0),\beta\wt \ell(0)) \geq c_\KK+c_{\arabic{const}}$.  By the
triangle inequality, for all $[\ga]\in \Ga_{\H}\bs \Ga/\Ga_{\wt
  D^-}\ssm F$ and $\beta\in\Ga_\H\ssm F'$, we thus have
\begin{equation}\label{eq:minocKK}
d'_\H(\wt\ell(0),\beta p_\ga)\geq d'_\H(\wt\ell(0),\beta \wt\ell(0))-
d'_\H(\beta\wt\ell(0),\beta p_\ga) \geq (c_\KK+c_{\arabic{const}})-
c_{\arabic{const}}=c_\KK\,.
\end{equation}
Therefore since $M$ is real hyperbolic, by the second claim of Lemma
\ref{lem:distancetodivline} \hyperlink{lemhypreelii}{(ii)} applied
with $D=\beta\ga\wt D^-$ and $x_0=\beta p_\ga$, the point
$z_{\beta\ga\wt D^-,\,\wt D^+}$ belongs to the positive subray
$\wt\ell(\,\interval[open]{0}{+\infty}\,)$.

Figure
\addtocounter{fig}{1}\hyperlink{figstep3}{\arabic{fig}}\addtocounter{fig}{-1}
represents in red the common perpendicular between $\beta\ga\wt D^-$
and $\wt D^+$ in the generic situation when $[\ga]\notin F$ and
$\beta\notin F'$.

\begin{center}
  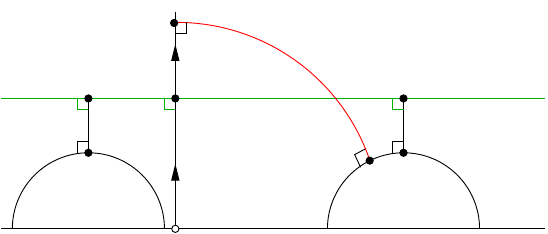
\end{center}
\vspace*{-0.8cm}\addtocounter{fig}{1}
\begin{center}
  {\small Figure \arabic{fig}~: \hypertarget{figstep3}{Generic
      common perpendiculars between $D^-$ and $D^+$}.}
\end{center}

\medskip\noindent{\bf Step \hypertarget{Step4}{4}. } In this rather
independent fourth step in the proof of Equation \eqref{eq:NOmegaplus},
we study the orbital growth of the parabolic subgroup $\Ga_\H$.

For every $[\ga]\in F$, we fix any representative $\ga\in\Ga$ of
$[\ga]$, and we fix $p_\ga\in \partial \H$ such that $d(\ga\wt
D^-,\partial \H)=d(\ga\wt D^-, p_\ga)$. See Figure
\addtocounter{fig}{1}\hyperlink{figstep5}{\arabic{fig}}\addtocounter{fig}{-1}
when $\ga\wt D^-$ meets $\partial \H$, though $\ga\wt D^-$ could be
contained in the interior of $\H$. For all $[\ga]\in\Ga_{\H}\bs\Ga/
\Ga_{\wt D^-}$ and $t\geq 0$, let
\[
\Phi_{[\ga]}(t)=\card\{\beta\in\Ga_\H:d'_\H(\wt\ell(0),\beta p_\ga)\leq t\}\,.
\]
The group of isometries of $\wt M$ preserving $\H$ acts transitively
on $\partial\H$ and preserves the measure $(\wt p_\bullet)_*
\wt\sigma^-_{\H}=(\wt p_\bullet)_*\wt\sigma^+_{\H}$ on $\partial \H$
(see \cite[Prop.~20 (3)]{ParPau17ETDS}, \cite[Lem.~12
  (iv)]{ParPau17MA}, \cite[Lem.~7$\cdot$2]{ParPau22MPCPS} for
details). Furthermore, using the definition of $\Xi_{\wt M}$ in
Equation \eqref{eq:defiXi}, it satisfies the following homogeneity
property: for every $x\in \partial \H$ and $r>0$, we have
\[
(\wt p_\bullet)_*\wt\sigma^-_{\H}(B_{d'_\H}(x,r))=\Xi_{\wt M}\,r^\delta\,.
\]
Recall that $\Ga_\H$ is a uniform lattice in the isometry group of
$(\partial\H,d'_\H)$. By the standard Gauss counting argument
(covering the ball with center $\wt\ell(0)$ and radius $r$ by
translates by elements of $\Ga_\H$ of a given compact fundamental
domain with measure zero boundary and measure $\|\sigma^-_{\V_+}\|$
for the measure $(\wt p_\bullet)_*\wt\sigma^-_{\V_+}$, with a
$\bigO(\cdot)$ which is uniform in $[\ga]$ since $p_\ga$ varies in a
compact subset of $\partial\H$ when $[\ga]\notin F$ and since $F$ is
finite when $[\ga]\in F$, we have
\begin{equation}\label{eq:horocount}
\Phi_{[\ga]}(t)=\frac{\Xi_{\wt M}}{\|\sigma^-_{\V_+}\|}\;t^\delta
+\bigO(t^{\delta-1})\,.
\end{equation}

\medskip\noindent{\bf Step \hypertarget{Step5}{5}. } In this fifth
step in the proof of Equation \eqref{eq:NOmegaplus}, we prove that the
contribution to the counting that will occur in Step \hyperlink{Step6}{6},
Equation \eqref{eq:defcountfunct3}, of the two exceptional finite
subsets $F$ in $\Ga_{\H}\bs\Ga/ \Ga_{\wt D^-}$ and $F'$ in $\Ga_\H$
defined in Step \hyperlink{Step3}{3} is negligible.

\begin{center}
  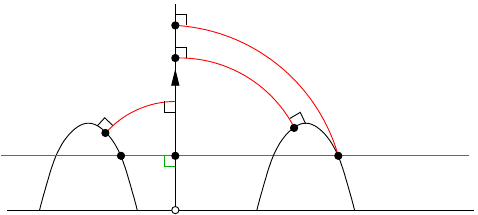
\end{center}
\vspace*{-0.8cm}\addtocounter{fig}{1}
\begin{center}
  {\small Figure \arabic{fig}~: \hypertarget{figstep5}{Exceptional
      common perpendiculars between $D^-$ and $D^+$}.}
\end{center}


Let $[\ga]\in F$ and $\beta\in\Ga_\H$ be such that $d(\beta\ga\wt
D^-,\wt D^+)>0$.  By the invariance under $\Ga_\H$ of the distance
function to $\partial \H$, by the convexity of $\beta\ga\wt D^-$,
since two geodesic rays converging to the same point at infinity $\wt
\ell(+\infty$) come (exponentially) close one of the other, and by the
finiteness of $F$, there exists a constant \addtocounter{const}{1}$
c_{\arabic{const}}'>0$ such that ${\displaystyle \max_{x\in\ga\wt
    D^-\cap \H}}\; d(x,\partial \H)\leq c_{\arabic{const}}'$.  By the
strict convexity of horoballs, since $\ga\wt D^-\cap \H$ is contained
in the $c_{\arabic{const}}'$-neighborhood of $\partial\H$, and by the
convexity of $\ga\wt D^-\cap \H$, there exists a constant
$c_{\arabic{const}}>0$ such that the diameter of $\ga\wt D^-\cap \H$
is at most $c_{\arabic{const}}$. Let $q^+_{\beta\ga}$ (respectively
$q_{\beta\ga}$) be the closest point to $\beta\ga\wt D^-$
(respectively $\beta p_\ga$) on $\wt D^+$, and let $q^-_{\beta\ga}$ be
the closest point to $\wt D^+$ on $\beta\ga\wt D^-$, see Figure
\hyperlink{figstep5}{\arabic{fig}}. Then by the inverse triangle
inequality and since closest point projections do not increase the
distances, we have
\[
d(\beta\ga\wt D^-,\wt D^+)\geq d(\beta p_\ga,\wt D^+) -
d(\beta p_\ga,q^-_{\beta\ga})-d(q^+_{\beta\ga},q_{\beta\ga})\geq
d(\beta p_\ga,\wt D^+) - 2\;c_{\arabic{const}}\,.
\]
Since the Hamenstädt distance and the cuspidal distance are
equivalent, there exists a constant
\addtocounter{const}{1}$c_{\arabic{const}}>0$ (with actually
$c_{\arabic{const}}=1$ when $M$ is real hyperbolic) such that we have
\[
d'_\H(\beta p_\ga,\wt\ell(0))\leq
c_{\arabic{const}}\;d_\H(\beta p_\ga,\wt\ell(0))\,.
\]
Therefore, by Lemma \ref{lem:majodistham} for the second inequality,
we have
\[
e^{d(\beta\ga\wt D^-,\wt D^+)}\geq e^{d(\beta p_\ga,\wt D^+)-
\addtocounter{const}{-1}2\,c_{\arabic{const}}}\geq e^{-2\,c_{\arabic{const}}}\;
d_\H(\beta p_\ga,\wt\ell(0))\geq
\addtocounter{const}{1}\frac{1}{c_{\arabic{const}}}\addtocounter{const}{-1}
\;e^{-2\,c_{\arabic{const}}}\;d'_\H(\beta p_\ga,\wt\ell(0))\,.
\]
Thus by Step \hyperlink{Step4}{4}, we have the following negligible
estimate for $F$:
\begin{align}
&\card\Big\{([\ga],\beta)\in F\times\Ga_\H: \begin{array}{l}
    0<d(\beta\ga \wt D^-,\wt D^+)\leq s\\z_{\beta\ga\wt D^-,\,\wt D^+}
    \in \wt\ell(\,\interval[open]{0}{+\infty}\,)\end{array}\Big\}
  \nonumber\\& \leq (\card \;F)\;\max_{[\ga]\in F}
  \card\big\{\beta\in \Ga_\H: d'_\H(\beta p_\ga,\wt\ell(0))\leq
  \addtocounter{const}{1}c_{\arabic{const}}\addtocounter{const}{-1}\,
  e^{s+2\,c_{\arabic{const}}}\big\}\nonumber \\&=(\card \;F)\;\max_{[\ga]\in F}
  \Phi_{[\ga]}(\addtocounter{const}{1}c_{\arabic{const}}\addtocounter{const}{-1}
  \,e^{s+2\,c_{\arabic{const}}})=\bigO(e^{\delta s})\,. \label{eq:countfunctF}
\end{align}

Let us prove a similar estimate for $F'$. Let $([\ga],\beta) \in
(\Ga_{\H}\bs\Ga/\Ga_{\wt D^-}\ssm F)\times F'$ be such that the
closest point $z_{\beta\ga\wt D^-,\,\wt D^+}$ on $\wt D^+$ to
$\beta\ga\wt D^-$ belongs to $\wt\ell(\,\interval[open]{0}
{+\infty}\,)$. Since $[\ga]\notin F$, the closed convex subsets
$\ga\wt D^-$ and $\H$, hence $\beta\ga\wt D^-$ and $\H$, are disjoint.
By the intermediate value theorem, the common perpendicular between
$\beta\ga\wt D^-$ and $\wt D^+$ meets $\partial \H$ (see Figure
\addtocounter{fig}{-1}\hyperlink{figstep3}{\arabic{fig}}\addtocounter{fig}{1}
before Step \hyperlink{Step4}{4}). Hence we have
\[
d(\ga \wt D^-,\H)=d(\beta\ga \wt
D^-,\H)\leq d(\beta\ga \wt D^-,\wt D^+)\,.
\]
Since the Riemannian orbifold $M$ is locally symmetric with finite
volume, the group $\Ga$ contains finitely many conjugacy classes of
finite subgroups. Hence there exists a constant
\addtocounter{const}{2}$c_{\arabic{const}}>0$ such that for every
$[\ga]\in \Ga_{\H}\bs\Ga/\Ga_{\wt D^-}\ssm F$, we have $m_{\ga\wt D^-,\,\H}
\geq c_{\arabic{const}}$.

Since both $D^-$ and $\V_+$ have finite skinning measures, by 
\cite[Theo.~1]{ParPau17ETDS}, we have $\N_{D^-,\;\V_+}(s)
=\bigO(e^{\delta s})$. Thus
\begin{align}
  &\card\Big\{([\ga],\beta)\in
(\Ga_{\H}\bs\Ga/\Ga_{\wt D^-}\ssm F)\times F': \begin{array}{l}
    0<d(\beta\ga \wt D^-,\wt D^+)\leq s\\z_{\beta\ga\wt D^-,\,\wt D^+}
    \in \wt\ell(\,\interval[open]{0}{+\infty}\,)\end{array}\Big\}
  \nonumber\\& \leq (\card \;F')\; \card\{[\ga]\in
  \Ga_{\H}\bs\Ga/\Ga_{\wt D^-}: 0< d(\ga \wt D^-,\H)\leq
  s\}\nonumber
  \\&\leq \frac{\card \;F'}{c_{\arabic{const}}}\N_{D^-,\;\V_+}(s)
  =\bigO(e^{\delta s})\,. \label{eq:countfunctFprime}
\end{align}

\medskip\noindent{\bf Step \hypertarget{Step6}{6}. } The sixth step in
the proof of Equation \eqref{eq:NOmegaplus} is to disintegrate the
counting defining $\N_{\wt D^-,\Omega_+}(s)$ in Equation
\eqref{eq:defcountfunct2} along the orbits of the parabolic subgroup
$\Ga_\H$ of $\Ga$ fixing $\wt \ell(+\infty)$.

Let $\Ga'_{\wt D^+}=\Ga_{\wt D^+} \cap\Ga_\H$.  Since $\ell$ is not
reversible in Case \hyperlink{Case1}{1}, we have $\Ga'_{\wt
  D^+} =\Ga_{\wt D^+}$, and in particular Equation
\eqref{eq:defcountfunct2} can be rewritten
\begin{equation}\label{eq:defcountfunct2.5}
  \N_{D^-,\,\Omega_+}(s)=\card\Big\{\ga\in\Ga'_{\wt D^+}\bs\Ga/\Ga_{\wt D^-}:
   \begin{array}{l}0<d(\ga \wt D^-,\wt D^+)\leq s\\z_{\ga\wt D^-,\,\wt D^+}
    \in \wt\ell(\,\interval[open]{0}{+\infty}\,)\end{array}\Big\}+
   \bigO(e^{\delta\,s})\;.
   \end{equation}
But what follows will also be useful for Case \hyperlink{Case2}{2},
hence the generality. Since the point at infinity $\wt\ell(+\infty)$ is a
parabolic fixed point of $\Ga$, the group $\Ga'_{\wt D^+}$ is the
pointwise stabilizer of $\wt D^+$, hence has order $m^+=m(D^+)$.

We use the representatives of double classes in $\Ga_{\H} \bs \Ga/
\Ga_{\wt D^-}$ defined in Step \hyperlink{Step3}{3}, though any choice
would work as well. 
The map from
$\Ga'_{\wt D^+} \bs\Ga/\Ga_{\wt D^-}$ to $\Ga_{\H}\bs\Ga/ \Ga_{\wt
  D^-}$ given by $\Ga'_{\wt D^+}\ga'\Ga_{\wt D^-} \mapsto \Ga_{\H}
\ga'\Ga_{\wt D^-}$ is well defined since $\Ga'_{\wt D^+}$ is contained
in $\Ga_{\H}$. Its fiber over the element $[\ga]\in \Ga_{\H} \bs\Ga/
\Ga_{\wt D^-}$ is the subset $\{ \Ga'_{\wt D^+} \beta \ga \Ga_{\wt
  D^-}:\beta \in \Ga_\H\}$ of $\Ga'_{\wt D^+}\bs\Ga/\Ga_{\wt D^-}$.
Let us fix $[\ga] \in\Ga_{\H} \bs\Ga/\Ga_{\wt D^-}\ssm F$, so that
$\H$ and $\ga\wt D^-$ have a common perpendicular. Given two distinct
elements $\Ga'_{\wt D^+} \beta,\Ga'_{\wt D^+} \beta'$ in $\Ga'_{\wt
  D^+}\bs \Ga_{\H}$, we have $\Ga'_{\wt D^+}\beta \ga\Ga_{\wt D^-} =
\Ga'_{\wt D^+} \beta' \ga\Ga_{\wt D^-}$ if and only if there exists
$\alpha\in\Ga_{\wt D^-}$ such that $\beta' \ga\alpha\ga^{-1}
\beta^{-1} \in\Ga'_{\wt D^+}$, hence if and only if $\ga\Ga_{\wt D^-}
\ga^{-1}\cap (\beta')^{-1} \Ga'_{\wt D^+}\beta$ is nonempty. As 
the classes $\Ga'_{\wt D^+} \beta$ and $\Ga'_{\wt D^+}\beta'$ are
distinct, this implies  $\Ga_{\ga\wt D^-}\cap \Ga_\H\neq
\{\id\}$. Since $[\ga]\notin F$, the multiplicity $m_{\ga\wt D^-,\,\H}$
defined in Equation \eqref{eq:defmultipli} of the common perpendicular
between $\ga\wt D^-$ and $\H$ is different from $1$. By Equation
\eqref{eq:decaymultnon1} applied with $(\H,\wt D^-)$ instead of $(\wt
D^+,\wt D^-)$, outside a number of elements $[\ga]\in\Ga_{\H}
\bs\Ga/\Ga_{\wt D^-}\ssm F$ that is a $\bigO(e^{\delta s})$, the map
$\Ga'_{\wt D^+}\beta \mapsto \Ga'_{\wt D^+}\beta \ga\Ga_{\wt D^-}$ is
injective. Note that the canonical map $\Ga_\H\ra \Ga'_{\wt D^+}
\bs\Ga_\H$ is $m^+$-to-1.

By Equation \eqref{eq:countfunctF} that controls the contribution of
the double classes $[\ga]\in \Ga_{\H} \bs\Ga/\Ga_{\wt D^-}$ that are
in $F$, and since the two conditions below on $([\ga],\beta)$ are
invariant under multiplying $\beta$ on the left by any element of
$\Ga'_{\wt D^+}$, Equation \eqref{eq:defcountfunct2.5} hence becomes
\begin{align}
  \N_{D^-,\,\Omega_+}(s)&=\frac{1}{m^+}\card\Big\{([\ga],\beta)\in
  (\Ga_{\H}\bs\Ga/\Ga_{\wt D^-})\times\Ga_\H: \!\!\begin{array}{l}
    0<d(\beta\ga \wt D^-,\wt D^+)\leq s\\z_{\beta\ga\wt D^-,\,\wt D^+}
    \in \wt\ell(\,\interval[open]{0}{+\infty}\,)\end{array}\!\Big\}
  \nonumber\\ &\;\;+\bigO(e^{\delta\,s}) \,.
  \label{eq:defcountfunct3}
\end{align}

\medskip\noindent{\bf Step \hypertarget{Step7}{7}. } In this final
 step in the proof of Equation \eqref{eq:NOmegaplus}, we
compute the contribution to the counting in Equation
\eqref{eq:defcountfunct3} of the elements in the main domain
$(\Ga_{\H}\bs\Ga/ \Ga_{\wt D^-}\ssm F)\times(\Ga_\H\ssm F')$, and we
conclude the proof of Equation \eqref{eq:NOmegaplus}. Most of the
technical work is devoted to getting an error term.

For every $s> 1$, let
\[
\Sigma_s=\card\Big\{([\ga],\beta)\in
(\Ga_{\H}\bs\Ga/\Ga_{\wt D^-}\ssm F)\times (\Ga_\H\ssm F'): \begin{array}{l}
    s-1<d(\beta\ga \wt D^-,\wt D^+)\leq s\\z_{\beta\ga\wt D^-,\,\wt D^+}
    \in \wt\ell(\,\interval[open]{0}{+\infty}\,)\end{array}\Big\}\,.
\]
The second assumption above is superfluous, since by the definition of
the set $F'$ in Step \hyperlink{Step3}{3}, we have $z_{\beta\ga\wt
  D^-,\,\wt D^+} \in \wt\ell(\,\interval[open]{0} {+\infty}\,)$
whenever  $[\ga]\in \Ga_{\H}\bs
\Ga/\Ga_{\wt D^-}\ssm F$ and $\beta\in \Ga_\H\ssm F'$.

Let $\eta\in\interval[open]{0}{1}\,$ (that will tend to $0$ at the end
of the proof). Recall that $d(\ga\wt D^-,\H)>0$ as $[\ga]\in
\Ga_{\H}\bs\Ga/\Ga_{\wt D^-} \ssm F$. By summing over thin slices with
width $\eta$ of the first factor elements, we have
\begin{align}\label{eq:doublesomm1}
\Sigma_s=\sum_{k=1}^{+\infty}\sum_{\substack{[\ga]\in
    \Ga_{\H}\bs\Ga/\Ga_{\wt D^-}\ssm\, F\\
    (k-1)\eta<d(\ga\wt D^-,\,\H)\leq k\eta}}
\quad\sum_{\substack{\beta\in\Ga_\H\ssm\, F'\\
    s-1<d(\beta\ga \wt D^-,\wt D^+)\leq s}} 1\,.
\end{align}
Let $k\in\NN\ssm\{0\}$ and $([\ga],\beta)\in
(\Ga_{\H}\bs\Ga/\Ga_{\wt D^-}\ssm F)\times (\Ga_\H\ssm F')$ be such
that
\begin{align}\label{eq:condigabeta}
(k-1)\eta<d(\ga\wt D^-,\H)\leq k\eta\quad\text{and}\quad
s-1 < d(\beta\ga \wt D^-,\wt D^+)\leq s\,.
\end{align}

Since $M$ is real hyperbolic (so that $d_\H=d'_\H$), by the first
claim of Lemma \ref{lem:distancetodivline}
\hyperlink{lemhypreelii}{(ii)} applied with $D=\beta\ga\wt D^-$ and
$x_0=\beta p_\ga$, whose assumptions are satisfied by Equation
\eqref{eq:minocKK}, and since $\H$ is invariant under $\beta^{-1}\in
\Ga_\H$, we have
\begin{align}
d(\beta\ga\wt D^-,\wt D^+) &=
d(\ga\wt D^-,\H)+\ln (2\,d'_\H(\beta p_\ga,\wt\ell(0)))\nonumber\\ &\quad +
\bigO\big(d'_\H(\beta p_\ga,\wt\ell(0))^{-2}\,e^{-2\,d(\ga\wt D^-,\,\H)}\big)\,.
\label{eq:prepabootstrap}
\end{align}
In particular, up to increasing the finite set $F'$, we may assume
that $d'_\H(\beta p_\ga,\wt\ell(0)))$ is large enough so that
$d(\ga\wt D^-,\H)\leq d(\beta\ga\wt D^-,\wt \ell)$. Let $N=\lfloor
\frac{s}{\eta} \rfloor$, so that we have $N\eta \leq s< (N+1)\eta$ and
in the summation \eqref{eq:doublesomm1}, we may restrict $k$ to vary
between $1$ and $N+1$ for a majoration and between $1$ and $N$ for a
minoration.

Since $d'_\H(\wt\ell(0),\beta p_\ga)\geq c_\KK$ by Equation
\eqref{eq:minocKK}, and since $e^{-2\,d(\ga\wt D^-,\,\H)}\leq 1$,
Equations \eqref{eq:prepabootstrap} and \eqref{eq:condigabeta} give
\[
d'_\H(\beta p_\ga,\wt\ell(0))=\frac{1}{2}\;
e^{d(\beta\ga\wt D^-,\,\wt D^+)-d(\ga\wt D^-,\,\H)+\bigO(1)}
= e^{s -k\eta+\bigO(1)}\;.
\]
Thus $d'_\H(\beta p_\ga,\wt\ell(0))^{-2}\,e^{-2\,d(\ga\wt D^-,\,\H)}=
\bigO(e^{-2\,s +2k\eta})\bigO(e^{-2k\eta})=\bigO(e^{-2\,s})$.
Bootstrapping this in Equation \eqref{eq:prepabootstrap}, and by
Equation \eqref{eq:condigabeta}, we have
\begin{equation}\label{eq:doubleencadmino}
\frac{1}{2}\;e^{s-1 -k\eta +\bigO(e^{-2\,s})}\leq
d'_\H(\beta p_\ga,\wt\ell(0))  \leq 
\frac{1}{2}\;e^{s -(k-1)\eta +\bigO(e^{-2\,s})}\;.
\end{equation}
Conversely (this will be used only at the end of Step
\hyperlink{Step7}{7}), if we had
\begin{equation}\label{eq:doubleencadmajo}
\frac{1}{2}\;e^{s-1 -(k-1)\eta +\bigO(e^{-2\,s})}\leq
d'_\H(\beta p_\ga,\wt\ell(0))  \leq 
\frac{1}{2}\;e^{s -k\eta +\bigO(e^{-2\,s})}\;,
\end{equation}
for an appropriate function $\bigO(\cdot)$ that is independent of
$\eta$, while still having the inequalities $(k-1)\eta<d(\ga\wt
D^-,\H)\leq k\eta$, then by Equation \eqref{eq:prepabootstrap}, we
would have the inequalities
\[
s-1 +\bigO(e^{-2\,s})< d(\beta\ga \wt D^-,\wt D^+)\leq s+\bigO(e^{-2\,s})\,.
\]
Note that the right hand sides of the two Equations
\eqref{eq:doubleencadmino} and \eqref{eq:doubleencadmajo} differ by a
multiplicative factor $e^{\bigO(\eta)}=1+\bigO(\eta)$ as $\eta$ tends
to $0$. The same holds for their left hand sides.

It follows from Equation \eqref{eq:doubleencadmino}, by Step
\hyperlink{Step4}{4} and Equation \eqref{eq:horocount}, that we have
(with a function $\bigO(\cdot)$ that is independent of $\eta$ and
$[\ga]$)
\begin{align}
&\card\{\beta\in\Ga_\H\ssm F':
s-1<d(\beta\ga \wt D^-,\wt D^+)\leq s\}\nonumber\\\leq\;&
  \Phi_{[\ga]}\big(\frac{1}{2}\;e^{s -(k-1)\eta +\bigO(e^{-2\,s})}\big)-
  \Phi_{[\ga]}\big(\frac{1}{2}\;e^{s-1 -k\eta +\bigO(e^{-2\,s})}\big)
  \nonumber\\=\;&\frac{\Xi_{\wt M}}{2^\delta\,\|\sigma^-_{\V_+}\|}
  e^{\delta s -\delta k\eta}\big(e^{\delta\eta
  +\bigO(e^{-s})}-e^{-\delta+\bigO(e^{-s})}\big)
  +\bigO(e^{(s -k\eta)(\delta -1)})\,.\label{eq:cotebeta}
\end{align}
Let $C_1=\frac{\Xi_{\wt M}}{2^\delta\,\|\sigma^-_{\V_+}\|}\,
(e^{\delta\eta+\bigO(e^{-s})}-e^{-\delta+\bigO(e^{-s})})$, let
$C'_1>0$ be a large enough constant (to be precised later on), and let
$f_\pm:[0,+\infty[\;\ra\RR$ be the function
\[
f_\pm:t\mapsto C_1\,e^{\delta s-\delta t\eta}\pm C'_1\,e^{(s-t\eta)(\delta-1)}\,.
\]
Since $s=N\eta+\bigO(\eta)$, we have $f_\pm(N+1)=\bigO(1)$.  For every
$k\in\NN$, let
\[
a_k=\card\{[\ga]\in \Ga_{\H}\bs\Ga/\Ga_{\wt D^-}\ssm F:
(k-1)\eta<d(\ga\wt D^-,\H)\leq k\eta\}\,.
\]
Note that $a_0=0$ by the definition of $F$. By \cite[Theo.~15
  (2)]{ParPau17ETDS}, which can be applied since its assumption on the
exponential decay of correlations is satisfied as recalled in Section
\ref{sec:geomback} since $M$ is locally symmetric with finite volume,
and by a further regularisation process in order to remove the
smoothness assumption on $\partial D^-$, there exists $\kappa''>0$ and
a function $\bigO(\cdot)$ that is independent of $\eta$ such that for
every $t\geq 0$, we have
\begin{align}
\sum_{k=0}^t\;a_k=
\frac{\|\sigma^+_{D^-}\|\;\|\sigma^-_{\V_+}\|}{\delta\;\|m_{\rm BM}\|}\;
e^{\delta t\eta}+\bigO(e^{(\delta -\kappa'')t\eta})\,.\label{eq:cotega}
\end{align}

Let $C_2=\frac{\|\sigma^+_{D^-}\|\,\|\sigma^-_{\V_+}\|}{\delta\,
  \|m_{\rm BM}\|}$.  By Equations \eqref{eq:doublesomm1} and
\eqref{eq:cotebeta}, for a large enough constant $C'_1$ in the
definition of the function $f_+$, by Abel's summation formula, by
Equation \eqref{eq:cotega} and again since $N\eta=s+\bigO(\eta)$, we
have
\begin{align*}
  \Sigma_s&\leq \sum_{k=0}^{N+1} a_k\,f_+(k)= \big(\sum_{k=0}^{N+1}
  a_k\,\big)f_+(N+1)- \int_0^{N+1}(\sum_{k=0}^{t}
  a_k\,\big)f'_+(t)\,dt\\&=\bigO(e^{\delta\,s})+
  \int_0^{N+1}\big(C_2\,e^{\delta t\eta}+\bigO(e^{(\delta -\kappa'')t\eta})\big)
  \big(\delta\,\eta\,C_1\, e^{\delta s-\delta t \eta}+
  \bigO(\eta\, e^{(s-t\eta)(\delta-1)})\big)\,dt\\&=\bigO(e^{\delta\,s})+
  \delta\,C_2\,C_1\,(N+1)\,\eta\,e^{\delta s}\\&\qquad\qquad\;+
  \bigO\big(e^{\delta\,s}\int_0^{N+1}\eta\,e^{-\kappa''t\eta}\,dt\big)
  +\bigO\big(e^{s(\delta-1)}\int_0^{N+1}\eta\,e^{t\eta}\,dt\big)
  \\&=\delta\,C_2\,C_1\,(s+\bigO(\eta))\,e^{\delta s}+\bigO(e^{\delta\,s})\,.
\end{align*}
Replacing $C_1$ and $C_2$ by their values, and letting $\eta$ tend to
$0$, we hence have
\[
  \Sigma_s\leq
\frac{\|\sigma^+_{D^-}\|\;\Xi_{\wt M}}{2^\delta\;\|m_{\rm BM}\|}
\;s\;e^{\delta \,s}(1-e^{-\delta})+\bigO(e^{\delta\, s})\,.
\]

The same lower bound is obtained similarly, replacing Equation
\eqref{eq:doubleencadmino} by Equation \eqref{eq:doubleencadmajo}
(which only replaces $s$ by $s+\bigO(e^{-2\,s})$ on the first line of
Equation \eqref{eq:cotebeta}, and this does not change its second
line, nor therefore its third line) and the function $f_+$ by the
function $f_-$.  By a summation, we have
\begin{align}
  &\card\Big\{([\ga],\beta)\in
(\Ga_{\H}\bs\Ga/\Ga_{\wt D^-}\ssm F)\times (\Ga_\H\ssm F'): \begin{array}{l}
    0<d(\beta\ga \wt D^-,\wt D^+)\leq s\\z_{\beta\ga\wt D^-,\,\wt D^+}
    \in \wt\ell(\,\interval[open]{0}{+\infty}\,)\end{array}\Big\}\;
  \nonumber\\&=\frac{\|\sigma^+_{D^-}\|\;\Xi_{\wt M}}{2^\delta\,\|m_{\rm BM}\|}\;
  s\,e^{\delta\,s} +\bigO(e^{\delta s})\,. \label{eq:countfunctmain}
\end{align}
By separating the counting domain $(\Ga_{\H}\bs\Ga/\Ga_{\wt D^-})
\times \Ga_\H$ as the disjoint union of $F\times\Ga_\H$, of
$(\Ga_{\H}\bs\Ga/\Ga_{\wt D^-}\ssm F)\times F'$ and of
$(\Ga_{\H}\bs\Ga/\Ga_{\wt D^-}\ssm F)\times (\Ga_\H\ssm F')$, Equation
\eqref{eq:NOmegaplus} finally follows from Equations
\eqref{eq:defcountfunct3}, \eqref{eq:countfunctF},
\eqref{eq:countfunctFprime} and \eqref{eq:countfunctmain}.

\medskip\noindent{\bf Case \hypertarget{Case2}{2}. }  Let us now
assume that $\ell$ is reversible. We then have $\Omega_-
=\Omega_+$. Hence $\N_{D^-,\Omega_-}(t)= \N_{D^-,\Omega_+} (t)$ and
\begin{equation}\label{eq:NDDOmegapmrec}
\N_{D^-,\Omega_+}(t)\leq \N_{D^-,D^+}(t)\leq \N_{D^-,\Omega_+}(t) +
\N_{D^-,\Omega_0}(t)
\end{equation}
for every $t\geq 0$. Let us prove that Equation \eqref{eq:NOmegaplus}
is still satisfied. Since $\iota_{\rm rev}^+=1$ when $\ell$ is
reversible, this will prove as in Case \hyperlink{Case1}{1},
replacing the call to Equation \eqref{eq:NDDOmegapm} by a call to
Equation \eqref{eq:NDDOmegapmrec}, that Equation
\eqref{eq:NDDabstrait} is still satisfied. Then by the same
computations as in Case \hyperlink{Case1}{1}, Theorem
\ref{theo:1geoddivreal} when $\ell$ is reversible will follow.

Since the divergent geodesic $\ell$ is reversible, there exists an element $\iota_{\wt
  D^+}\in\Ga$ such that 
 $\iota_{\wt D^+}\,\wt\ell
(\,\interval[open]{0}{+\infty}\,)=\wt\ell (\,\interval[open]{-\infty}
{t_-}\,)$. By the definition of $\Omega_+$, by the commutativity of
the diagram \eqref{eq:ppbullet} and since the family $(H_\xi)_{\xi\in
  \operatorname{Par}_\Ga}$ is precisely invariant, we have
\begin{align*}
  \partial^1_-\wt D^+\cap \wt p^{\,-1}(\Omega_+)&=
  \partial^1_-\wt D^+\cap \wt p^{\,-1}\big({p_\bullet}^{-1}
  (p\circ\wt\ell(\,\interval[open]{0}{+\infty}\,))\big)\\
  &=
  \partial^1_-\wt D^+\cap
  \Ga \,{\wt p_\bullet}^{\;-1}(\wt\ell(\,\interval[open]{0}{+\infty}\,))\\
  &=
  \partial^1_-\wt D^+\cap \big({\wt p_\bullet}^{\;-1}
  (\wt\ell(\,\interval[open]{0}{+\infty}\,))
  \cup\iota_{\wt D^+} {\wt p_\bullet}^{\;-1}(\wt\ell(\,\interval[open]{0}
  {+\infty}\,))\big)\\&=\partial^1_-\wt D^+\cap {\wt p_\bullet}^{\;-1}
    (\wt\ell(\,\interval[open]{-\infty}{t_-}\,\cup\,\interval[open]{0}
  {+\infty}\,))\,.
\end{align*}
Hence, as in Step \hyperlink{Step1}{1}, we now have
\begin{align}
  \N_{D^-,\,\Omega_+}(s)=
  \sum_{[\ga]\in\Ga_{\wt D^+}\bs\Ga/\Ga_{\wt D^-} \,:\;
    0<d(\ga \wt D^-,\wt D^+)\leq s,\; z_{\ga\wt D^-,\,\wt D^+}\,
    \in\, \wt\ell(\,\interval[open]{-\infty}{t_-}\,
    \cup\,\interval[open]{0}{+\infty}\,)} m_{\ga\wt D^-,\wt D^+}\,.
  \label{eq:defcountfunct1rec}
\end{align}
Since $\ell$ is reversible, the intersection $\Ga'_{\wt D^+} =
\Ga_{\wt D^+} \cap\Ga_\H$
now has index $2$ in $\Ga_{\wt D^+}$.  Given a double class
$[\ga]\in\Ga_{\wt D^+}\bs\Ga/\Ga_{\wt D^-}$ such that $0<d(\ga \wt
D^-,\wt D^+)\leq s$, its preimage by the canonical projection
$\Ga'_{\wt D^+}\bs\Ga/\Ga_{\wt D^-}\ra \Ga_{\wt D^+}\bs\Ga/\Ga_{\wt
  D^-}$ is the set $\{\,\Ga'_{\wt D^+} \ga \Ga_{\wt D^-}, \Ga'_{\wt
  D^+}\iota_{\wt D^+}\ga\Ga_{\wt D^-} \,\}$. We have $\Ga'_{\wt D^+}
\ga \Ga_{\wt D^-} = \Ga'_{\wt D^+} \iota_{\wt D^+} \ga \Ga_{\wt D^-}$
if and only if there exists $\alpha\in \Ga_{\wt D^-}$ such that
$\iota_{\wt D^+} \ga \alpha\ga^{-1}\in \Ga'_{\wt D^+}$, hence if and
only if $\Ga_{\ga \wt D^-} \cap \iota_{\wt D^+} ^{-1} \Ga'_{\wt D^+}$
is nonempty. Since $\iota_{\wt D^+}^{-1}\Ga'_{\wt D^+}$ is not the
trivial class in $\Ga_{\wt D^+} /\Ga'_{\wt D^+}$, this implies that
$m_{\ga \wt D^-,\wt D^+}$ is different from $1$. Hence by Equation
\eqref{eq:decaymultnon1}, the canonical map $\Ga'_{\wt
  D^+}\bs\Ga/\Ga_{\wt D^-}\ra \Ga_{\wt D^+} \bs \Ga/\Ga_{\wt D^-}$ is
$2$-to-$1$ outside a number $\bigO(e^{\delta s})$ of elements, and
exactly one $\Ga'_{\wt D^+} \ga \Ga_{\wt D^-}$ of the two preimages
satisfies that $z_{\ga \wt D^-,\wt D^+}$ belongs to
$\wt\ell(\,\interval[open] {0}{+\infty}\,)$. Thus, as in Steps
\hyperlink{Step1}{1} and \hyperlink{Step2}{2}, we have
\begin{equation}
  \N_{D^-,\,\Omega_+}(s)=\card\Big\{\ga\in\Ga'_{\wt D^+}\bs\Ga/\Ga_{\wt D^-}:
   \begin{array}{l}0<d(\ga \wt D^-,\wt D^+)\leq s\\z_{\ga\wt D^-,\,\wt D^+}
    \in \wt\ell(\,\interval[open]{0}{+\infty}\,)\end{array}\Big\}+
    \bigO(e^{\delta\,s}) \,, 
  \label{eq:defcountfunct2rec}
\end{equation}
that is, Equation \eqref{eq:defcountfunct2.5} is still valid. As in
Step \hyperlink{Step6}{6}, we therefore have
\begin{align}
  \N_{D^-,\,\Omega_+}(s)&=\frac{1}{m^+}\card\Big\{([\ga],\beta)\in
  (\Ga_{\H}\bs\Ga/\Ga_{\wt D^-})\times\Ga_\H: \!\!\begin{array}{l}
    0<d(\beta\ga \wt D^-,\wt D^+)\leq s\\z_{\beta\ga\wt D^-,\,\wt D^+}
    \in \wt\ell(\,\interval[open]{0}{+\infty}\,)\end{array}\!\Big\}
  \nonumber\\ &\;\;+\bigO(e^{\delta\,s}) \,,
\label{eq:defcountfunct3rec}
\end{align}
that is, Equation \eqref{eq:defcountfunct3} is still valid. The
remainder of the proof of Equation \eqref{eq:NOmegaplus}, that is, its
Step \hyperlink{Step7}{7}, now proceeds exactly as in Case
\hyperlink{Case1}{1}.
\cqfd

\btheo\label{theo:divergentperp2} Let $M$ be a noncompact finite
volume complete connected real hyperbolic good orbifold of dimension
$n$. Let $D^+$ and $D^-$ be the images of two divergent geodesics in
$M$. Then, as $s\to+\infty$, we have
\[
\N_{D^-,D^+}(s)=\frac{(n-1)\;\pi^{\frac{n}{2}-1}\;\Ga(\frac{n}{2})
  \,\iota_{\rm rev}(D^-)\,\iota_{\rm rev}(D^+)}
  {2^{n+1}\;\Ga(\frac{n+1}{2})^2\;m(D^-)\;m(D^+)\;\Vol M}\;
  s^2\, e^{(n-1) \,s}+ \bigO\big(s\;e^{(n-1) \,s}\,\big)\,.
\]
\etheo

\dem
The strategy is similar to the one we used in the proof of Theorem
\ref{theo:1geoddivreal}, except that we will now disintegrate the
study of $\N_{D^-,D^+}$ over the study of the number $\N_{\V_\pm,D^+}$
of common perpendiculars starting from a Margulis neighbourhood
$\V_\pm$ of an end of $D^-$ and arriving at $D^+$, and replace the
call to \cite{ParPau17ETDS} in Equations \eqref{eq:NOmegazero},
\eqref{eq:countfunctFprime} and \eqref{eq:cotega} by a call to Theorem
\ref{theo:1geoddivreal} that we just proved.

The notation is now the following one (and differs from the one at the
beginning of the proof of Theorem \ref{theo:1geoddivreal}). Let $\ell$
be a divergent geodesic in $M$ whose image is $D^-$.  Let $\V_-$ and
$\V_+$ be the initial and terminal Margulis neighbourhoods of
$\ell$. Let $t_-$ be the first exit time of $\ell$ from $\V_-$ and let
$t_+$ be the last entry time of $\ell$ into $\V_+$. We may assume that
$t_+=0$. Let $\wt\ell$ be a lift of $\ell$ in $\wt M$, and let $\wt
D^-$ be the image of $\wt\ell$. For simplicity, let $m^\pm=m(D^\pm)$
and $\iota_{\rm rev}^\pm =\iota_{\rm rev}(D^\pm)$.  Let
\begin{align*}
  \Omega_- & =\big\{v\in\partial^1_+D^-:
  p_\bullet(v)\in\ell(\,\interval[open]{-\infty}{t_-}\,)\big\}\,,\\
  \Omega_0\, & =\big\{v\in\partial^1_+D^-:
  p_\bullet(v)\in\ell(\interval{t_-}{0})\big\}\quad\textrm{and}\\
  \Omega_+ & =\big\{v\in\partial^1_+D^-:
  p_\bullet(v)\in\ell(\,\interval[open]{0}{+\infty}\,)\big\}\;.
\end{align*}
We denote by $\H=H_{\wt\ell(+\infty)}$ the horoball of the 
family $(H_\xi)_{\xi\in\operatorname{Par}_\Ga}$ with point at infinity
$\wt\ell(+\infty)$, that is a lift of $\V_+$. Let $\wt D^+$ be a
geodesic line in $\wt M$ whose image in $M$ is $D^+$.

\medskip\noindent{\bf Case \hypertarget{Case1bis}{1}. } Let us first
assume that $\ell$ is not reversible.
As for Equation \eqref{eq:NDDOmegapm}, the subsets $\Omega_-$ and
$\Omega_+$ are disjoint and for every $s\geq 0$, we have
\begin{equation}\label{eq:NDDOmegapmthm4}
  \N_{\Omega_-,D^+}(s)+\N_{\Omega_+,D^+}(s)\leq
\N_{D^-,D^+}(s)\leq\N_{\Omega_-,D^+}(s)+\N_{\Omega_0,D^+}(s)+\N_{\Omega_+,D^+}(s)\;.
\end{equation}
We shall prove that as $s\ra+\infty$, we have
\begin{equation}\label{eq:NOmegaplusthm4}
  \N_{\Omega_+,D^+}(s)=
  \frac{\delta\;\iota_{\rm rev}^+\;{\Xi_{\wt M}}^2}
       {2^{2\delta+1}\;m^-\,m^+\,\|m_{\rm BM}\|}\;
  s^2\,e^{\delta\,s}+\bigO(s\;e^{\delta\,s})\,.
\end{equation}
By the same argument as in the proof of Equation
\eqref{eq:NOmegamoins}, we will also have
\begin{equation}\label{eq:NOmegamoinsthm4}
  \N_{\Omega_-,D^+}(s)=
  \frac{\delta\;\iota_{\rm rev}^+\;{\Xi_{\wt M}}^2}
       {2^{2\delta+1}\;m^-\,m^+\,\|m_{\rm BM}\|}\;
  s^2\,e^{\delta\,s}+\bigO(s\;e^{\delta\,s})\,.
\end{equation}
Let $D_0=\ell([t_-,t_+])$, which is a compact nonempty properly
immersed locally convex subset of $M$, hence has a nonzero finite
outer skinning measure. By Theorem \ref{theo:1geoddivreal}, as
$s\ra+\infty$, we therefore have
\begin{equation}\label{eq:NOmegazerothm4}
  \N_{\Omega_0,D^+}(s)\leq \N_{D_0,D^+}(s)= \bigO(s\;e^{\delta\,s})\,.
\end{equation}
Thus by Equations \eqref{eq:NDDOmegapmthm4},
\eqref{eq:NOmegaplusthm4}, \eqref{eq:NOmegamoinsthm4} and
\eqref{eq:NOmegazerothm4}, since $\iota_{\rm rev}^-=2$ as $\ell$ is
not reversible, as $s\ra+\infty$, we will have
\begin{equation}\label{eq:NDDabstraitthm4}
  \N_{D^-,D^+}(s)=
  \frac{\delta\;\iota_{\rm rev}^-\;\iota_{\rm rev}^+\;{\Xi_{\wt M}}^2}
       {2^{2\delta+1}\;m^-\,m^+\,\|m_{\rm BM}\|}\;
  s^2\,e^{\delta\,s}+\bigO(s\;e^{\delta\,s})\,.
\end{equation}

As $M$ is finite volume real hyperbolic, we have $\delta=n-1$, and
Theorem \ref{theo:divergentperp2} will follow from Equation
\eqref{eq:NDDabstraitthm4} using Equations \eqref{eq:Xirealhyp} and
\eqref{eq:PP17ETDS20_1}, once we have proven Equation
\eqref{eq:NOmegaplusthm4}.

The remainder of the proof is devoted to proving Equation
\eqref{eq:NOmegaplusthm4}. For every element $\ga\in\Ga$ such that
$d(\wt D^-, \ga\wt D^+)>0$, we now denote by $z_{\wt D^-,\ga\wt D^+}
\in\wt D^+$ the origin of the common perpendicular from $\wt D^-$
to $\ga\wt D^+$. As in Steps \hyperlink{Step1}{1} and
\hyperlink{Step2}{2} in the proof of Theorem \ref{theo:1geoddivreal},
since $\bigO(s^2e^{(\delta-1)s})=\bigO(e^{\delta \,s})$, for every
$s\geq 0$, we have
\begin{align}
  \N_{\Omega_+,D^+}(s)&=
  \sum_{[\ga]\in\Ga_{\wt D^-}\bs\Ga/\Ga_{\wt D^+}\,:\;
    0<d(\wt D^-,\ga\wt D^+)\leq s\,,\; z_{\wt D^-,\ga\wt D^+}\,
    \in\, \wt\ell(\,\interval[open]{0}{+\infty}\,)} m_{\wt D^-,\ga\wt D^+}
  \nonumber\\ &=
  \card\bigg\{\ga\in\Ga_{\wt D^-}\bs\Ga/\Ga_{\wt D^+}:
   \begin{array}{l}0<d(\wt D^-,\ga \wt D^+)\leq s\\z_{\wt D^-,\,\ga\wt D^+}
    \in \wt\ell(\,\interval[open]{0}{+\infty}\,)\end{array}\bigg\}+
    \bigO(e^{\delta\,s}) \,.
  \label{eq:defcountfunct2thm4}
\end{align}
As in the first part of Step \hyperlink{Step3}{3} in the proof of
Theorem \ref{theo:1geoddivreal}, we now define
\[
F=\{[\ga]\in\Ga_{\H}\bs\Ga/\Ga_{\wt D^+}:\H\cap \ga \wt D^+\neq\emptyset\}
\]
and for every $\ga\in\Ga$ such that $\Ga_{\H}\ga\Ga_{\wt D^+}\notin
F$, we now denote by $p_\ga\in\partial \H$ and $p'_\ga\in\partial
(\ga\wt D^+)$ the two endpoints of the common perpendicular between
$\H$ and $\ga \wt D^+$, see Figure
\addtocounter{fig}{1}\hyperlink{figstep3pm}{\arabic{fig}}\addtocounter{fig}{-1}.
\begin{center}
  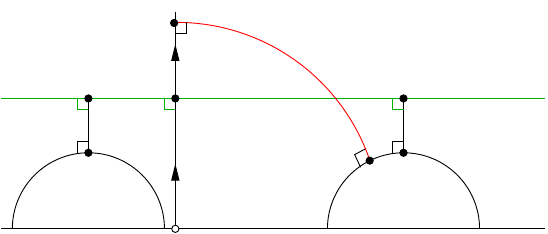
\end{center}
\vspace*{-0.6cm}\addtocounter{fig}{1}
\begin{center}
{\small Figure \arabic{fig}~: \hypertarget{figstep3pm}{Generic
      common perpendiculars between $D^+$ and $D^-$}.}
\end{center}
As in the second part of Step \hyperlink{Step3}{3} in the proof of
Theorem \ref{theo:1geoddivreal} (recalling that $c_\KK=2$ in the real
hyperbolic case), there exists
\addtocounter{const}{1}$c_{\arabic{const}}>0$ and a choice of
representatives $\ga$ in $[\ga] \in\Ga_{\H} \bs\Ga/\Ga_{\wt D^+}$ such
that $d'_\H(p_\ga,\wt\ell(0)) \leq c_{\arabic{const}}$. We now define
\[
F'=\big\{\beta\in\Ga_\H:d'_\H(\wt \ell(0),\beta\,\wt
\ell(0)) < c_\KK+c_{\arabic{const}}\big\}\,.
\]
Since $M$ is real hyperbolic and by the second claim of Lemma
\ref{lem:distancetodivline} \hyperlink{lemhypreelii}{(ii)} applied
with $D=\beta\ga\wt D^+$ and $x_0=\beta p_\ga$, if $[\ga]\in\Ga_{\H}
\bs\Ga/\Ga_{\wt D^+}\ssm F$ and $\beta\in\Ga_\H\ssm F'$, then $z_{\wt
  D^-,\,\beta\ga\wt D^+}\in \wt\ell(\,\interval[open]{0}{+\infty}\,)$
(see Figure \hyperlink{figstep3pm}{\arabic{fig}}).  For all
$[\ga]\in\Ga_{\H}\bs\Ga/\Ga_{\wt D^+}\ssm F$ and $t\geq 0$, let us now
define
\[
\Phi_{[\ga]}(t)=\card\{\beta\in\Ga_\H:d'_\H(\wt\ell(0),\beta p_\ga)\leq t\}\,.
\]
As in Step \hyperlink{Step4}{4} in the proof of Theorem
\ref{theo:1geoddivreal}, we have
\begin{equation}\label{eq:horocountthm4}
\Phi_{[\ga]}(t)=\frac{\Xi_{\wt M}}{\|\sigma^+_{\V_+}\|}\;t^\delta
+\bigO(t^{\delta-1})\,.
\end{equation}
As in the first part of Step \hyperlink{Step5}{5} in the proof of
Theorem \ref{theo:1geoddivreal}, we have
\begin{align}
&\card\bigg\{([\ga],\beta)\in F\times\Ga_\H: \begin{array}{l}
    0<d(\wt D^-,\beta\ga \wt D^+)\leq s\\z_{\wt D^-,\,\beta\ga \wt D^+}
    \in \wt\ell(\,\interval[open]{0}{+\infty}\,)\end{array}\bigg\}
  =\bigO(e^{\delta s})\,. \label{eq:countfunctFthm4}
\end{align}
Since $\V_+$ has finite outer skinning measure and since $D^+$ is a
divergent geodesic, by Theorem \ref{theo:1geoddivreal}, we have
\begin{equation}\label{eq:NVD}
\N_{\V_+,D^+}(s)=\bigO(se^{\delta\,s})\,.
\end{equation}
Hence as in the second part of Step \hyperlink{Step5}{5} in the proof of
Theorem \ref{theo:1geoddivreal}, we have
\begin{align}
  &\card\bigg\{([\ga],\beta)\in
(\Ga_{\H}\bs\Ga/\Ga_{\wt D^+}\ssm F)\times F': \begin{array}{l}
    0<d(\wt D^-,\beta\ga \wt D^+)\leq s\\z_{\wt D^-,\,\beta\ga\wt D^+}
    \in \wt\ell(\,\interval[open]{0}{+\infty}\,)\end{array}\bigg\}
\nonumber\\=\;&\bigO(s\,e^{\delta s})\,. \label{eq:countfunctFprimethm4}
\end{align}
As in Step \hyperlink{Step6}{6} in the proof of Theorem
\ref{theo:1geoddivreal}, since $\ell$ is not reversible, the
stabilizer $\Ga_{\wt D^-}$ of $\wt D^-$ coincides with its pointwise
stabilizer $\Ga'_{\wt D^-}=\Ga_{\wt D^-}\cap\Ga_\H$, hence has order
$m^-$ and is contained in $\Ga_\H$.  As for Equation
\eqref{eq:defcountfunct3}, disintegrating the counting function of
Equation \eqref{eq:defcountfunct2thm4} by the canonical map $\Ga'_{\wt
  D^-} \bs\Ga/\Ga_{\wt D^+}\ra \Ga_{\H} \bs\Ga/\Ga_{\wt D^+}$,
Equation \eqref{eq:defcountfunct2thm4} gives
\begin{align}
  \N_{\Omega_+,D^+}(s)&=\frac{1}{m^-}\card\bigg\{([\ga],\beta)\in
  (\Ga_{\H}\bs\Ga/\Ga_{\wt D^+})\times\Ga_\H: \begin{array}{l}
    0<d(\wt D^-,\beta\ga \wt D^+)\leq s\\z_{\wt D^-,\,\beta\ga\wt D^+}
    \in \wt\ell(\,\interval[open]{0}{+\infty}\,)\end{array}\bigg\}
  \nonumber\\ &\;\;+\bigO(e^{\delta\,s}) \,.
  \label{eq:defcountfunct3thm4}
\end{align}

As in the beginning of Step \hyperlink{Step7}{7} in the proof of
Theorem \ref{theo:1geoddivreal}, with $s\in
\interval[open]{1}{+\infty}$ large enough, $\eta \in
\interval[open]{0}{1}\,$ small enough and $N=\lfloor \frac{s}{\eta}
\rfloor$, if
\[
\Sigma_s=\card\bigg\{([\ga],\beta)\in
(\Ga_{\H}\bs\Ga/\Ga_{\wt D^+}\ssm F)\times (\Ga_\H\ssm F'): \begin{array}{l}
    s-1<d(\wt D^-,\beta\ga \wt D^+)\leq s\\z_{\wt D^-,\,\beta\ga\wt D^+}
    \in \wt\ell(\,\interval[open]{0}{+\infty}\,)\end{array}\bigg\}\,,
\]
then
\begin{align}\label{eq:doublesomm1thm4}
\Sigma_s=\sum_{k=1}^{N+1}\sum_{\substack{[\ga]\in
    \Ga_{\H}\bs\Ga/\Ga_{\wt D^+}\ssm \,F\\
    (k-1)\eta<d(\ga\wt D^+,\,\H)\leq k\eta}}
\quad\sum_{\substack{\beta\in\Ga_\H\ssm\, F'\\
    s-1<d(\wt D^-,\,\beta\ga \wt D^+)\leq s}} 1\,.
\end{align}
Let $k\in\llbracket1,N+1\rrbracket$ and $([\ga],\beta)\in
(\Ga_{\H}\bs\Ga/\Ga_{\wt D^+}\ssm F)\times (\Ga_\H\ssm F')$ be such
that
\begin{align*}
(k-1)\eta<d(\ga\wt D^+,\H)\leq k\eta\quad\text{and}\quad
s-1 < d(\wt D^-,\beta\ga \wt D^+)\leq s\,.
\end{align*}
Then since $M$ is real hyperbolic, by the first claim of Lemma
\ref{lem:distancetodivline} \hyperlink{lemhypreelii}{(ii)} applied
with $D=\beta\ga\wt D^+$ and $x_0=\beta p_\ga$, we have
\begin{align}
d(\wt D^-,\beta\ga\wt D^+) &=
d(\ga\wt D^+,\H)+\ln (2\,d'_\H(\beta p_\ga,\wt\ell(0)))\nonumber\\ &\quad +
\bigO\big(d'_\H(\beta p_\ga,\wt\ell(0))^{-2}\,e^{-2\,d(\ga\wt D^+,\H)}\big)\,.
\label{eq:prepabootstrapthm4}
\end{align}
As in the middle part of  Step \hyperlink{Step7}{7} in the proof of
Theorem \ref{theo:1geoddivreal},  up to increasing $F'$, we have
\begin{equation}\label{eq:doubleencadminothm4}
\frac{1}{2}\;e^{s-1 -k\eta +\bigO(e^{-2\,s})}\leq
d'_\H(\beta p_\ga,\wt\ell(0))  \leq 
\frac{1}{2}\;e^{s -(k-1)\eta +\bigO(e^{-2\,s})}\;.
\end{equation}
By Equations \eqref{eq:doubleencadminothm4} and
\eqref{eq:horocountthm4}, with functions $\bigO(\cdot)$ independent of
$\eta$ and $[\ga]$, we have
\begin{align}
&\card\{\beta\in\Ga_\H\ssm F':
s-1<d(\wt D^-,\beta\ga \wt D^+)\leq s\}\nonumber\\\leq\;&
  \frac{\Xi_{\wt M}}{2^\delta\|\sigma^+_{\V_+}\|}
  e^{\delta \,s -\delta k\eta}\big(e^{\delta\eta
  +\bigO(e^{-s})}-e^{-\delta+\bigO(e^{-s})}\big)
  +\bigO(e^{(s -k\eta)(\delta -1)})\,.\label{eq:cotebetathm4}
\end{align}
Since $M$ is real hyperbolic and since $\V_+$ has finite outer
skinning measure, let us apply Theorem \ref{theo:1geoddivreal} with
$D^-=\V_+$, and more precisely Equation \eqref{eq:NDDabstrait} with the
help of Step \hyperlink{Step2}{2} in the proof of Theorem
\ref{theo:1geoddivreal} in order to deal with the multiplicities
$m_{\ga \V_+,\ga'\wt D^+}$ not equal to $1$. Then, with a function
$\bigO(\cdot)$ that is independent of $\eta$, if we now define, for
every $k\in\NN$,
\[
a_k=\card\big\{[\ga]\in \Ga_{\H}\bs\Ga/\Ga_{\wt D^+}\ssm F:
  (k-1)\eta<d(\ga\wt D^+,\H)\leq k\eta\big\}\,,
\]
then for every $t\geq 0$, we have
\begin{equation}\label{eq:cotegathm4}
\sum_{k=0}^t\;a_k=
\frac{\iota_{\rm rev}^+\,\|\sigma^+_{\V_+}\|\;\Xi_{\wt M}}
     {2^\delta\;m^+\;\|m_{\rm BM}\|}\;
t\,\eta\; e^{\delta t\eta} +\bigO\big(e^{\delta t\eta}\big)\,.
\end{equation}

Let us now define 
\[
C_1=\frac{\Xi_{\wt M}}{2^\delta\,
  \|\sigma^+_{\V_+}\|} \, (e^{\delta\eta+\bigO(e^{-s})}-
e^{-\delta+\bigO(e^{-s})})\, , \hspace{.7cm}
C_2=\frac{\iota_{\rm rev}^+\,\|\sigma^+_{\V_+}\|\,\Xi_{\wt
    M}}{2^\delta\,m^+\, \|m_{\rm BM}\|}\;,
\] 
and a function $f:t\mapsto C_1\, e^{\delta s-\delta t\eta}+
\bigO(e^{(s-t\eta)(\delta-1)})$ with an appropriately chosen
$\bigO(\cdot)$. As in the middle part of Step \hyperlink{Step7}{7} in
the proof of Theorem \ref{theo:1geoddivreal}, by Equations
\eqref{eq:doublesomm1thm4} and \eqref{eq:cotebetathm4}, by Abel's
summation formula, since $f(N+1)=\bigO(1)$ and by Equation
\eqref{eq:cotegathm4}, using again that $N\eta=s+\bigO(\eta)$ and
since $\int_0^{(N+1)\eta}u\,e^u\,du=\bigO(s\,e^s)$, we have
\begin{align*}
  \Sigma_s&\leq \sum_{k=0}^{N+1} a_k\,f(k)= \big(\sum_{k=0}^{N+1}
  a_k\,\big)f(N+1)- \int_0^{N+1}(\sum_{k=0}^{t}
  a_k\,\big)f'(t)\,dt\\&=\bigO(s\,e^{\delta\,s})+
  \int_0^{N+1}\big(C_2\,t\,\eta\,e^{\delta t\eta}+\bigO(e^{\delta t\eta})\big)
  \big(\delta\,\eta\,C_1\, e^{\delta s-\delta t \eta}+
  \bigO(\eta\, e^{(s-t\eta)(\delta-1)})\big)\,dt\\&=\bigO(s\,e^{\delta\,s})+
  \delta\,C_2\,C_1\,\frac{(N+1)^2}{2}\,\eta^2\,e^{\delta s}+
  \bigO\big((N+1)\,\eta\,e^{\delta\,s}\big)\\&\qquad\qquad\;
  +\bigO\big(e^{s(\delta-1)}\int_0^{N+1}\eta\,t\,e^{t\eta}\,\eta\,dt\big)
  +\bigO\big(e^{s(\delta-1)}\int_0^{N+1}e^{t\eta}\,\eta\,dt\big)
  \\&=\frac{\delta\,C_2\,C_1}{2}\,(s+\bigO(\eta))^2\,e^{\delta s}
  +\bigO(s\,e^{\delta\,s})\,.
\end{align*}
Replacing $C_1$ and $C_2$ by their values, and letting $\eta$ tend to
$0$, we hence have
\[
\Sigma_s\leq
\frac{\delta\;\iota_{\rm rev}^+\;{\Xi_{\wt M}}^2}
     {2^{2\delta+1}\:m^+\;\|m_{\rm BM}\|}
\;s^2\;e^{\delta \,s}(1-e^{-\delta})+\bigO(s\,e^{\delta\, s})\,.
\]

The same lower bound is obtained as at the end of Step
\hyperlink{Step7}{7} in the proof of Theorem \ref{theo:1geoddivreal},
and by a summation, we have
\begin{align}\label{eq:countfunctmainthm4}
  &\card\Big\{([\ga],\beta)\in
(\Ga_{\H}\bs\Ga/\Ga_{\wt D^+}\ssm F)\times (\Ga_\H\ssm F'): \begin{array}{l}
    0<d(\wt D^-,\beta\ga \wt D^+)\leq s\\z_{\wt D^-,\,\beta\ga\wt D^+}
    \in \wt\ell(\,\interval[open]{0}{+\infty}\,)\end{array}\Big\}\;.
  \nonumber\\&=\frac{\delta\;\iota_{\rm rev}^+\;{\Xi_{\wt M}}^2}
{2^{2\delta+1}\;m^+\;\|m_{\rm BM}\|}\;s^2\;e^{\delta \,s}+
  \bigO\big(s\,e^{\delta\, s}\big)\,. 
\end{align}
By separating the counting domain $(\Ga_{\H}\bs\Ga/\Ga_{\wt D^+})
\times \Ga_\H$ as the disjoint union of $F\times\Ga_\H$, of
$(\Ga_{\H}\bs\Ga/\Ga_{\wt D^+}\ssm F)\times F'$ and of
$(\Ga_{\H}\bs\Ga/\Ga_{\wt D^+}\ssm F)\times (\Ga_\H\ssm F')$, Equation
\eqref{eq:NOmegaplusthm4} finally follows from Equations
\eqref{eq:defcountfunct3thm4}, \eqref{eq:countfunctFthm4},
\eqref{eq:countfunctFprimethm4} and \eqref{eq:countfunctmainthm4}.

\medskip\noindent{\bf Case \hypertarget{Case2bis}{2}. }  Let us now
assume that $\ell$ is reversible. As in Case
\hyperlink{Case2}{2} of Theorem \ref{theo:1geoddivreal}, we then have
$\Omega_- =\Omega_+$ and 
\begin{equation}\label{eq:NDDOmegapmrecbis}
\N_{\Omega_+,D^+}(t)\leq \N_{D^-,D^+}(t)\leq \N_{\Omega_+,D^+}(t) +
\N_{\Omega_0,D^+}(t)
\end{equation}
for every $t\geq 0$. Let us prove that Equation
\eqref{eq:NOmegaplusthm4} is still satisfied. Since $\iota_{\rm rev}^-
=1$ as $\ell$ is reversible, this will prove that Equation
\eqref{eq:NDDabstraitthm4} is still satisfied, hence Theorem
\ref{theo:divergentperp2} when $\ell$ is reversible will
follow.

As in Case \hyperlink{Case2}{2} of Theorem \ref{theo:1geoddivreal},
with $\Ga'_{\wt D^-} = \Ga_{\wt D^-} \cap\Ga_\H$, which has index $2$
in $\Ga_{\wt D^-}$ and order $m^-$, as $s\ra+\infty$, we have
\begin{align}
 & \N_{\Omega_+,\,D^+}(s)=
  \sum_{[\ga]\in\Ga_{\wt D^-}\bs\Ga/\Ga_{\wt D^+} \,:\;
    0<d(\wt D^-,\ga \wt D^+)\leq s,\; z_{\wt D^-,\,\ga \wt D^+}\,
    \in\, \wt\ell(\,\interval[open]{-\infty}{t_-}\,
\cup\,\interval[open]{0}{+\infty}\,)} m_{\wt D^-,\ga\wt D^+}\nonumber\\=\;&
  \card\Big\{\ga\in\Ga'_{\wt D^-}\bs\Ga/\Ga_{\wt D^+}:
   \begin{array}{l}0<d(\wt D^-,\ga \wt D^+)\leq s\\z_{\wt D^-,\,\ga\wt D^+}
    \in \wt\ell(\,\interval[open]{0}{+\infty}\,)\end{array}\Big\}+
    \bigO(e^{\delta\,s})\nonumber\\=\;&\frac{1}{m^-}\card\Big\{([\ga],\beta)\in
  (\Ga_{\H}\bs\Ga/\Ga_{\wt D^+})\times\Ga_\H: \begin{array}{l}
    0<d(\wt D^-,\beta\ga \wt D^+)\leq s\\z_{\wt D^-,\,\beta\ga \wt D^+}
    \in \wt\ell(\,\interval[open]{0}{+\infty}\,)\end{array}\Big\}
  +\bigO(e^{\delta\,s})\,,
  \label{eq:defcountfunct1recbis}
\end{align}
that is, Equation \eqref{eq:defcountfunct3thm4} is still valid. The
remainder of the proof of Equation \eqref{eq:NOmegaplusthm4} now
proceeds exactly as in Case \hyperlink{Case1bis}{1}.
\cqfd

\section{Common perpendiculars of divergent geodesics in
  non-real hyperbolic geometry}
\label{sec:complexhypgeom}

In this section, we prove Theorem \ref{theo:1+2geoddivcompquat}, which
is a complex and quaternionic hyperbolic version of Theorems
\ref{theo:1geoddivreal} and \ref{theo:divergentperp2}.  This result
will be applied in \cite{ParPauHD} to study the distribution of
Heisenberg Farey neighbours.

In what follows, we denote by $\KK$ either the field of complex
numbers $\CC$ endowed with the conjugation $x=x_0+ix_1\mapsto
\overline{x}=x_0-ix_1$ or the skew field of Hamiltonian numbers $\HH$
(with standard basis $1,i,j,k$ over $\RR$) endowed with the
conjugation $x=x_0+x_1i+x_2j+x_3k\mapsto \overline{x}=x_0-x_1i-x_2j
-x_3k$.  We refer to \cite{Vigneras80} for background on $\HH$. We
denote by $\Re:x\mapsto \frac{1}{2} (x+\overline{x})$ and
$\Im:x\mapsto \frac{1}{2} (x-\overline{x})$ the real and imaginary
part maps of $\KK$,\footnote{Note the nonstandard definition of $\Im$
when $\KK=\CC$.} so that $\Im \CC=i\RR$ and $\Im\HH=\RR i+\RR j+\RR k$
are the imaginary subspaces of $\CC$ and $\HH$ respectively. We endow
$\KK$ and $\Im\KK$ with the Euclidean scalar product that makes their
canonical basis orthonormal and with its associated Lebesgue
measure. Let $n \in\NN\ssm\{0,1\}$. We endow $\KK^{n-1}$ with the
product Euclidean scalar product and product Lebesgue measure.

For all $w,w'$ in the right vector space $\KK^{n-1}$ over $\KK$, we
denote their standard Hermitian product by $\overline{w}\cdot w'= \sum_{i=1}^{n-1} \overline{w_i}\,
w'_i$, and we define $|w|=
\sqrt{\overline w\cdot w}$.  Recall that the {\em Siegel domain} model
of the hyperbolic $n$-space $\hnk$ over $\KK$ is the open subset
\[
\big\{(w_0,w)\in\KK\times\KK^{n-1}\;:\; 
2\,\Re w_0 -|w|^2>0\big\}\,,
\]
endowed with the Riemannian metric
\begin{equation}\label{eq:siegelmetric}
ds^2_{\,\hnk}=\frac{1}{(2\,\Re w_0 -|w|^2)^2}
\big(\,|\,dw_0-\,\overline{dw}\cdot w\,|^2+
(2\,\Re w_0 -|w|^2)\;|dw|^2\,\big)\,.
\end{equation}
The metric is normalized so that its sectional curvatures are in
$[-4,-1]$, instead of in $[-1,-\frac{1}{4}]$ as in \cite{Goldman99}
when $\KK=\CC$. The boundary at infinity of $\hnk$ is
$$
\partial_\infty\hnk=\big\{(w_0,w)\in
\KK\times \KK^{n-1} \;:\; 2\,\Re w_0 -|w|^2=0\big\}\cup\{\infty\}\,.
$$

As in \cite[\S 3]{ParPau17MA} when $\KK=\CC$ and \cite[\S
  6]{ParPau22MPCPS} when $\KK=\HH$, the {\it horospherical coordinates}
$(\zeta,u,t)\in \KK^{n-1}\times\Im\KK\times\,[0,+\infty[$ of a point
    $(w_0,w)\in\hnk\cup (\partial_\infty \hnk\ssm \{\infty\})$
    are
\begin{equation}\label{eq:horosphecoord}
(\zeta,u,t)  =(w,\;w_0-\overline{w_0}, \;2\,\Re
  w_0-|w|^2)\quad\text{hence}\quad(w_0,w) 
  =\Big(\frac{|\zeta|^2+t+u}{2} ,\zeta\Big)\,.
\end{equation}

In horospherical coordinates, the subset
\[
\H_\infty=\big\{(\zeta,u,t)\in\hnk\;:\;t\geq 1\big\},
\]
is a horoball in $\hnk$ centred at $\infty\in\partial_\infty
\hnk$. The geodesic line in $\hnk$ from $\infty$ to the point at
infinity $(\zeta,u,0)\in\partial_\infty \hnk\ssm \{\infty\} $, through
$\partial \H_\infty$ at time $s=0$, is the map $s\mapsto
(\zeta,u,e^{-2s})$.

The following lemma is a convenient replacement in hyperbolic geometry
over $\KK$ of the classical angle of parallelism formula in real
hyperbolic geometry. See also \cite[\S 3.2.4]{Goldman99} for a
different presentation. The proof follows ideas from \cite[Lemma 10]
{ParPau17MA} when $\KK=\CC$ and \cite[Lemma 6$\cdot$2]{ParPau22MPCPS}
when $\KK=\HH$. See also \cite[Prop.~7.1]{Parker10} for an expression
of the distance from a point to a geodesic line in the projective
model of $\hnc$.

\blemm \label{lem:orthprojgeod} The orthogonal projection map from
$\hnk$ to the geodesic line $\interval[open]0\infty$ in $\hnk$
with points at infinity $(0,0)$ and $\infty$ is, in horospherical
coordinates, the map
\[
(\zeta,u,t)\mapsto (0,0,|\,|\zeta|^2+t+u\,|)\,. 
\]
The distance from $(\zeta,u,t)\in\hnk$ to $\interval[open]0\infty$ is 
$\frac 12\arcosh\big(\frac{|\zeta|^2+|\,|\zeta|^2+t+u\,|}t\big)$.
\elemm

\dem
Let $\BB^n$ be the open unit sphere in the standard right Hermitian space
$\KK^n$ over $\KK$.  The Cayley transform $\Phi:\BB^n\to \hnk$, defined by
\[
\Phi:(z_1,\dots,z_n)\mapsto
\Big(\frac{1-z_n}2,z_1,z_2,\dots,z_{n-1}\Big)(1+z_n)^{-1}\,,
\]
is easily seen to be a smooth bijection, with inverse 
\begin{equation}\label{eq:ballsiegelisom}
(w_0,w)\mapsto (2w,1-2w_0)(1+2w_0)^{-1}\,.
\end{equation}
The {\em ball model} of the hyperbolic $n$-space over $\KK$ is the
open subset $\BB^n$ endowed with the pull-back of the Riemannian
metric \eqref{eq:siegelmetric} by $\Phi$.

Let $\rho>0$. In this ball model, the metric sphere $S(0,\rho)$ of
radius $\rho$ centered at the origin $0$ coincides with the Euclidean
sphere of radius $\tanh \rho$ centered at $0$ by \cite[page 78, see
  also \S 3.3.4]{Goldman99} when $\KK=\CC$, taking into account the
different normalization of the curvatures.

The isometry $\Phi$ maps $0\in\BB^n$ to $(0,0,1)\in\hnk$ in the
horospherical coordinates. By Equation \eqref{eq:ballsiegelisom}, for
all $z'\in\KK^{n-1}$ and $z_n\in\KK$, writing $(w_0,w)$ the point
$\Phi(z',z_n) $ and denoting by $(\zeta,u,t)$ its horospherical
coordinates, we have $|z'|^2+|z_n|^2 =\tanh^2 \rho$ if and only if
\[
|\,2w\,|^2+|\,1-2w_0\,|^2=|\,1+2w_0\,|^2\tanh^2 \rho
\]
that is, using Equation \eqref{eq:horosphecoord} and an easy
computation, if and only if
\begin{equation}\label{eq:sphererho}
|\,1+|\zeta|^2+t+u\,|^2-4t\cosh^2 \rho=0\,.
\end{equation}

The Riemannian metric of $\hnk$ given by Equation
\eqref{eq:siegelmetric} becomes in the horospherical coordinates
\begin{equation}\label{eq:siegelmetrichc}
ds^2_{\,\hnk}=\frac{1}{4\,t^2}
\big(dt^2+ |\,du-2\Im\,\overline{d\zeta}\cdot \zeta\,|^2+
4\,t\;|\,d\zeta\,|^2\big)\,.
\end{equation}
Hence for all $\lambda>0$ and $(\zeta',u')\in\KK^{n-1}\times\Im \KK$, the
{\it Heisenberg dilation}
\[
h_\lambda:(\zeta,u,t)\mapsto (\lambda\zeta,\lambda^2 u,\lambda^2t)\,,
\]
whose inverse is $h_{\lambda^{-1}}$, and the {\it Heisenberg translation}
\[
\tau_{(\zeta',u')}:(\zeta,u,t)\mapsto
(\zeta'+\zeta,\,u'+u+ 2\,\Im\,\overline{\zeta'}\cdot\zeta,\,t)\,,
\]
whose inverse is $\tau_{(-\zeta',-u')}$, are isometries of $\hnk$
fixing $\infty$.

Using the horospherical coordinates, let us fix $(\zeta_0,u_0,t_0)\in
\hnk$, and let us compute its orthogonal projection on the geodesic
line $\interval[open]0\infty$.  Let us consider $\lambda=\sqrt{t_0}$,
$\zeta'=\frac{\zeta_0} {\sqrt{t_0}}$ and $u'=\frac{u_0}{t_0}$. The
isometry $h_\lambda\circ \tau_{(\zeta',u')}$ maps $(0,0,1)$ to
$(\zeta_0,u_0,t_0)$, hence maps the sphere of radius $\rho$ centered
at $(0,0,1)$ to the sphere of radius $\rho$ centered at
$(\zeta_0,u_0,t_0)$ in $\hnk$. For every $(\zeta,u,t)\in \hnk$, using
Equation \eqref{eq:sphererho} (multiplied by $t_0^2$) for the last
equivalence, we have
\begin{align}
  &(\zeta,u,t)\in S((\zeta_0,u_0,t_0),\rho)\nonumber
  \\\Leftrightarrow \quad& \tau_{(\zeta',u')}^{-1}\circ h_\lambda^{-1}
  (\zeta,u,t)\in S((0,0,1),\rho)\nonumber\\\Leftrightarrow \quad&
  \big(-\zeta'+\lambda^{-1}\zeta,-u'+\lambda^{-2}u+
  2\Im (\,\overline{-\zeta'}\,)\cdot(\lambda^{-1}\zeta),\lambda^{-2}t\,\big)
  \in S((0,0,1),\rho)\nonumber\\\Leftrightarrow \quad&
  \big|\,t_0+|\zeta-\zeta_0|^2 +t +(u-u_0 -
  2\Im \,\overline{\zeta_0}\cdot\zeta)\,\big|^2=4t_0t\cosh^2 \rho\,.
  \label{eq:formulesphere}
\end{align}

The closest point to $(\zeta_0,u_0,t_0)$ on $\interval[open]0\infty$
is attained when the parameter $\rho$ gives a double point of
intersection $(0,0,t)$ between this sphere and
$\interval[open]0\infty$\,. Taking $\zeta=0$ and $u=0$ in Equation
\eqref{eq:formulesphere} gives the following quadratic equation in $t$
\[
t^2+2t(|\zeta_0|^2+t_0-2t_0\cosh^2 \rho)+|\,|\zeta_0|^2+t_0+u_0\,|^2=0\,.
\]
It has a double solution if and only if its reduced discriminant
\[
(|\zeta_0|^2+t_0-2t_0\cosh^2 \rho)^2- |\,|\zeta_0|^2 +t_0+u_0\,|^2
\]
is equal to zero, that is, if and only if
\begin{equation}\label{ew:prepavaleurrho}
  |\zeta_0|^2 +t_0-2t_0\cosh^2 \rho=-|\,|\zeta_0|^2+t_0+u_0\,|\,.
\end{equation}
The double solution of the above quadratic equation is then $t=
|\,|\zeta_0|^2 +t_0+u_0\,|$. This proves the first claim of Lemma
\ref{lem:orthprojgeod}.  The second claim follows from Equation
\eqref{ew:prepavaleurrho} by using the fact that
$2\cosh^2\rho=\cosh(2\rho)+1$.  \cqfd

\medskip
Let us now make explicit the Hamenstädt distance and cuspidal distance
on horospheres in $\hnk$. The {\it Heisenberg group $\Heis_{n,\KK}$}
is the real Lie group $\KK^{n-1}\times \Im \KK$ with law
\[
(\zeta',u')(\zeta,u)=
(\zeta'+\zeta,\,u'+u+ 2\,\Im\,\overline{\zeta'}\cdot\zeta)\,.
\]
As defined for instance in \cite[page 160]{Goldman99} when $\KK=\CC$,
the {\it Cygan distance} $d_{\rm Cyg}$ on $\Heis_{n,\KK}$ is the
unique left-invariant distance on $\Heis_{n,\KK}$ such that
\begin{equation}\label{eq:defidCyg}
d_{\rm Cyg}\big((\zeta,u),(0,0)\big)=\sqrt{|\,|\zeta|^2+u\,|\,}
=\sqrt[4]{|\zeta|^4+|u|^2}\,.
\end{equation}
As introduced in \cite[page 372]{ParPau10GT} when $\KK=\CC$, the
{\it modified Cygan distance} $d'_{\rm Cyg}$ on $\Heis_{n,\KK}$ is the unique
left-invariant distance on $\Heis_{n,\KK}$ such that
\begin{equation}\label{eq:defidprimCyg}
d'_{\rm Cyg}\big((\zeta,u),(0,0)\big)=\sqrt{|\zeta|^2+|\,|\zeta|^2+u\,|\,}
=\sqrt{|\zeta|^2+\sqrt{|\zeta|^4+|u|^2}\,}\,.
\end{equation}
It is easy to check that $d_{\rm Cyg}\le d'_{\rm Cyg}\le \sqrt 2\;
d_{\rm Cyg}$.  As in \cite[page 216]{HerPau02b} and \cite[page
  370]{ParPau10GT} both when $\KK=\CC$, for every $t'>0$, using
Equation \eqref{eq:horosphecoord} for the second equality, let
\begin{equation}\label{eq:defihoroHtprim}
H_{t'}=
\{(w_0,w)\in\hnk:2\Re w_0-|w|^2\geq t'\}=\{(\zeta,u,t)\in\hnk:t\geq t'\}\,,
\end{equation}
which is a horoball in $\hnk$ centered at $\infty$, so that
$H_1=\H_\infty$.  The Heisenberg group $\Heis_{n,\KK}$ acts, by the
map $(\zeta',u')\mapsto {\tau_{(\zeta',u')}}_{\mid H_{t'}}$, simply
transitively on the horosphere $H_{t'}$ for every $t'>0$, as well as
on $\partial_\infty\hnk \ssm \{\infty\}$. Let us prove the claim that for all
$t''\geq t'$ and $(\zeta,u), (\zeta',u')\in\KK^{n-1}\times\Im \KK$, we
have
\begin{equation}\label{eq:variatdisthamen}
d_{H_{t'}}((\zeta,u,t'),(\zeta',u',t'))=
\sqrt{\frac{t''}{t'}}\;d_{H_{t''}}((\zeta,u,t''),(\zeta',u',t''))\,.
\end{equation}
\dem Since $t''\geq t'$, the horoball $H_{t''}$ is contained in
$H_{t'}$.  An immediate computation using the geodesic line
$s\mapsto (0,0,e^{2s})$ yields $d(\partial H_{t''},\partial
H_{t'})=\frac{1}{2}\ln\frac{t''}{t'}$.  By the definition of the
Hamenstädt distance in Equation \eqref{eq:defihamendist}, we have as
wanted
\[
d_{H_{t'}}((\zeta,u,t'),(\zeta',u',t'))=e^{d(\partial H_{t''},\partial H_{t'})}
\;d_{H_{t''}}((\zeta,u,t''),(\zeta',u',t''))\,.\quad\Box
\]

\medskip
By Equation \eqref{eq:variatdisthamen} applied with $t'=1$ and $t''=2$
and by \cite[Prop.~3.12]{HerPau02b} (which uses the horosphere
$\partial H_2$ instead of the horosphere $\partial
H_1=\partial\H_\infty$) when $\KK=\CC$, and by a similar computation
when $\KK=\HH$, we have
\begin{equation}\label{eq:vald}
  d_{\H_\infty}((\zeta,u,1),(\zeta',u',1))=
  \sqrt{2}\;d_{H_2}((\zeta,u,2),(\zeta',u',2))=
d_{\rm Cyg}((\zeta,u),(\zeta',u'))\,.
\end{equation}
By \cite[Prop.~6.2]{ParPau10GT} applied with $s_0=1$ (so that the
horoball $H_1$ of loc. cit. is equal to our horoball $\H_\infty$) when
$\KK=\CC$, and by a similar computation
when $\KK=\HH$, we have
\begin{equation}\label{eq:valdprim}
d'_{\H_\infty}((\zeta,u,1),(\zeta',u',1))=\frac{1}{\sqrt{2}}\;
d'_{\rm Cyg}((\zeta,u),(\zeta',u'))\,.
\end{equation}

\blemm\label{lem:distancetodivlineKK} Let $\H$ be a horoball in $\hnk$
and let $\wt\ell$ be a geodesic line in $\hnk$ that enters $\H$
perpendicularly at $\wt\ell(0) \in\partial\H$.

\smallskip\noindent \hypertarget{lemhypKKi}{(i)} Let $\wt\ell'$ be a
geodesic line in $\hnk$ that exits $\H$ perpendicularly at
$\wt\ell'(0)\in\partial\H$ such that $d_\H(\wt\ell'(0),\wt\ell(0))
\geq 1$.  For every $s\ge 0$, we have
\[
d(\wt\ell'(s),\wt\ell\,)=s+\ln d'_\H(\wt\ell'(0),\wt\ell(0))+\ln 2
+\bigO\big(d'_\H(\wt\ell'(0),\wt\ell(0))^{-2}e^{-2s}\big)\,.
\]

\smallskip\noindent\hypertarget{lemhypKKii}{(ii)} Let $D$ be a closed
convex subset of $\hnk$ disjoint from $\H$ and let $x_0\in\partial\H$
be the closest point to $D$ in $\H$. There exists a constant $c_\KK\geq 1$
such that if  $d_\H(x_0,\wt\ell(0))\geq c_\KK$, then
\[
d(D,\wt\ell\,)=d(D,\H)+\ln d'_\H(x_0,\wt\ell(0)) +\ln 2 +
\bigO\big(d'_\H(x_0,\wt\ell(0))^{-2}e^{-2\,d(D,\,\H)}\big)\,,
\]
and furthermore, the closest point to $D$ on the image of $\wt\ell$
belongs to $\H$.
\elemm

\dem We use the horospherical coordinates of $\hnk$. The isometry
group of $\hnk$ acts transitively on the set of horoballs of $\hnk$
and the stabilizer of each horoball acts transitively on its boundary
horosphere. Hence we may assume that $\H=\H_\infty$ and that
$\wt\ell(0) = (0,0,1)$. Therefore the geodesic line $\wt \ell$ is the
map $s\mapsto (0,0,e^{2s})$ and its image is $\interval[open]0\infty$
with the notation of Lemma \ref{lem:orthprojgeod}.

\medskip
\hyperlink{lemhypKKi}{(i)} Let us define $(\zeta,u)\in\KK^{n-1}\times
\Im \KK$ such that the geodesic line $\wt \ell'$ is the map given by
$s\mapsto (\zeta,u,e^{-2s})$. Using Equations \eqref{eq:defidCyg} and
\eqref{eq:vald}, let
\[
D=|\,|\zeta|^2+u\,|=d_{\rm Cyg}((\zeta,u),(0,0))^2= d_\H(\wt\ell'(0),
\wt\ell(0))^2\,.
\]
Using Equations \eqref{eq:defidprimCyg} and \ref{eq:valdprim}, let
\begin{equation}\label{eq:defiDprim}
D'=|\zeta|^2+|\,|\zeta|^2+u\,|=d'_{\rm Cyg}((\zeta,u),(0,0))^2=
2\;d'_\H(\wt\ell'(0),\wt\ell(0))^2\,.
\end{equation}
Since $|\zeta|^2\leq\sqrt{|\zeta|^4+|u|^2}=D$ and since $D\geq 1$ by
the assumption of Assertion \hyperlink{lemhypKKi}{(i)}, we have
\begin{align*}
  |\,|\zeta|^2+e^{-2s}+u\,|&
  =\sqrt{|\zeta|^4+2e^{-2s}|\zeta|^2+e^{-4s}+|u|^2}
  \\&=\sqrt{|\zeta|^4+|u|^2}\;\sqrt{1+\frac{2e^{-2s}|\zeta|^2}
    {|\zeta|^4+|u|^2}+\frac{e^{-4s}}{|\zeta|^4+|u|^2}}
  \\&= D\sqrt{1+\bigO(e^{-2s}D^{-1})}=D(1+\bigO(e^{-2s}D^{-1}))\,.
\end{align*}
Since $D\leq |\zeta|^2+D= D'\leq 2D$, we hence have
\begin{align*}
  x&=\frac{|\zeta|^2+|\,|\zeta|^2+e^{-2s}+u\,|}{e^{-2s}}
  =e^{2s}D'\frac{|\zeta|^2+D + D\bigO(e^{-2s}D^{-1})}{D'}\\&=
  e^{2s}D'\big(1+\bigO(e^{-2s}{D'}^{-1})\big)\,.
\end{align*}
Recall that as $x\in\interval[open right]{1}{+\infty}$ tends to
$+\infty$, we have
\[
\arcosh x=\ln(x+\sqrt{x^2-1}\,)=
\ln(2\,x)+\bigO\big(\frac{1}{x^2}\big)\,.
\]
By the last claim of Lemma \ref{lem:orthprojgeod}, we therefore have
\begin{align*}
  d(\wt\ell'(s),\wt\ell\,)&=\frac 12\arcosh
  \Big(\frac{|\zeta|^2+|\,|\zeta|^2+e^{-2s}+u\,|}{e^{-2s}}\Big)
  =s+\frac 12\log (2D') +\bigO(e^{-2s}{D'}^{-1})\,.
\end{align*}
Assertion \hyperlink{lemhypKKi}{(i)} then follows from Equation
\eqref{eq:defiDprim}.

\medskip
\hyperlink{lemhypKKii}{(ii)} As in the proof of Lemma
\ref{lem:distancetodivline} \hyperlink{lemhypreelii}{(ii)}, let
$x_{\wt\ell}\in D$ be the closest point in $D$ to $\wt\ell$, and let
$x_{\H}\in D$ be the closest point in $D$ to $\H$. Let $\wt\ell'$ be
the geodesic line exiting $\H$ perpendicularly at $x_0$ at time
$0$. Let $s=d(D,\H)$, so that $x_0=\wt\ell'(0)$ and $x_{\H}=
\wt\ell'(s)$, and let $D''=d'_\H(x_0,\wt \ell(0))$. Finally, let
$p_{x_{\wt\ell}}$ (respectively $p_{x_{\H}}$) be the closest point to
$D$ (respectively to $x_{\H}$) on the image of $\wt\ell$. We have the
upper bound
\[
d(x_{\wt\ell},\wt\ell\,)=d(D,\wt\ell\,)\leq d(x_\H,\wt \ell\,)
=s+ \ln D''+\ln 2+\bigO({D''}^{-2}e^{-2s})
\]
by Assertion \hyperlink{lemhypKKi}{(i)} In order to obtain the similar
lower bound on $d(x_{\wt\ell},\wt\ell\,)$, as in the proof of Lemma
\ref{lem:distancetodivline} \hyperlink{lemhypreelii}{(ii)} (except
that the union of the geodesic lines perpendicular to $\wt\ell'$ at
$x_\H$ is no longer totally geodesic), we may replace $D$ by a
geodesic line $D'$ through $x_{\wt\ell}$ and $x_{\H}$ perpendicular to
$\wt\ell'$ at $x_\H$. See Figure
\addtocounter{fig}{1}\hyperlink{pointprocheKK}{\arabic{fig}}\addtocounter{fig}{-1}.

\begin{center}
  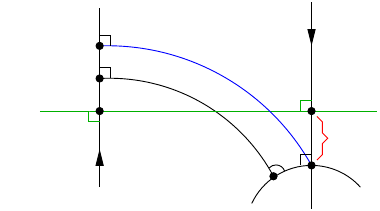
\end{center}
\vspace*{-0.8cm}\addtocounter{fig}{1}
\begin{center}
  {\small Figure \arabic{fig}~:
    \hypertarget{pointprocheKK}{Computations of distances to a geodesic
      line}.}
\end{center}
By Assertion \hyperlink{lemhypKKi}{(i)}, we have $d(x_{\H},\wt\ell\,)
=s+\ln D'' +\ln 2+\bigO(e^{-2(s+\ln D'')})$.  Since the closest point
projection to a convex subset does not increase the distances, we have
the inequality $d(p_{x_{\wt\ell}}, p_{x_{\H}})\leq d(x_{\wt\ell},
x_{\H})$.  By the triangle inequality, we have
\[
d(x_{\H},p_{x_{\H}})-d(p_{x_{\wt\ell}}, p_{x_{\H}})- d(x_{\wt\ell}, x_{\H})
\leq d(D,\wt \ell)= d(x_{\wt\ell},p_{x_{\wt\ell}})
\leq d(x_{\H},p_{x_{\H}})\,.
\]
To obtain Assertion \hyperlink{lemhypKKii}{(ii)}, we only have to
prove that $d(x_{\wt\ell},x_{\H})=
\bigO(e^{-2\,d(x_{\H},\,\wt\ell\,)})$.  In $\hnr$, this is the new
input of Lemma \ref{lem:distancetodivline}
\hyperlink{lemhypreelii}{(ii)} with respect to Lemma
\ref{lem:distancetodivline} \hyperlink{lemhypreeli}{(i)}. In $\hnk$,
the result follows by comparison using the ideal quadrangle with
vertices $\infty$, $p_{x_{\wt\ell}}$, $x_{\wt\ell}$ and $x_{\H}$ with
right angles at $p_{x_{\wt\ell}}$ and $x_{\H}$, and angle at least
$\frac{\pi}{2}$ at $x_{\wt\ell}$.
\cqfd

\medskip
We now prove analogs of Theorems \ref{theo:1geoddivreal} and
\ref{theo:divergentperp2} in the complex or quaternionic hyperbolic
case. In order to simplify the notation, we define $d_\KK=\dim_\RR\KK$
for $\KK=\CC$ and $\KK=\HH$. 

\btheo\label{theo:1+2geoddivcompquat} Let $M$ be a noncompact finite
volume complete connected complex or quaternionic hyperbolic good
orbifold, with dimension $n\geq 2$ over $\KK=\CC$ or $\KK=\HH$.

\smallskip\noindent \hypertarget{theocompquat1}{(1)} Let $D^-$ be a
nonempty properly immersed closed locally convex subset of $M$ with
nonzero finite outer skinning measure and let $D^+$ be the image of
a divergent geodesic in $M$. As $s\ra+\infty$, we have
\[
\N_{D^-,D^+}(s)= \frac{\prod_{i=1}^{\frac{d_\KK}{2}}(\frac{nd_\KK}{2}-i)\;
  \iota_{\rm rev}(D^+)\;\|\sigma^+_{D^-}\|}{4^{d_\KK-1}\;\sqrt{\pi}\;
  \Ga(\frac{d_\KK-1}{2})\,m(D^+)\Vol(M)}\;
  s\, e^{(d_\KK(n+1)-2) \,s} + \bigO(e^{(d_\KK(n+1)-2) s})\,.
\]
  
\smallskip\noindent \hypertarget{theocompquat2}{(2)} Let $D^+$
and $D^-$ be the images of two divergent geodesics in $M$.
As $s\to+\infty$, we have
\begin{align*}
\N_{D^-,D^+}\!(s)\!=&\frac{(d_\KK(n+1)-2)\;\pi^{\frac{nd_\KK}{2}-1}
  \prod\limits_{i=1}^{\frac{d_\KK}{2}}(\frac{nd_\KK}{2}-i)\;
  \iota_{\rm rev}(D^-)\;\iota_{\rm rev}(D^+)}
  {2^{d_\KK(n+3)-4}\,\Ga(\frac{d_\KK-1}{2})^2\,(\frac{d_\KK(n-1)}{2}-1)!\,
    m(D^-)m(D^+)\Vol M} s^2\, e^{(d_\KK(n+1)-2) \,s}\\ &
  + \bigO\big(s\;e^{(d_\KK(n+1)-2) s}\big)\,.
\end{align*}
\etheo

We believe that a similar statement is valid also for the octonionic
hyperbolic plane case, but we leave the proof to the readers.

\medskip
\dem This proof follows closely the proofs of Theorem
\ref{theo:1geoddivreal} for Assertion \hyperlink{theocompquat1}{(1)}
and of Theorem \ref{theo:divergentperp2} for Assertion
\hyperlink{theocompquat2}{(2)}, that were written for this purpose. We
only indicate the changes, that are the ones involving specifically
the fact that $M$ was assumed to be real hyperbolic, and no longer is.
We start with a lemma, that will replace Equation
\eqref{eq:Xirealhyp}.

\blemm\label{lem:computXinonreel}
We have $\Xi_{\wt M}=\frac{2\;\pi^{\frac{n\,d_\KK-1}{2}}}
{\Ga(\frac{d_\KK-1}{2})\;(\frac{d_\KK(n-1)}{2}-1)!}$.
\elemm

\dem Since the definition \eqref{eq:defiXi} of $\Xi_{\wt M}$ is
independent of the choices of a horoball $\H$ and of a point
$x\in\partial \H$, we may assume that $\H=\H_\infty$ and that
$x=(0,0,1)$ in the horospherical coordinates of $\hnk$.  By
\cite[Lemma 12 (iv)]{ParPau17MA} when $\KK=\CC$ and by \cite[Lemma
  7$\cdot$2 (iv)]{ParPau22MPCPS} when $\KK=\HH$, the measure $(\wt
p_\bullet)_*\,\wt\sigma^-_{\H_\infty}$ is $2^{d_\KK-1}$ times the
Riemannian measure $\vol_{\partial \H_\infty}$ of the induced
Riemannian metric on $\partial \H_\infty$. Since $\H_\infty=H_1$ with
the notation of Equation \eqref{eq:defihoroHtprim}, by \cite[Equation
  (15)]{ParPau17MA} when $\KK=\CC$ and by \cite[Equation
  (7$\cdot$11)]{ParPau22MPCPS} when $\KK=\HH$, we have
\[
d\vol_{\partial \H_\infty}(\zeta,u,1)=\frac{1}{2^{d_\KK-1}}\;d\zeta\,du\,.
\]
By Equations \eqref{eq:valdprim} and \eqref{eq:defidprimCyg}, we have
\begin{align*}
  B_{d'_{\H_\infty}}((0,0,1),1)&=\{(\zeta,u,1)\in \partial \H_\infty:
  d'_{\rm Cyg}((\zeta,u),(0,0))\leq \sqrt{2}\}\\&
  =\{(\zeta,u,1)\in \partial \H_\infty:
  |\zeta|^2+\sqrt{|\zeta|^4+|u|^2}\leq 2\}\,.
\end{align*}
Hence by the definition \eqref{eq:defiXi} of $\Xi_{\wt M}$ and by
using the formulas $du=\rho \,d\rho\,d\vol_{\SSS^{d_\KK-2}}$ and
$d\zeta = \rho \,d\rho\,d\vol_{\SSS^{d_\KK(n-1)-1}}$ of the Lebesgue
measures in polar coordinates of the Euclidean spaces $\Im\KK$ and
then $\KK^{n+1}$, we have
\begin{align*}
  \Xi_{\wt M}& =(\wt p_\bullet)_*\,\wt\sigma^-_{\H_\infty}
  \big(B_{d'_{\H_\infty}}((0,0,1),1))\big)
  =\int_{(\zeta,u,1)\in B_{d'_{\H_\infty}}((0,0,1),1)} \;d\zeta\,du\\&=
  \int_{|\zeta|\leq 1}\Big(\int_{|u|\leq 2\sqrt{1-|\zeta|^2}}du\Big)d\zeta=
  \Vol(\SSS^{d_\KK-2})\int_{|\zeta|\leq 1}2(1-|\zeta|^2)\;d\zeta\\&=
  \frac{1}{2}\,\Vol(\SSS^{d_\KK-2})\,\Vol(\SSS^{d_\KK(n-1)-1})\,.
\end{align*}
Lemma \ref{lem:computXinonreel} follows from this equation, as $\Vol(\SSS^{m-1})=
\frac{2\;\pi^{\frac{m}{2}}} {\Ga(\frac{m}{2})}$ and $\Ga(m')=(m'-1)!$
for all $m,m'\in\NN\ssm \{0\}$.
\cqfd

\medskip
By \cite[Lemma 12 (iii)]{ParPau17MA} when $\KK=\CC$ and by \cite[Lemma
  7$\cdot$2 (iii)]{ParPau22MPCPS} when $\KK=\HH$, we now have
\begin{equation}\label{eq:valmBM}
  \|m_{\rm BM}\|= \frac{1}{2^{d_\KK(n-1)}}\;\Vol(\SSS^{n\,d_\KK-1})\;\Vol(M)
  = \frac{\pi^{\frac{n\,d_\KK}{2}}}
  {2^{d_\KK(n-1)-1}(\frac{n\,d_\KK}{2}-1)!}\;\Vol(M)\,.
\end{equation}

In order to prove Assertion \hyperlink{theocompquat1}{(1)} of
Theorem \ref{theo:1+2geoddivcompquat}, we use the
same notation as in the beginning of the proof of Theorem
\ref{theo:1geoddivreal}. The discussion on whether $\ell$ is reversible or not is the same one as in the proof of Theorem
\ref{theo:1geoddivreal}. If Equation \eqref{eq:NOmegaplus} is still
valid, then Equation \eqref{eq:NDDabstrait} is also still valid, by
the same proof. Assertion \hyperlink{theocompquat1}{(1)} of
Theorem \ref{theo:1+2geoddivcompquat} follows from Equation
\eqref{eq:NDDabstrait} by using Lemma \ref{lem:computXinonreel}
instead of Equation \eqref{eq:Xirealhyp}, by using Equation
\eqref{eq:valmBM} instead of Equation \eqref{eq:PP17ETDS20_1}, and by
using Equation \eqref{eq:valexpcrit} instead of $\delta=n-1$.

The proof of Equation \eqref{eq:NOmegaplus} when $\wt M=\hnk$ follows
the same seven steps as in the real hyperbolic case, except that in
Steps \hyperlink{Step3}{3} and \hyperlink{Step7}{7}, the use of Lemma
\ref{lem:distancetodivline} \hyperlink{lemhypreelii}{(ii)} with the
constant $c_\RR=2$ is replaced by the use of Lemma
\ref{lem:distancetodivlineKK} \hyperlink{lemhypKKii}{(ii)} with the
constant $c_\KK$.

\medskip
In order to prove Assertion \hyperlink{theocompquat2}{(2)} of Theorem
\ref{theo:1+2geoddivcompquat}, we use the same notation as in the
beginning of the proof of Theorem \ref{theo:divergentperp2}. The
discussion on whether $\ell$ is reversible or not is the same
one as in the proof of Theorem \ref{theo:divergentperp2}. If Equation
\eqref{eq:NOmegaplusthm4} is still valid, then Equation
\eqref{eq:NDDabstraitthm4} is also still valid, by the same
proof. Assertion \hyperlink{theocompquat2}{(2)} of Theorem
\ref{theo:1+2geoddivcompquat} follows from Equation
\eqref{eq:NDDabstraitthm4} by using Lemma \ref{lem:computXinonreel}
instead of Equation \eqref{eq:Xirealhyp}, by using Equation
\eqref{eq:valmBM} instead of Equation \eqref{eq:PP17ETDS20_1}, and by
using Equation \eqref{eq:valexpcrit} instead of $\delta=n-1$.

The proof of Equation \eqref{eq:NDDabstraitthm4} when $\wt M=\hnk$ is
similar to its proof when $\wt M=\hnr$, replacing the call to Theorem
\ref{theo:1geoddivreal} in the proofs of Equations
\eqref{eq:NOmegazerothm4}, \eqref{eq:NVD} and \eqref{eq:cotegathm4} by a
call to Equation \eqref{eq:NDDabstrait} that we just proved during the
proof of Assertion \hyperlink{theocompquat1}{(1)} of Theorem
\ref{theo:1+2geoddivcompquat}.  \cqfd

\section{Ambiguous geodesics}
\label{sec:ambi}

In this section, we first show that the ambiguous conjugacy classes of
hyperbolic elements of the modular group $\PSL_2(\ZZ)$ discussed in
\cite{Sarnak07} correspond to common perpendiculars of divergent
geodesics in the modular orbifold $\PSL_2(\ZZ)\bs \hdr$.  We then use
Theorem \ref{theo:1geoddivreal} and Theorem \ref{theo:divergentperp2}
to recover, by hyperbolic geometry methods, asymptotic counting
results of special conjugacy classes in $\PSL_2(\ZZ)$, due to Sarnak
\cite{Sarnak07} by arithmetic methods.  We start by recalling standard
facts on the modular orbifold, and on the images of the imaginary axis
by the modular group.

As in Section \ref{sec:realhypgeom}, let $\hdr\subset \CC$ be the
upper halfplane model of the real hyperbolic plane, so that
$\partial_\infty\hdr= \PP_1(\RR)= \RR\cup \{\infty\}$. Given $\xi\neq
\eta$ in $\partial_\infty\hdr$, we denote by
$\interval[open]{\xi}{\eta}$ the geodesic line in $\hdr$ with points
at infinity $\xi$ and $\eta$. We denote by $\begin{bsmallmatrix}
a&b \\ c&d \end{bsmallmatrix}$ the image in $\PGL_2(\CC)$ of $
\begin{psmallmatrix}a&b\\c&d \end{psmallmatrix}\in\GL_2(\CC)$.
The group $\PSLR$ acts isometrically and faithfully by homographies on
$\hdr$, by the map $(\ga,z)\mapsto \ga\cdot z= \frac{az+b}{cz+d}$ for
$z\in\hdr$ and $\ga= \begin{bsmallmatrix} a&b\\ c&d
\end{bsmallmatrix}\in\PSLR$. Let $\Ga_\ZZ= \PSL_2(\ZZ)$ be the
{\em modular group}, which is a nonuniform arithmetic lattice in
$\PSLR$ with set of parabolic fixed points $\operatorname{Par}
_{\Ga_\ZZ} = \Ga_\ZZ\cdot \infty=\QQ\cup \{\infty\}$.  Let
$M_\ZZ=\Ga_\ZZ\bs\hdr$ be the {\em modular orbifold}, which is a
non\-compact complete connected real hyperbolic good orbifold of
volume $\frac\pi 3$ with only one cusp.

Let $\wt\ell:t\mapsto i\,e^t$ be the geodesic line in $\hdr$ through
$i\in\hdr$ at time $t=0$, with endpoints at infinity $0$ and $\infty$,
and let $\wt\Delta=\wt\ell(\RR)=i\RR\cap \hdr=\;
\interval[open]{0}{\infty}$ be its image. Then $\ell=\Ga_\ZZ\,\wt
\ell$ is a divergent geodesic in $M_\ZZ$, converging at $\pm\infty$ to
the only cusp $\Ga_\ZZ\cdot \infty$ of $M_\ZZ$. Note that $\ell$ is
strongly reversible since $\wt \Delta$ is preserved by the involution $\iota=
\begin{bsmallmatrix} 0&-1\\1&\ \; 0\end{bsmallmatrix} \in
\Ga_\ZZ$ with fixed point set $\{i\}$ in $\hdr$.  The stabilizer in
$\Ga_\ZZ$ of $\wt\Delta$ is $\Ga_{\wt \Delta}= \{\id,\iota\}$.

Similarly, let $\wt\ell_1:t\mapsto \frac{2} {1-2ie^{-t}}$, which is
the geodesic line in $\hdr$ image of $\wt\ell$ by
$\begin{bsmallmatrix} 1&0\\1/2&1\end{bsmallmatrix}$, with
endpoints at infinity $0$ and $2$, and let $\wt\Delta_1=\wt\ell_1(\RR)
=\;]0,2[$ be its image. Then $\ell_1= \Ga_\ZZ\,\wt \ell_1$ is a also a
divergent geodesic in $M_\ZZ$, which is also strongly reversible since $1+i
\in(\Ga_\ZZ \cdot i)\cap \wt\Delta_1$.  The stabilizer in $\Ga_\ZZ$ of
$\wt\Delta_1$ is $\Ga_{\wt \Delta_1}= \big\{\id,\begin{bsmallmatrix}
-1&2\\-1&1\end{bsmallmatrix} \big\}$.  But $\wt\Delta_1$ is not
the image of $\wt\Delta$ by any element of $\Ga_\ZZ$.

Let $\Delta=\ell(\RR)$ be the image of $\ell$ and $\Delta_1= \ell_1
(\RR)$ the one of $\ell_1$. Let $\wt D^-,\wt D^+\in \{\wt\Delta,
\wt\Delta_1\}$ and $D^-,D^+$ their images in $M_\ZZ$.  The action of
$\Ga_\ZZ$ on $T^1\hdr$ is free, hence the multiplicities of the common
perpendiculars from $D^-$ to $D^+$ (defined in Equation
\eqref{eq:defmultipli}) are all equal to $1$. Thus Equation
\eqref{eq:defcountingfunct} gives that
\[
\N_{D^-,D^+}(t)=\card\{[\ga]\in\Ga_{\wt D^-}\bs\Ga_\ZZ/\Ga_{\wt D^+}:
0< d(\wt D^-,\ga\wt D^+)\leq t\}
\]
is the number of images under $\Ga_\ZZ$ of $\wt D^+$ that are at
positive distance at most $t$ from $\wt D^-$, modulo the left action
of $\Ga_{\wt D^-}$.  Let ${\hdr}^\pm=\{z\in\hdr: \pm \Re(z)>0\}$.
Noting that $\iota{\hdr}^\pm={\hdr}^\mp$, we see that
$\N_{\Delta,\Delta} (t)$ equals the number of images under $\Ga_\ZZ$
of $\wt\Delta$ that are contained in ${\hdr}^+$ and at a positive
distance at most $t$ from the imaginary axis.

\blemm \label{lem:travaildivgeodmodular} Let $\ga=\begin{bsmallmatrix}
a&b\\c&d\end{bsmallmatrix}\in\Ga_\ZZ$.

\smallskip\noindent \hypertarget{tradivgeodmod1}{(i)} The geodesic
lines $\wt\Delta$ and $\ga\wt\Delta$ have a common point at infinity if
and only if $a\,b\,c\,d=0$. They have a common perpendicular if and
only if $a\,b\,c\,d\ne0$.

\smallskip\noindent\hypertarget{tradivgeodmod2}{(ii)} The geodesic line
$\ga\wt\Delta$ is contained in the right halfplane ${\hdr}^+$ and has
no common point at infinity with $\wt \Delta$ if and only if
$\ga\Ga_{\wt\Delta}$ has a representative $\begin{psmallmatrix}
  a&b\\c&d \end{psmallmatrix} \in\SL_2(\ZZ)$ with $a,b,c,d>0$.  This
representative is then unique.
\elemm

\dem \hyperlink{tradivgeodmod1}{(i)} Note that
$\ga\cdot\infty=\frac{a}{c}$ and $\ga\cdot0= \frac{b}{d}$. The
intersection $\partial_\infty(\ga\wt\Delta\,)\cap \partial_\infty
\wt\Delta$ is nonempty if and only if $\ga\cdot\infty= \infty$ or
$\ga\cdot 0=0$ or $\ga\cdot0=\infty$ or $\ga\cdot\infty=0$, that is, if
and only if $c=0$ or $b=0$ or $d=0$ or $a=0$ respectively, which
proves the first part of Assertion \hyperlink{tradivgeodmod1}{(i)}.

The geodesic lines $\wt\Delta$ and $\ga\wt\Delta$ have a common
perpendicular if and only if they do not have a common point at
infinity and do not have a common point inside $\hdr$.  Let
$\Ga_\infty$ be the stabilizer in $\Ga_\ZZ$ of the horoball
$\H_\infty$ defined in Equation \eqref{eq:defiHinftyreal} for
$n=2$. It is well known that the $\Ga_\ZZ$-equivariant family
$(\ga\H_\infty)_{\ga\in\Ga_\ZZ/\Ga_\infty}$ is precisely invariant.
The only horoballs in this family containing $i$ are $\H_\infty$ and
$\iota\H_\infty$, and $\wt \Delta$ is contained in $\H_\infty\cup
\iota\H_\infty$. Hence for every $\ga\in\Ga_\ZZ$, the geodesic lines
$\wt \Delta$ and $\ga\wt\Delta$ meet if and only if
$\ga\in\Ga_{\wt\Delta}= \{\id,\iota\}$, in which case $\wt \Delta$ and
$\ga\wt\Delta$ have a common point at infinity.

Thus $\wt\Delta$ and $\ga\wt\Delta$ have a common perpendicular if and
only if the points $\ga\cdot \infty = \frac{a}{c}$ and $\ga\cdot0=
\frac{b}{d}$ are in the same connected component of $\RR\ssm\{0\}$.
This yields the inequality $\frac{a}{c}\frac{b}{d}>0$, which proves
the second part of Assertion \hyperlink{tradivgeodmod1}{(i)} by
multiplying by $(cd)^2>0$.

\smallskip\noindent\hyperlink{tradivgeodmod2}{(ii)} The above
computations show that $\ga\wt\Delta$ is contained in the right
halfplane ${\hdr}^+$ and has no common point at infinity with $\wt
\Delta$ if and only if $ac>0$ and $bd>0$. Since $a\neq 0$, the element
$\ga$ has a unique representative in $\SL_2(\ZZ)$ such that $a>0$, and
then $c>0$. Assume from now on that $a>0$. Note that $
\begin{psmallmatrix} a&b\\c&d \end{psmallmatrix}\begin{psmallmatrix}
  \ \; 0&1\\-1&0\end{psmallmatrix} = \begin{psmallmatrix} -b&a\\-d&c
\end{psmallmatrix}$. Hence if $b>0$, then $d>0$ and $
  \begin{psmallmatrix} a&b\\c&d \end{psmallmatrix}$ is the unique
representative in $\SLZ$ with positive coefficients of an element of
$\ga\Ga_{\wt \Delta}$.  And if conversely $b<0$, then $d<0$ and
$\begin{psmallmatrix} -b&a\\-d&c \end{psmallmatrix}$ is the unique
representative in $\SLZ$ with positive coefficients of an element of
$\ga\Ga_{\wt \Delta}$.
\cqfd

\medskip
Now, let us define the special conjugacy classes in $\Ga_\ZZ$ that we
will study. Let
\[
w=\begin{bmatrix} 1&\ \ 0\\0&-1\end{bmatrix}
\quad\text{and}\quad
w_1=\begin{bmatrix} 1&\ \ 0\\1&-1\end{bmatrix}\,,
\]
which are involutions (elements of order $2$) in $\PGL_2(\ZZ)$. Recall
that $\PGL_2(\ZZ)$ acts by conjugation on its normal subgroup $\Ga_\ZZ
= \PSL_2(\ZZ)$. An element $\ga\in \Ga_\ZZ$ is {\em ambiguous of the
  first kind}, respectively {\it ambiguous of the second kind}, if
\begin{equation}\label{eq:ambi1or2}
w\ga w=\ga^{-1}\,, \quad\textrm{ respectively }\quad w_1\ga w_1=\ga^{-1}\;, 
\end{equation}
and {\em ambiguous} if it is conjugated in $\Ga_\ZZ$ to an element in
$\Ga_\ZZ$ which is ambiguous of the first kind or ambiguous of the
second kind.  Such elements, when hyperbolic, are automorphs of Gauss'
ambiguous integral binary quadratic forms, see \cite{Sarnak07} and
\cite[Sect.~14.4]{Cassels78} for details and background. Recall that a
hyperbolic element $\ga\in \Ga_\ZZ$ is {\em primitive} if it is not a
proper power in $\Ga_\ZZ$. Every hyperbolic element $\ga\in\Ga_\ZZ$
has a unique {\it root}, i.e.~a primitive element $\ga_0\in\Ga_\ZZ$
such that there exists $n\in\NN \ssm \{0\}$ with $\ga=\ga_0^n$.  For a
hyperbolic element of $\Ga_\ZZ$, being ambiguous, ambiguous of the
first kind or ambiguous of the second kind is invariant by taking
nonzero powers and roots.

\medskip
We say that conjugacy classes of hyperbolic elements of $\Ga_\ZZ$ that
contain ambiguous elements of both the first and the second kinds are
{\em biambiguous}.  The existence of such conjugacy classes is not
immediate from the arithmetic definitions but it follows naturally
using reflections, as seen in the proof of Proposition
\ref{lem:ambichar} \hyperlink{ambichariii}{(iii)} below. If a
conjugacy class of ambiguous elements is not biambiguous, it is a {\em
  monoambiguous} conjugacy class.

We will compute the asymptotic growth of the number of biambiguous
conjugacy classes with translation length at most $s\ra+\infty$ during
the proof of Theorem \ref{theo:ambi}, see Lemma \ref{lem:1non2}
(\hyperlink{lem:1non2_i}{i}) and Equation \eqref{eq:growthprimesansi}:
It is equal to $\frac{3} {8\,\pi^2}\; s^2\, e^{\frac{s}{2}}+
\bigO\big(s\; e^{\frac{s}{2}} \big)$, hence it is not negligible with
respect to the number of those containing only ambiguous elements of a
given kind.

\medskip
The normalizer of $\Ga_\ZZ$ in the full isometry group of $\hdr$
contains the reflection
\[
W:z\mapsto -\,\overline  z
\]
in the geodesic line $\wt\Delta$.  The {\em extended modular group}
$\Ga_\ZZ^+$ is the group generated by $\Ga_\ZZ$ and $W$. It
contains $\Ga_\ZZ$ as a normal subgroup of index $2$.  The two
extensions $\PGL_2(\ZZ)$ and $\Ga_\ZZ^+$ of $\Ga_\ZZ$ are actually
isomorphic, see \cite{Beardon15} for a detailed discussion. Let
$\overline{\stackrel{\;\;}{\cdot}}: z\mapsto \overline z$ be the
complex conjugation. The map $\Phi : \PGL_2(\ZZ)\to\Ga_\ZZ^+$, which
is the identity on $\Ga_\ZZ$ and maps $\eta\in\PGL_2(\ZZ)\ssm\Ga_\ZZ$
to the anti-homography $\eta\circ \overline{\stackrel{\;\;}{\cdot}} :
\hdr\ra\hdr$, is a group isomorphism, that is compatible with the
actions on $\Ga_\ZZ$ by conjugation of the two groups: For all
$\ga\in\Ga_\ZZ$ and $\eta\in\PGL_2(\ZZ)$, we have
\[
\Phi(\eta)\,\ga\,\Phi(\eta)^{-1}=\eta\,\ga\,\eta^{-1}\,.
\]
If $\eta\in\PGL_2(\ZZ)\ssm\Ga_\ZZ$ is an involution, then $\Phi(\eta)$
is a reflection in the geodesic line whose endpoints are the two fixed
points of $\eta$ in $\PP_1(\RR)$. In particular, we have $\Phi(w)=W$,
and $W_1=\Phi(w_1)$ is the reflection in the geodesic line $\wt
\Delta_1 =\interval[open]02$. The group $\Ga_\ZZ^+$ has exactly three
conjugacy classes of involutions, which are the conjugacy classes of
$W$, of $W_1$ and of the orientation-preserving involution $\iota$.
Given an involution $\tau\in\Ga_\ZZ^+$, as in \cite{ErlParPau}, we say
that $\ga\in\Ga_\ZZ^+$ is $\tau${\em -reversible} in $\Ga_\ZZ^+$ if
$\tau\ga\tau= \ga^{-1}$. We denote by $\Ax_\ga$ the translation axis
of every hyperbolic element $\ga\in\Ga_\ZZ$.

\blemm\label{lem:ambichar} Let  $\ga=\begin{bsmallmatrix}a&b\\c&d
\end{bsmallmatrix}\in\Ga_\ZZ$.

\smallskip\noindent\hypertarget{ambichari}{(i)} The element $\ga$ is
ambiguous if and only if there exists an involution $\tau\in\Ga_\ZZ^+
\ssm\Ga_\ZZ$ such that $\ga$ is $\tau$-reversible in $\Ga_\ZZ^+$.

\smallskip\noindent\hypertarget{ambicharii}{(ii)} The element $\ga\in
\Ga_\ZZ$ is ambiguous of the first kind (respectively ambiguous of the
second kind) if and only if $a=d$ (respectively $a+b=d$). If $\ga\in
\Ga_\ZZ$ is hyperbolic, then $\ga$ is ambiguous of the first kind
(respectively ambiguous of the second kind) if and only if $\Ax_\ga$
meets perpendicularly $\wt \Delta=\Fix(W)$ (respectively $\wt \Delta_1
= \Fix(W_1)$), and $\ga$ is $\iota$-reversible if and only if $\Ax_\ga$
contains $\{i\}=\Fix(\iota)$.

\smallskip\noindent\hypertarget{ambichariii}{(iii)} The only elements
of $\Ga_\ZZ$ that are both ambiguous of the first kind and ambiguous
of the second kind are $\begin{bsmallmatrix}1&0\\ c&1
\end{bsmallmatrix}$ for $c\in\ZZ$, which are not hyperbolic. 
There exist biambiguous conjugacy classes of hyperbolic elements of
$\Ga_\ZZ$

\smallskip\noindent\hypertarget{ambichariv}{(iv)} A primitive
hyperbolic element of a monoambiguous conjugacy class of $\Ga_\ZZ$ has
exactly $4$ conjugates that are ambiguous of the given kind. A
primitive hyperbolic element of a biambiguous conjugacy class of
$\Ga_\ZZ$ has exactly $2$ conjugates that are ambiguous of the first
kind and $2$ conjugates that are ambiguous of the second kind.
\elemm

\dem \hyperlink{ambichari}{(i)} Let $\ga\in\Ga_\ZZ$ be ambiguous. Then
there exist $\nu\in\Ga_\ZZ$ and $w'\in\{w,w_1\}$ such that
$w'(\nu\ga\nu^{-1}) w'=(\nu\ga\nu^{-1})^{-1}$. Thus, we have
\[
\ga^{-1}=(\nu^{-1}w'\nu)^{-1}\ga (\nu^{-1}w'\nu)=
\Phi(\nu^{-1}w'\nu)^{-1}\ga \,\Phi(\nu^{-1}w'\nu)\,,
\]
and $\ga$ is $\tau$-reversible, with $\tau=\Phi(\nu^{-1}w'\nu)\in
\Ga_\ZZ^+\ssm\Ga_\ZZ$.  The converse  is proven similarly.

\medskip
\hyperlink{ambicharii}{(ii)} The first claim follows by an easy
computation. A reflection $\tau$ of $\hdr$ fixes pointwise a unique
geodesic line $L_\tau$ of $\hdr$ and acts by $v\mapsto -v$ on the unit
normal bundle of $L_\tau$. Hence a geodesic line $L$ of $\hdr$ is
preserved by $\tau$ if and only if $L$ meets $L_\tau$, otherwise
$\tau(L)$ and $L$ are in different connected components of $\hdr\ssm
L_\tau$. For every $x\in L\cap L_\tau$, denoting by $\pm v_x$ the
two unit tangent vectors to $L$ at $x$, either $v_x$ is also tangent
to $L_\tau$ at $x$ so that $\tau(\pm v_x)=\pm v_x$ or $v_x$ is normal
to $L_\tau$ so that $\tau(\pm v_x)=\mp v_x$. Hence an involution
$\tau\in\Ga_\ZZ^+$ preserves a geodesic line $L$ in $\hdr$ if and only
if either $\tau\in\Ga_\ZZ$ and $L$ contains the singleton
$\Fix(\tau)$, or $\tau\in\Ga_\ZZ^+\ssm\Ga_\ZZ$ and either $L$
intersects the geodesic line $\Fix(\tau)$ perpendicularly or
$L=\Fix(\tau)$.  The second claim follows since for all
$\beta,\ga\in\Ga_\ZZ^+$ with $\ga$ hyperbolic, we have
$\Ax_{\beta\ga\beta^{-1}}=\beta\Ax_\ga$, and $\ga$ and $\ga^{-1}$
translate in opposite directions on $\Ax_{\ga^{-1}}=\Ax_\ga$.

\medskip
\hyperlink{ambichariii}{(iii)} The first claim follows easily from the
first claim of Assertion \hyperlink{ambicharii}{(ii)}.  We prove the
second claim by giving an explicit example. Let $p= 
\begin{bsmallmatrix}2&1\\3&2 \end{bsmallmatrix}\in\Ga_\ZZ$. Then
$p$ is hyperbolic, and ambiguous of the first kind by the first claim
of Assertion \hyperlink{ambicharii}{(ii)}. The element
$\nu=\begin{bsmallmatrix} \ \;0&1\\-1&3\end{bsmallmatrix} \in\Ga_\ZZ$
maps the geodesic line $\wt\Delta_1= \interval[open] 02$ to
$\interval[open] {\frac 13}1$. The reflection in the geodesic line
$\interval[open] {\frac 13}1$ is hence $\nu W_1\nu^{-1}\in \Ga_\ZZ^+
\ssm\Ga_\ZZ$. We have $(\nu W_1\nu^{-1})\, p\, (\nu W_1 \nu^{-1})
=p^{-1}$ by an easy computation, so that $\nu^{-1}p\,\nu$ is ambiguous
of the second kind.
\begin{center}
\begin{overpic}[width=.6\textwidth]{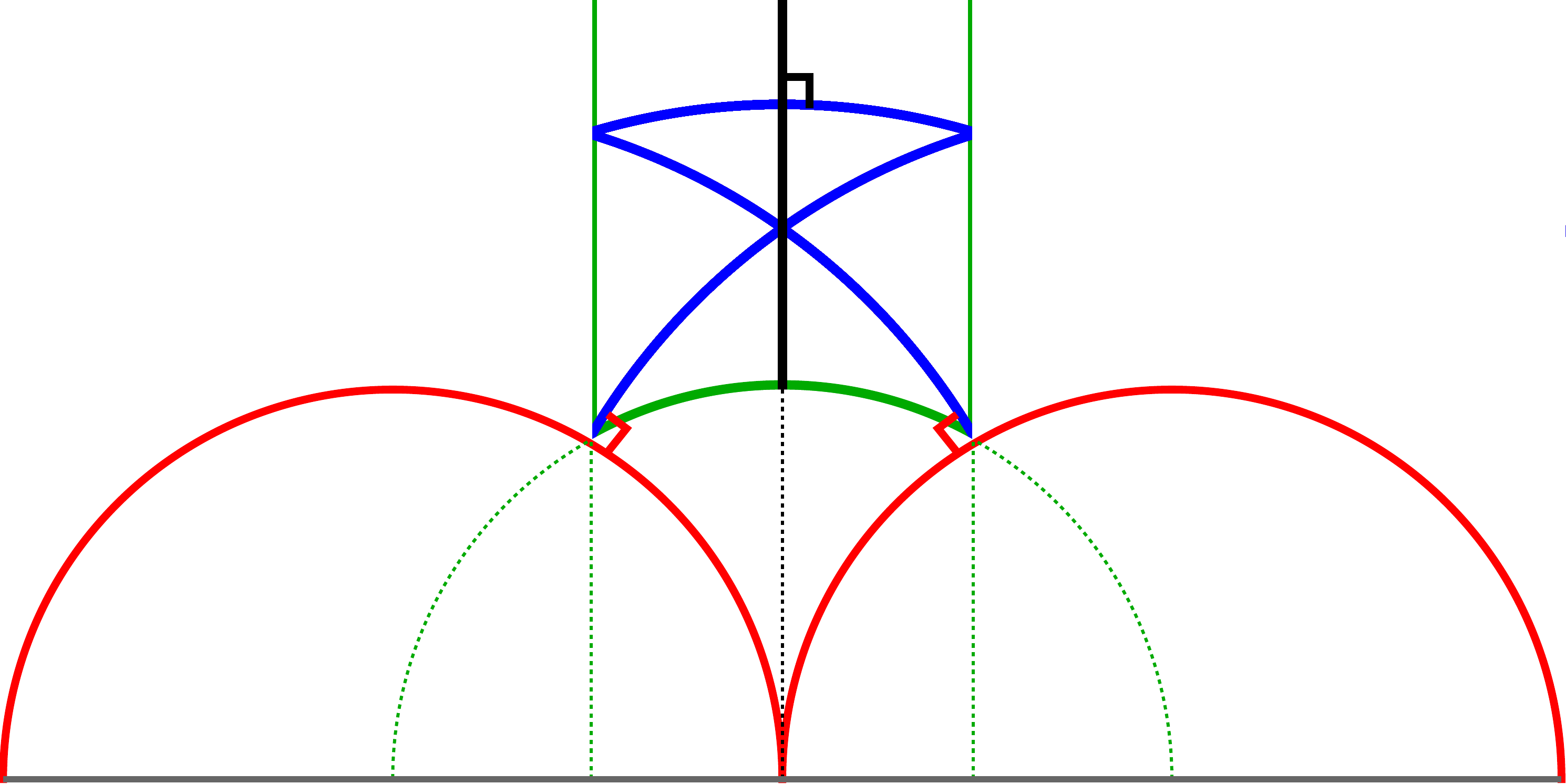}
 \put (-4,-4) {$-2$}
 \put (33,-4){$-\frac 12$}
  \put (48.5,-4) {$0$}
  \put (60.5,-4){$\frac 12$}
    \put (99,-4) {$2$}
      \put (45.5,49) {$\wt\Delta$}
      \put (73,26.5) {$\wt\Delta_1$}
\end{overpic}
\qquad\includegraphics[width=2.5cm]{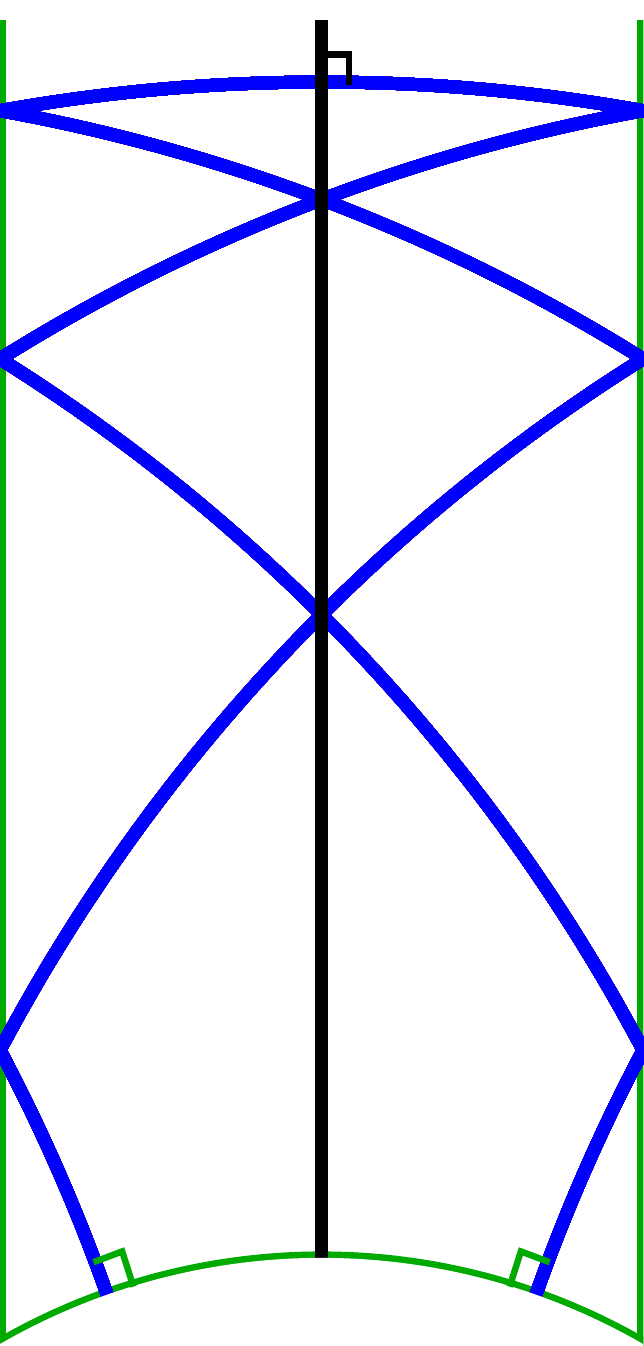}
\end{center}\vspace*{-0.1cm}\addtocounter{fig}{1}
\begin{center}
{\small Figure \arabic{fig}~: \hypertarget{biambiguous}{Biambiguous
    conjugacy classes}.}
\end{center}

\smallskip \noindent Figure
\hyperlink{biambiguous}{\arabic{fig}} on the left shows in blue the
image in $M_\ZZ$ of the translation axis of $p$, lifted to the
standard fundamental domain of $\Ga_\ZZ$ with its usual boundary
identifications, that is orthogonal in $M_\ZZ$ both to $\Delta_1$ (at
the cone point of angle $\frac{2\pi}{3}$ of the orbifold $M_\ZZ$) and
to $\Delta$. Figure \hyperlink{biambiguous}{\arabic{fig}} on the right
shows similarly the image in $M_\ZZ$ of the translation axis of
$p'=\begin{bsmallmatrix}3&8 \\ 1&3\end{bsmallmatrix}$, which is the
union of two common perpendiculars between $\Delta$ and the image of 
$\interval[open]{-1}{+1}$ in $M_\ZZ$.  Noting that $\interval[open]{-1}{+1}$ is
the image of $\wt\Delta_1$ by the element $z\mapsto z-1$ of $\Ga_\ZZ$,
the element $p'$ is hence also both (conjugated to) an ambiguous
element of the first kind and conjugated to an ambiguous element of
the second kind.

\medskip
\hyperlink{ambichariv}{(iv)} Let $\ga\in\Ga_\ZZ$ be primitive
hyperbolic, ambiguous of the first kind and not conjugated in
$\Ga_\ZZ$ to an element ambiguous of the second kind. By Lemma
\ref{lem:ambichar} \hyperlink{ambicharii}{(ii)}, let $x$ be the
perpendicular intersection point of $\wt \Delta$ with $\Ax_\ga$. The
element $\ga W\in\Ga_\ZZ^+\ssm\Ga_\ZZ$ is the reflection in the
perpendicular bisector $\wt M$ of the segment $[x,\ga\cdot x]$. Let $m$ be the
midpoint of $[x,\ga\cdot x]$. Since the involutions of
$\Ga_\ZZ^+\ssm\Ga_\ZZ$ are conjugated by elements of $\Ga_\ZZ$ to
either $W$ or $W_1$, let $\beta\in\Ga_\ZZ$ be such that $\wt
M=\beta\cdot \wt \Delta$ or $\wt M=\beta\cdot \wt \Delta_1$. The
second possibility does not occur, otherwise $\beta^{-1}\Ax_\ga$ would
meet perpendicularly $\wt \Delta_1$, and $\beta^{-1}\ga\beta$ would be
ambiguous of the second kind. Note that $\beta$ is unique up to right
multiplication by $\iota$ since $\Ga_{\wt \Delta}=\{\id,\iota\}$.

\begin{center}
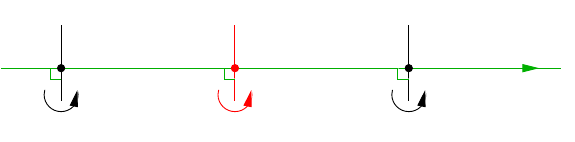
\end{center}
\vspace*{-0.8cm}\addtocounter{fig}{1}
\begin{center}
{\small Figure \arabic{fig}~:
  \hypertarget{multbiambiguous}{Multiplicities in biambiguous
    conjugacy classes}.}
\end{center}

Let $\alpha\in\Ga_\ZZ$ be such that $\alpha\ga\alpha^{-1}$ is
ambiguous of the first kind. Up to replacing $\alpha$ by a right
multiple by a power of $\ga$, which does not change $\alpha\ga
\alpha^{-1}$, we may assume that $\alpha^{-1}\cdot\wt \Delta$ meets
perpendicularly $\Ax_\ga$ in $z\in[x,\ga\cdot x[\,$. We claim that
$\alpha^{-1} \cdot\wt \Delta= \wt\Delta$ or $\alpha^{-1}\cdot\wt \Delta=
\beta\cdot\wt\Delta$.  Otherwise, if $z\in\; ]x,m[\,$, then
$(\alpha^{-1} W\alpha) W$ would be an hyperbolic element with same
translation axis as $\ga$, but with translation length $2\,d(x,z)<
d(x,\ga\cdot x)$, contradicting the fact that $\ga$ is primitive.  And
similarly if $z\in \;]m,\ga\cdot x[\,$, then $(\alpha^{-1}W\alpha)
(\beta W\beta^{-1})$ would be an hyperbolic element with same
translation axis as $\ga$ and translation length $2\,d(m,z)<d(x,
\ga\cdot x)$.  Therefore $\alpha^{-1}\in\{\id,\iota, \beta, \beta
\iota\}$, which proves the first claim, upon checking that $\ga$,
$\iota\ga\iota$, $\beta^{-1}\ga\beta$ and $\iota\beta^{-1}\ga\beta
\iota$ are pairwise distinct, since the centralizer of $\ga$ in
$\Ga_\ZZ$ is $\ga^\ZZ$.

The second claim is proven similarly, except that now $\wt M=\beta
\cdot \wt \Delta_1$, and the conjugates of $\ga$ that are ambiguous of
the first kind are $\ga$ and $\iota\ga\iota$, and the conjugates of
$\ga$ that are ambiguous of the second kind are $\beta^{-1}\ga\beta$
and $\iota\beta^{-1}\ga\beta\iota$.
\cqfd

\brema\label{rem:ambi2} Let $\ga=\begin{bsmallmatrix} a&b\\c&d
\end{bsmallmatrix} \in\Ga_\ZZ$ with $a,b,c,d>0$. The composition
\[
[\ga,W]= (\ga W\ga^{-1})W=\begin{bmatrix}
ad+bc&2ab\\2cd&ad+bc\end{bmatrix}\in\Ga_\ZZ
\]
of the reflections $W$ in $\wt\Delta$ and $\ga W\ga^{-1}$ in
$\ga\wt\Delta$ is ambiguous of the first kind by Lemma \ref{lem:ambichar}
\hyperlink{ambicharii}{(ii)} or since $W[\ga,W]W= [W,\ga]=
[\ga,W]^{-1}$.  By Lemma \ref{lem:travaildivgeodmodular}, the geodesic
lines $\wt\Delta$ and $\ga\wt\Delta$ are disjoint, hence $[\ga,W]$ is
hyperbolic. Its translation axis is the geodesic line that contains
the common perpendicular from $\wt\Delta$ to $\ga\wt\Delta=
\interval[open]{\frac bd}{\frac ac}$, which is $\interval[scaled,open]
{-\sqrt{\frac{ab}{cd}}\;}{\sqrt{\frac{ab}{cd}}\,}$ by an easy
computation (see Figure
\addtocounter{fig}{6}\hyperlink{computcommperp}{\arabic{fig}}\addtocounter{fig}{-6} in the proof of Lemma
\ref{lem:travaildivgeodmodularbis} with $x=\frac{b}{d}$ and
$y=\frac{a}{c}$). Geometrically, this implies that any common
perpendicular from $\Delta$ to itself in $M_\ZZ$ can be extended to a
closed geodesic in $M_\ZZ$, of length twice the length of the common
perpendicular.

\begin{center}
  \includegraphics[width=4.4cm]{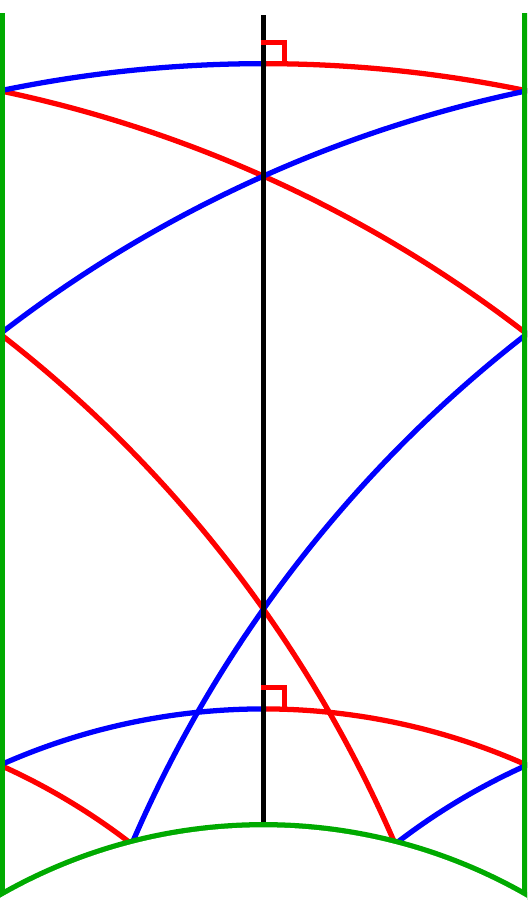}\hspace*{2cm}
  \includegraphics[width=2.2cm]{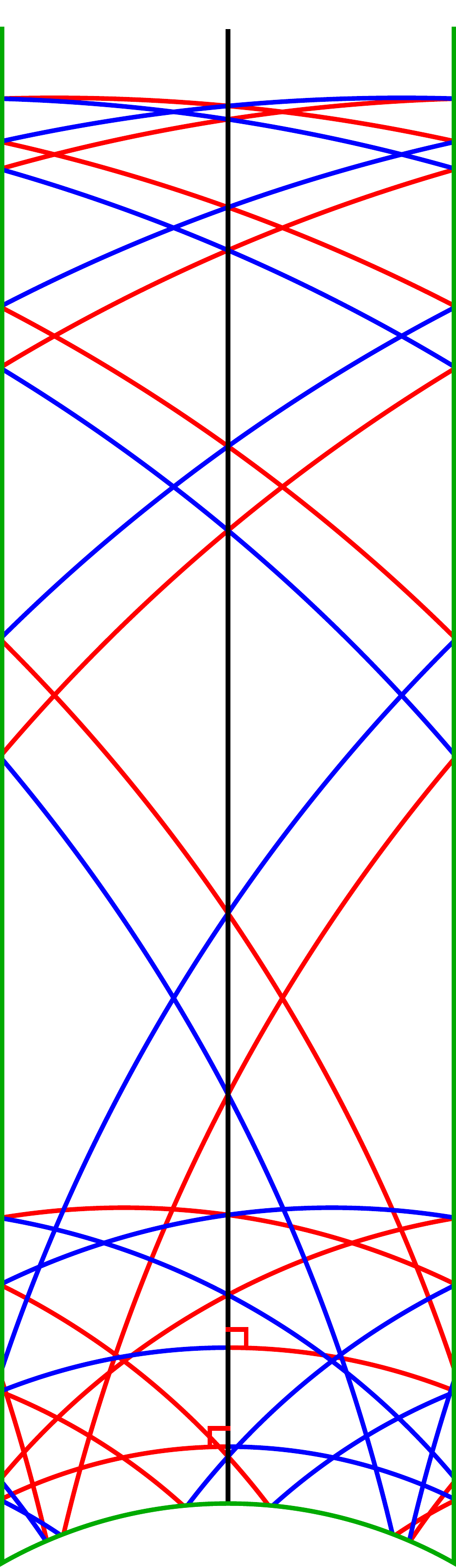}  
\end{center}
\vspace*{-0.8cm}\addtocounter{fig}{1}
\begin{center}
  {\small Figure \arabic{fig}~: \hypertarget{recipclosgeod}{Ambiguous
      geodesics of the first kind on the modular curve}.}
\end{center}
Figure \hyperlink{recipclosgeod}{\arabic{fig}} on the left
(respectively right) shows the common perpendicular from $\Delta$ to
itself in $M_\ZZ$, lifted to the standard fundamental domain of
$\Ga_\ZZ$ with its usual boundary identifications, constructed as
above with $\ga=\begin{bsmallmatrix} 3&1\\2&1\end{bsmallmatrix}$
(respectively $\ga=\begin{bsmallmatrix}
31&23\\35&26\end{bsmallmatrix}\,$), and its reflection along $\Delta$,
doubling it to a closed geodesic of $M_\ZZ$.

Figure
\addtocounter{fig}{-10}\hyperlink{commperp}{\arabic{fig}}\addtocounter{fig}{10}
in the introduction shows the common perpendiculars from $\Delta$ to
itself constructed as above by elements $\ga$ with $bc\le 300$ and
$2.05\le\sqrt{\frac{ab}{cd}}\le 2.1$. The initial tangent vectors of
these common perpendiculars have their footpoints in the standard
fundamental domain contained in the interval
$\interval{2.05\,i}{2.1\,i}$ in the positive imaginary axis.  The
mustard yellow geodesic defined as above by $\ga=\begin{bsmallmatrix}
25&23\\12&11\end{bsmallmatrix}\,$ meets
$\interval[open]{-\frac{1}{2}}{\infty}=\iota\wt\Delta_1$
perpendicularly at its highest point, indicating that the element
$[\ga,W]$ is also conjugate to an ambiguous element of the second kind
by Lemma \ref{lem:ambichar} \hyperlink{ambicharii}{(ii)}.  \erema

The following theorem is the main result of this section. Recall that
a closed geodesic (oriented but not pointed, understood in the
orbifold sense in particular concerning its length) in $M_\ZZ$ is
ambiguous if it meets perpendicularly $\Delta$ or $\Delta_1$, and
reciprocal if it contains $\Ga_\ZZ i$.  Geometrically, such a closed
geodesic then ``reflects'' in $\Delta$, $\Delta_1$, or $\Ga_\ZZ i$
when mapped to $\Ga^+_\ZZ\bs\hdr$.

\btheo\label{theo:ambi} The number $\N_{A}(s)$ of ambiguous primitive
closed geodesics of length at most $s$ satisfies, as $s\ra+\infty$,
that
\[
\N_{A}(s)=\frac{3}{4\,\pi^2}\;s^2\,e^{\frac s2}+\bigO(s\;e^{\frac s2})\,.
\]
The number $\N_{AR}(s)$ of primitive closed geodesics of length at
most $s$ that are both ambiguous and reciprocal satisfies, as
$s\ra+\infty$, that
\[
\N_{AR}(s)=\frac{3}{8\,\pi}\;s\;e^{\frac s4}+\bigO(e^{\frac s4})\,.
\] 
\etheo

The map from the set of conjugacy classes of primitive hyperbolic
elements of $\Ga_\ZZ$ to the set of primitive (oriented but not
pointed) closed geodesic of the modular orbifold $M_\ZZ$, which maps
such a conjugacy class $\llbracket\ga\rrbracket$ to the image in
$M_\ZZ$ of the oriented geodesic segment $[a,\ga a]$ for any $a\in
\Ax_\ga$, is well-known to be a bijection. By Lemma \ref{lem:ambichar}
\hyperlink{ambicharii}{(ii)}, it sends ambiguous/reciprocal conjugacy
classes to ambiguous/reciprocal closed geodesics. Recall that the
absolute value $x$ of the trace of any representative in $\SL_2(\ZZ)$
of any element of the class $\llbracket\ga\rrbracket$ and the length
$s$ of the associated closed geodesic satisfy, when large, that
$x=2\cosh\frac{s}{2}\sim e^{\frac{s}{2}}$. We hence recover, in the
two claims of Theorem \ref{theo:ambi}, respectively Equation (12) and
Equation (15) of \cite{Sarnak07}, up to a multiplicative constant,
possibly coming from the fact that the equality in Equation (61) in
\cite{Sarnak07} seems incorrect by Lemma \ref{lem:ambichar}
\hyperlink{ambichariii}{(iii)}.

\medskip
\dem For every $s> 0$, let us denote by $\operatorname{AC}(s)$ the set
of ambiguous conjugacy classes of primitive hyperbolic elements of
$\Ga_\ZZ$, by $\operatorname{AC}^{1\,\ssm\, 2}(s)$ the ones containing
an ambiguous element of the first kind, but no ambiguous element of
the second kind, by $\operatorname{AC} ^{2\,\ssm \,1} (s)$ the ones
containing an ambiguous element of the second kind, but no ambiguous
element of the first kind, and by $\operatorname{AC}^{1\&2}(s)$ the
set of biambiguous conjugacy classes.  Let $\operatorname{ARC}(s)$,
$\operatorname{ARC}^{1\,\ssm \,2}(s)$, $\operatorname{ARC}^{2\,\ssm
  \,1}(s)$, $\operatorname{ARC} ^{1\&2}(s)$ be the intersection of
these sets with the set of conjugacy classes of reciprocal primitive
hyperbolic elements of $\Ga_\ZZ$.

\blemm\label{lem:non12rec}
The set $\operatorname{ARC}^{1\&2}(s)$ is empty.
\elemm

\begin{center}
    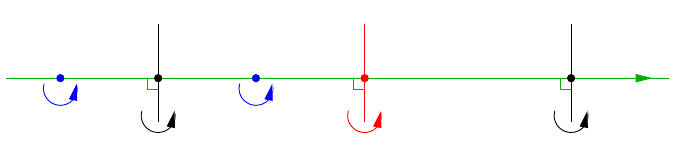
\end{center}
\vspace*{-0.8cm}\addtocounter{fig}{1}
\begin{center}
{\small Figure \arabic{fig}~: \hypertarget{ARCempty}{Contradicting the
    existence of biambiguous and reciprocal conjugacy classes}.}
\end{center}
\dem Assume for a contradiction that there exists a primitive
hyperbolic element $\ga\in\Ga_\ZZ$ which is ambiguous of the first
kind, conjugated to an element ambiguous of the second kind and
reciprocal. By Lemma \ref{lem:ambichar} \hyperlink{ambicharii}{(ii)},
there exists $\alpha,\beta\in\Ga_\ZZ$ such that $\Ax_\ga$ meets
perpendicularly $\wt \Delta$ at a point $x$, meets perpendicularly
$\alpha\cdot\wt \Delta_1$ at a point $y$ and contains the point
$z=\beta \cdot i$. Up to multiplying $\alpha$ and $\beta$ on the left
by powers of $\ga$, we may assume that $y,z\in[x,\ga\cdot x[\,$.  Note
that $x\neq y$ since $\wt \Delta$ and $\wt \Delta_1$ are not in the
same $\Ga_\ZZ$ orbit, that $x\neq z$ since $(\Ga_\ZZ\cdot i)\cap \wt
\Delta=\{i\}$ and the orthogonal geodesic line
$\interval[open]{-1}{1}$ to $\wt \Delta$ at $i$, whose endpoints are
rational, cannot be a translation axis of an element of $\Ga_\ZZ$, and
similarly that $y\neq z$. Let $\ell(\ga)$ be the translation length of
$\ga$, which is the minimal translation length of an hyperbolic
element of $\Ga_\ZZ$ whose translation axis is $\Ax_\ga$ since $\ga$
is primitive. We claim that $d(x,y)\geq \frac{\ell(\ga)}{2}$.
Otherwise, the element $(\alpha W_1\alpha^{-1}) W$, composition of the
reflections $W$ in $\wt \Delta$ and $\alpha W_1\alpha^{-1}$ in $\alpha
\wt \Delta_1$, would belong to $\Ga_\ZZ$, and be hyperbolic with
translation axis $\Ax_\ga$ and translation distance $2\,d(x,y)<
\ell(\ga)$, a contradiction. Similarly, we have $d(\ga \cdot x,y)\geq
\frac{\ell(\ga)}{2}$ and hence $y$ is the midpoint of the segment
$[x,\ga\cdot x]$.

As $W$ fixes $i\in\wt\Delta$ and $W_1$ fixes $1+i\in\wt\Delta_1$, and
since $W$ and $W_1$ normalize $\Ga_\ZZ$, we have $W\Ga_\ZZ\cdot i=
\Ga_\ZZ\cdot i$ and $W_1\Ga_\ZZ\cdot i= W_1\Ga_\ZZ\cdot(1+i)= \Ga_\ZZ
\cdot (1+i)=\Ga_\ZZ\cdot i$. Since $\alpha W_1\alpha^{-1}$ is the
reflection along $\alpha\cdot\Delta_1$, both segments
$\interval[open]{x}{y}$ and $\interval[open]{y}{\ga\cdot x}$ contain a
point of $\Ga_\ZZ\cdot i$. Hence we may assume that $z=\beta\cdot
i\in\interval[open]{x}{y}$, so that $d(x,z)< d(x,y)$.  Consider
$(\beta\iota \beta^{-1})(W\beta\iota \beta^{-1}W)$, the composition of
the angle $\pi$ hyperbolic rotations $W\beta\iota \beta^{-1}W$ around
$W\cdot z=W\beta\cdot i$ and $\beta\iota \beta^{-1}$ around
$z=\beta\cdot i$. It belongs to $\Ga_\ZZ$, and is hyperbolic with
translation axis $\Ax_\ga$ and translation distance $2\,d(x,z)< 2
\,d(x,y)=\ell(\ga)$, a contradiction.  \cqfd

\medskip
Let $\wt D^-,\wt D^+\in \{\wt\Delta,\wt\Delta_1,\{i\}\}$ and let
$D^-,D^+$ be their images in $M_\ZZ$. For every $s>0$, we denote by
$\Perp'(D^-,D^+,s)$ the set of common perpendiculars in $M_\ZZ$
between $D^-$ and $D^+$ that are {\it primitive}, i.e.~that do not
meet perpendicularly in their interior $\Delta,\Delta_1,\{\Ga_\ZZ\cdot
i\}$ (with the convention that an open geodesic segment meets
a point  perpendicularly if and only if it contains the point).

\blemm\label{lem:1non2} (\hypertarget{lem:1non2_i}{i}) The map
$\Psi^{1\&2}$ from $\operatorname{AC} ^{1\&2} (s)$ to $\Perp'(\Delta,
\Delta_1, \frac{s}{2})$ which maps the conjugacy class of an ambiguous
primitive hyperbolic element of the first kind $\ga\in\Ga_\ZZ$ to the
image in $M_\ZZ$ of the oriented geodesic segment $[x,m]$ where $x$ is
the perpendicular intersection point of $\Ax_\ga$ with $\wt \Delta$,
and where $m$ is the midpoint of $[x,\ga\cdot x]$, is a bijection.
  
\smallskip\noindent(ii) The map $\Phi^{1\,\ssm\, 2}$ from $\Perp'
(\Delta, \Delta, \frac{s}{2})$ to $\operatorname{AC}^{1\,\ssm\, 2}
(s)\ssm \operatorname{ARC} ^{1\,\ssm\, 2} (s)$ which maps the image in
$M_\ZZ$ of the common perpendicular between $\wt\Delta$ and a disjoint
image $\beta\cdot \wt\Delta$ with $\beta\in\Ga_\ZZ$ to the conjugacy
class of $\beta W\beta^{-1} W$ is a $2$-to-$1$ map.
\elemm

Similarly, $\operatorname{AC}^{2\,\ssm\, 1} (s)\ssm
\operatorname{ARC}^{2\,\ssm\, 1} (s)$ has half the cardinality of
$\Perp' (\Delta_1,\Delta_1,\frac{s}{2})$.

\medskip
\dem (i) If $\ga\in\Ga_\ZZ$ is primitive hyperbolic, ambiguous of the
first kind, with conjugacy class in $\operatorname{AC} ^{1\&2} (s)$,
let $x$ (which exists by Lemma \ref{lem:ambichar}
\hyperlink{ambicharii}{(ii)}) and $m$ be as in the statement. As seen
in the proof of the second claim of Lemma \ref{lem:ambichar}
\hyperlink{ambichariv}{(iv)}, note that $m$ is the perpendicular
intersection point of $\Ax_\ga$ with $\beta\cdot\wt\Delta_1$ for some
$\beta\in\Ga_\ZZ$. Hence $c_\ga=\Ga_\ZZ\cdot[x,m]$ is indeed a common
perpendicular between $\Delta$ and $\Delta_1$. By Lemma
\ref{lem:non12rec}, the interior of $[x,m]$ contains no point of the
orbit $\Ga_\ZZ\cdot i$. If the interior of $[x,m]$ was meeting
perpendicularly the image of $\Delta$ or of $\Delta_1$ by some
$\beta'\in\Ga_\ZZ$, then the element $\beta'W (\beta')^{-1} W$ or
$\beta'W_1(\beta')^{-1} W_1$, which belongs to $\Ga_\ZZ$ and is
hyperbolic, would have the same translation axis as $\ga$, and a
strictly shorter translation length, contradicting the fact that $\ga$
is primitive. Hence $c_\ga$ is primitive. Since $d(x,m)=\frac{1}{2}
d(x,\ga\cdot x)\leq \frac{s}{2}$, we have $c_\ga\in\Perp'(\Delta,
\Delta_1, \frac{s}{2})$.

If $\ga' \in\Ga_\ZZ$ is primitive hyperbolic, ambiguous of the first
kind, conjugated to $\ga$ and different from $\ga$, then by the second
claim of Lemma \ref{lem:ambichar} \hyperlink{ambichariv}{(iv)}, we
have $\ga'=\iota \ga \iota$. Hence the perpendicular intersection
point of $\Ax_{\ga'}$ with $\wt \Delta$ is $x'=\iota\cdot x$, and the
midpoint of $[x',\ga' \cdot x'] =\iota\cdot[x,\ga\cdot x]$ is $\iota
\cdot m$.  Thus we have the equality  $c_{\ga'}= \Ga_\ZZ\cdot[x',m'] = \Ga_\ZZ\,\iota
\cdot [x,m] =c_\ga$, and the map $\Psi^{1\&2}$ is well defined.  The
map $\Phi^{1\& 2}$ from $\Perp' (\Delta, \Delta_1, \frac{s}{2})$ to
$\operatorname{AC} ^{1\,\ssm\, 2} (s)$ which maps the image in $M_\ZZ$
of the common perpendicular between $\wt\Delta$ and a disjoint image
$\beta\cdot \wt\Delta_1$ for some $\beta\in\Ga_\ZZ$ to the conjugacy
class of $(\beta W_1\beta^{-1}) W$ is easily seen to be an inverse of
$\Psi^{1\&2}$.

\smallskip
(ii) Let $\beta\in\Ga_\ZZ$ be such that the intersection
$\wt\Delta\cap \beta\cdot\wt\Delta$ is empty. Let $\wt c=[x,y]$ be the
common perpendicular between $\wt\Delta$ and $\beta\cdot\wt\Delta$
with $x\in\wt \Delta$, and assume that its interior does not meet
perpendicularly an image of $\wt\Delta$, $\wt\Delta_1$ or $\{i\}$ by
an element of $\Ga_\ZZ$. Then the composition $\ga_{\wt c}=(\beta
W\beta^{-1}) W$ of the reflection $W$ in $\wt\Delta$ and the reflection
$\beta W\beta^{-1}$ in $\beta\cdot\wt\Delta$ is a hyperbolic element
of $\Ga_\ZZ$, with translation axis containing $\wt c$, and $y$ is the
midpoint of $x$ and $\ga_{\wt c}\cdot x$. Hence $d(x,\ga_{\wt c}\cdot
x)= 2\,d(x,y)\leq 2\,\frac{s}{2}=s$. If $\ga_{\wt c}$ is not
primitive, let $\ga_0$ and $k\geq 2$ be such that $\ga_{\wt c} =
\ga_0^k$. Then $\ga_0 \cdot x\in[x,y]$, and $\ga_0W$ would be a
reflection in $\Ga_\ZZ^+$ fixing the perpendicular bisector of $[x,\ga_0\cdot
  x]$. Therefore the interior of $[x,y]$ would meet perpendicularly an
image of $\wt\Delta$ or $\wt\Delta_1$ by an element of $\Ga_\ZZ$, a
contradiction.  Furthermore, $\ga_{\wt c}$ is ambiguous of the first
kind and neither reciprocal nor conjugated to an element ambiguous of the
second kind, by Lemma \ref{lem:ambichar} \hyperlink{ambicharii}{(ii)},
since $\Ax_{\ga_{\wt c}}$ meets perpendicularly $\wt \Delta$ at $x$,
and does not meet perpendicularly an image of $\{i\}$ or $\wt
\Delta_1$.  If $\alpha \in \Ga_\ZZ \ssm\{\id\}$ is such that
$\alpha\cdot \wt c= [x',y']$ is a common perpendicular between
$\wt\Delta$ and $\beta' \cdot\wt\Delta$ for some $\beta'\in\Ga_\ZZ$,
then $\alpha$ preserves $\wt \Delta$, hence $\alpha=\iota$, hence
$x'=\iota\cdot x$, $y'=\iota\cdot y$, and $\ga_{\alpha\,\wt c}
=\iota\,\ga_{\wt c}\,\iota $ is conjugated to $\ga_{\wt c}$. Hence the
map $\Phi^{1\,\ssm\, 2}$ is well defined.

Let us prove that the map $\Phi^{1\,\ssm\, 2}$ is onto. Let $\ga\in
\Ga_\ZZ$ be primitive hyperbolic, ambiguous of the first kind and neither 
reciprocal nor conjugated to an element ambiguous of the second
kind. Let $x$ be the perpendicular intersection point of $\wt \Delta$
and $\Ax_\ga$, and $m$ the midpoint of $[x,\ga\cdot x]$. Then $\wt c=
[x,m]$ is a common perpendicular between $\wt\Delta$ and $\beta \cdot
\wt\Delta$ for some $\beta\in\Ga_\ZZ$, that does not meet
perpendicularly in its interior an image of $\wt\Delta$ or
$\wt\Delta_1$ or $\{i\}$ by an element of $\Ga_\ZZ$, since $\ga$ is
primitive. By construction, we have $\ga=\ga_{\wt c}$. By the proof of
the first claim of Lemma \ref{lem:ambichar}
\hyperlink{ambichariv}{(iv)}, the preimages by $\Phi^{1\,\ssm\, 2}$ of
the conjugacy class of $\ga$ are the images of $[x,m]$ and
$[m,\ga\cdot x]$ in $M_\ZZ$. Note that the two segments $[x,m]$ and
$[m,\ga\cdot x]$ are not in the same orbit under $\Ga_\ZZ$, since
$\ga$ is primitive.  This proves that $\Phi^{1\,\ssm\, 2}$ is
$2$-to-$1$. \cqfd

\blemm\label{lem:1rec} The map $\Phi^{1R}$ from $\Perp'(\Delta,
\{\Ga_\ZZ \cdot i\},\frac{s}{4})$ to $\operatorname{ARC}^{1\,\ssm\, 2}
(s)$ which maps the image in $M_\ZZ$ of the common perpendicular
between $\wt\Delta$ and $\alpha\cdot \{i\}$, for some
$\alpha\in\Ga_\ZZ$, to the conjugacy class of $(\alpha
\iota\alpha^{-1} W)^2$ is a $2$-to-$1$ map.
\elemm

Similarly, the set $\operatorname{ARC}^{2\,\ssm\, 1} (s)$ has half the
cardinality of $\Perp' (\Delta_1,\{\Ga_\ZZ \cdot i\},\frac{s}{4})$.

\begin{center}
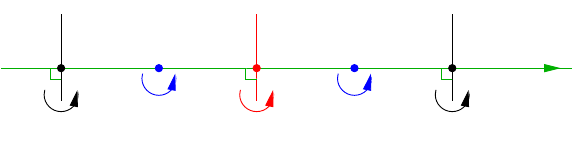
\end{center}
\vspace*{-0.8cm}\addtocounter{fig}{1}
\begin{center}
{\small Figure \arabic{fig}~: \hypertarget{monoambigcommperp}{Relating
    common perpendiculars and monoambiguous conjugacy classes}.}
\end{center}
\dem Every element $\Ga_\ZZ \cdot c$ in $\Perp'(\Delta, \{\Ga_\ZZ
\cdot i\}, \frac{s}{4})$ is the image in $M_\ZZ$ of the common
perpendicular $[x,z]$ between $\wt\Delta$ and $\alpha\cdot \{i\}$, for
some 
$\alpha\in\Ga_\ZZ$ with $z=\alpha\cdot i$. Since $z\notin
\wt\Delta$, the element $\ga=(\alpha \iota\alpha^{-1} W)^2$, which
belongs to $\Ga_\ZZ$, is 
hyperbolic 
with translation axis
containing $[x,z]$ and translation length $4\,d(x,z)\leq
4\,\frac{s}{4}=s$. Since $\Ax_\ga$ meets perpendicularly $\wt \Delta$
and contains $\alpha\cdot i$, the element $\ga$ is reciprocal and
ambiguous of the first kind by Lemma \ref{lem:ambichar}
\hyperlink{ambicharii}{(ii)}. As in the previous proof, $\ga$ is
primitive and the conjugacy class of $\ga$ in $\Ga_\ZZ$ does not
depend on the above choice of a representative $[x,z]$ of $\Ga_\ZZ
\cdot c$.

Let us prove that $\Phi^{1R}$ is onto. Let $\ga\in\Ga_\ZZ$ be
primitive hyperbolic, reciprocal and ambiguous of the first kind
(hence not conjugated to an element ambiguous of the second kind by
Lemma \ref{lem:non12rec}). As in the proof of the second claim of
Lemma \ref{lem:1non2}, let $x$ be the perpendicular intersection point
of $\wt \Delta$ and $\Ax_\ga$, let $m$ be the midpoint of $[x,\ga\cdot
  x]$ and let $\beta\in\Ga_\ZZ$ be such that $\beta\cdot\wt\Delta$ is
the perpendicular bisector of $[x,\ga\cdot x]$. By Lemma \ref{lem:ambichar}
\hyperlink{ambicharii}{(ii)}, since $\ga$ is reciprocal, the
translation axis $\Ax_\ga$ meets the orbit $\Ga_\ZZ\cdot i$. By
translating by powers of $\ga$, there is an orbit point $z=\alpha\cdot
i$ in $[x,\ga\cdot x[\,$. Since $\Ga_\ZZ^+$ preserves $\Ga_\ZZ\cdot i$,
up to replacing $z$ by its image by the reflection $\beta W\beta^{-1}$,
we may assume that $z\in [x,m]$.  Since $\Ax_\ga$ has irrational
endpoints, we have $z\neq x,m$. If $z$ is not the midpoint of $[x,m]$,
say $d(x,z)< \frac{1}{2}d(x,m)$, then $(\alpha \iota\alpha^{-1}W)^2$
would be an hyperbolic element in $\Ga_\ZZ$ with translation axis
$\Ax_\ga$ and translation length $4\,d(x,z)< 2\,d(x,m)=d(x,\ga \cdot
x)$, contradicting the fact that $\ga$ is primitive. Hence $d(x,z)=
\frac{1}{4}\,d(x,\ga \cdot x)\leq \frac{s}{4}$. Since $\ga$ is
primitive, the interior of $[x,z]$, which is a common perpendicular
between $\wt\Delta$ and $\alpha\cdot \{i\}$, does not meet
perpendicularly an image of $\wt\Delta,\wt\Delta_1$ or $\{i\}$.  Hence
$\Phi^{1R}$ is onto. Furthermore, using the reflection $\beta
W\beta^{-1}$, there also exists $\alpha'\in \Ga_\ZZ$ such that the
midpoint $z'$ of $[m, \ga\cdot x]$ is $z'=\alpha'\cdot i$. By the
proof of the first claim of Lemma \ref{lem:ambichar}
\hyperlink{ambichariv}{(iv)}, the preimages by $\Phi^{1R}$ of the
conjugacy class of $\ga$ are the images of $[x,z]$ and $[m,z']$ in
$M_\ZZ$. Note that these two segments are not in the same orbit under
$\Ga_\ZZ$, since $\ga$ is primitive. This proves that $\Phi^{1R}$ is
$2$-to-$1$.
\cqfd

\medskip
As $\Delta$ and $\Delta_1$ are strongly reversible, we have $\iota_{\rm rev}
(\Delta)=\iota_{\rm rev} (\Delta_1)=1$.  Note that the order of the
pointwise stabilizer in $\Ga_\ZZ$ of $\wt\Delta$ and $\wt\Delta_1$ is
$m(\Delta)=m(\Delta_1)=1$, and the one of $\{i\}$ is $m(\{\Ga_\ZZ
\cdot i\})=2$.  Furthermore, by the normalisation of the
Patterson-Sullivan measures in Section \ref{sec:geomback}, we have
$\|\sigma^+_{\{\Ga_\ZZ \cdot i\}}\| =\frac{1}{m(\{\Ga_\ZZ \cdot i\})}
\|\mu_i\|= \frac{1}{2} \Vol(\SSS^1) =\pi$. Recall that $\Vol(M_\ZZ)=
\frac{\pi}{3}$. By the standard argument comparing the growth of
primitive closed geodesics and the nonprimitive ones, see for instance
Step 2 of the proof of \cite[Theorem 9.11]{PauPolSha15}, as
$s\ra+\infty$, if $D^-,D^+\in \{\Delta,\Delta_1\}$, by Theorem
\ref{theo:divergentperp2} applied with $n=2$ and $M=M_\ZZ$, we have
then
\begin{equation}\label{eq:growthprimesansi}
  \card\Perp'(D^-,D^+,s)=\frac{3}{2\pi^2}\;
  s^2\, e^{s}+ \bigO\big(s\;e^{ s}\big)\,.
\end{equation}
Similarly, by Theorem \ref{theo:1geoddivreal}, if $D^+\in \{\Delta,
\Delta_1\}$, then
\begin{equation}\label{eq:growthprimesaveci}
  \card\Perp'(\{\Ga_\ZZ \cdot i\},D^+,s)=
  \frac{3} {2\pi}\; s\, e^{s} + \bigO(e^{s})\,,
\end{equation}

Now, by Lemmas \ref{lem:non12rec} and \ref{lem:1rec}, and
Equation \eqref{eq:growthprimesaveci}, we have
\begin{align*}
  &\card\operatorname{ARC}(s)=\card\operatorname{ARC}^{1\,\ssm\, 2}(s)
  + \card\operatorname{ARC}^{2\,\ssm\, 1}(s)
  \\=\;& \frac{1}{2}\card\Perp'(\Delta,\{\Ga_\ZZ \cdot i\},\frac{s}{4})
  +\frac{1}{2}\card\Perp'(\Delta_1,\{\Ga_\ZZ \cdot i\},\frac{s}{4})
  =\frac{3}{8\,\pi}\;s\,e^{\frac s4}+\bigO(e^{\frac s4})\,.
\end{align*}
Similarly, using Lemmas \ref{lem:non12rec} and \ref{lem:1non2}, the
previous computation that implies that we have $\card
\operatorname{ARC}(s)= \bigO(s\;e^{\frac s4})$, and Equation
\eqref{eq:growthprimesansi}, we have
\begin{align*}
  \card\operatorname{AC}(s)&=\card\big(\operatorname{AC}^{1\,\ssm\, 2}(s)
  \ssm\operatorname{ARC}^{1\,\ssm\, 2}(s)\big)
  + \card\big(\operatorname{AC}^{2\,\ssm\, 1}(s)\ssm
  \operatorname{ARC}^{2\,\ssm\, 1}(s)\big)\\&\quad
  + \card\operatorname{AC}^{1\&2}(s)+\card\operatorname{ARC}(s)
\\&=\frac{1}{2}\card\Perp'(\Delta,\Delta,\frac{s}{2})
  +\frac{1}{2}\card\Perp'(\Delta_1,\Delta_1,\frac{s}{2})
  \\&\quad+\card\Perp'(\Delta,\Delta_1,\frac{s}{2})+
  \bigO(s\,e^{\frac s4})\\ &=
  2\Big(\frac{3}{2\pi^2}\;\big(\frac{s}{2}\big)^2\, e^{\frac{s}{2}}+
  \bigO\big(\frac{s}{2}\;e^{\frac{s}{2}}\big)\Big)=
  \frac{3}{4\,\pi^2}\;s^2\,e^{\frac s2}+\bigO(s\;e^{\frac s2})\,.
\end{align*}
This concludes the proof of Theorem \ref{theo:ambi}.
\cqfd

\begin{center}
\begin{overpic}[width=.95\textwidth]{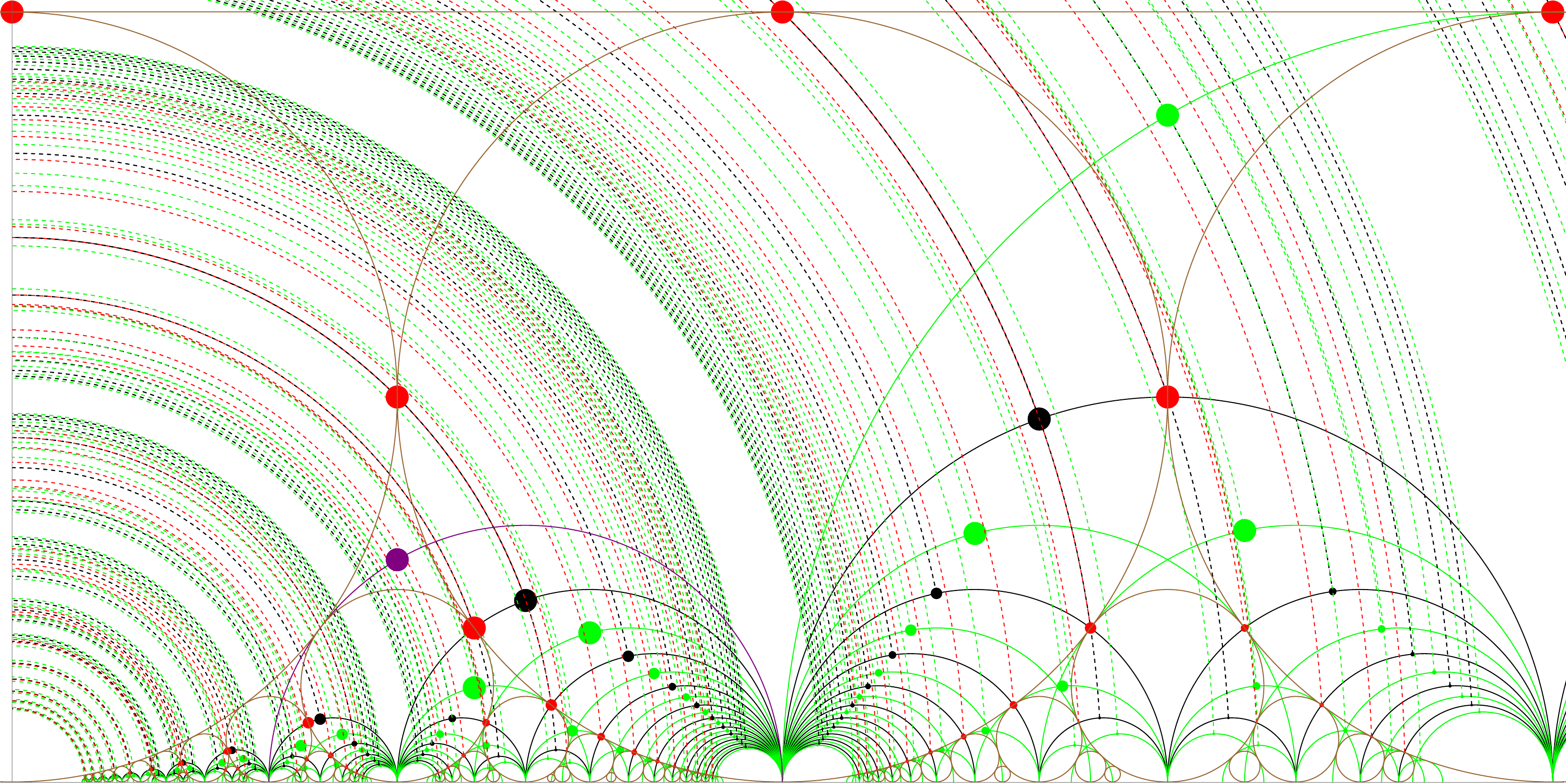} 
  \put (0,-3) {$0$}
  \put (16.25,-3) {$\frac{1}{3}$}
  \put (24.5,-3) {$\frac{1}{2}$}
  \put (49.5,-3) {$1$} 
  \put (99,-3) {$2$}
  \put (73.75,-3) {$\frac{3}{2}$}
      \put (-2,49) {$i$}
\end{overpic}
\end{center}
\vspace*{-0.2cm}\addtocounter{fig}{1}
\begin{center}
{\small Figure \arabic{fig}~: \hypertarget{creosote}{An anthology of
    common perpendiculars}.}
\end{center}

Figure \hyperlink{creosote}{\arabic{fig}} shows images in the right
halfplane ${\hdr}^+$ by elements of $\Ga_\ZZ$ of the positive
imaginary axis $\wt\Delta$ in continuous black, and images of
$\wt\Delta_1$ in continuous green (except the geodesic line
$\interval[open]{\frac 13}1=\nu\wt \Delta_1$ used in the proof of
Lemma \ref{lem:ambichar} \hyperlink{ambichariii}{(iii)} which is drawn
in purple).  Images by elements of $\Ga_\ZZ$ of the horosphere
$\partial\H_\infty$ are drawn in brown. Images by elements of
$\Ga_\ZZ$ of $i$ are drawn as red points.  The common perpendiculars
starting from $\wt \Delta$ and ending at images of $\wt\Delta$, $\wt
\Delta_1$ or $\{i\}$ by elements of $\Ga_\ZZ$ are drawn with dashed
lines in black, green and red correspondingly to the color of their
arrival point, that are marked by black, green and red dots
respectively (except the purple one on $\interval[open]{\frac
  13}1$). Note that there are (nonprimitive) common perpendiculars
passing through black and green points (again giving examples for
Lemma \ref{lem:ambichar} \hyperlink{ambichariii}{(iii)}), or through
black and red points (corresponding to elements of
$\operatorname{ARC}^{1\,\ssm \,2}(s)$ for some large enough $s$, with
the notation of the proof of Theorem \ref{theo:ambi}), but none
through black and green and red dots (in accordance with Lemma
\ref{lem:non12rec}).

\section{On the binary additive divisor problem for integers}
\label{sec:ingham}

In this section, we discuss the connection of Theorem
\ref{theo:divergentperp2} with the binary additive divisor problem in
$\ZZ$ and use this connection to show that the error term obtained in
Theorem \ref{theo:divergentperp2} is optimal.

Let $\nod:\NN\ssm\{0\}\ra \NN\ssm\{0\}$ be the {\em number of divisors
  function} of the natural numbers, defined by $n\mapsto \card\{d\in
\NN\ssm\{0\}:d\mid n\}= \frac{1}{2}\,\card\{d\in \ZZ\ssm\{0\}:d\mid
n\}$. The {\em binary additive divisor problem}\footnote{See for
instance \cite{Ingham27,Estermann31,HeathBrown79, Motohashi94}.} in
$\ZZ$ studies the asymptotic properties as $n\ra+\infty$ of the sums
$\sum_{k=1}^n\nod(k)\,\nod(k+f)$ for any positive integer $f$.  The
link between our counting problem of common perpendiculars between
divergent geodesics and the binary additive divisor problem in $\ZZ$
will be given by Proposition \ref{prop:relatcommperpbinadddiv} below.

We keep the notation $\Ga_\ZZ$, $M_\ZZ$, $\wt \ell$, $\wt \Delta$,
$\ell$, $\Delta$ and ${\hdr}^\pm$ of Section \ref{sec:ambi}. We start
by proving a quantitative complement to Lemma
\ref{lem:travaildivgeodmodular}.

\blemm \label{lem:travaildivgeodmodularbis} Let
$\ga=\begin{bsmallmatrix} a&b\\c&d\end{bsmallmatrix}
\in\Ga_\ZZ$ with $a,b,c,d>0$. Then the length of the common
perpendicular between $\wt\Delta$ and $\ga\wt\Delta$ is $\lambda=
\arcosh(1+2bc)$.
\elemm

Figure
\addtocounter{fig}{1}\hyperlink{autocommperpDelta}{\arabic{fig}}\addtocounter{fig}{-1}
shows in black some of the $\Ga_\ZZ$-translates of $\wt\Delta$ in the
right halfplane ${\hdr}^+$, and the six closest points to the vertical
geodesic $\wt\Delta$ on the $\Ga_\ZZ$-translates at distance
$\arcosh(9)$, corresponding to $bc=4$ with the above notation, and in
green the corresponding common perpendiculars.

\begin{center}
\begin{overpic}[width=0.95\textwidth,tics=10]{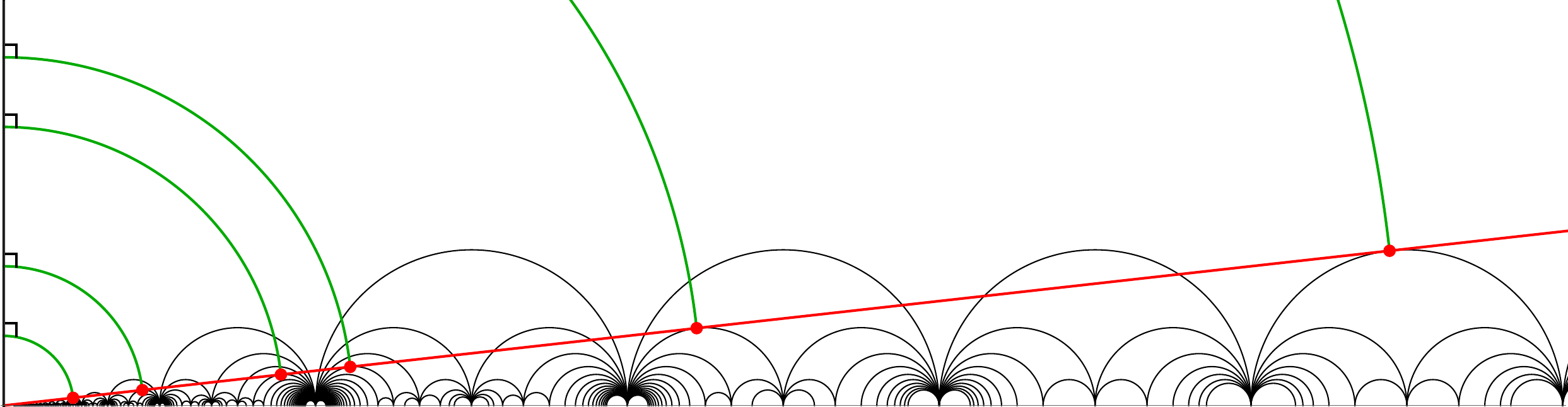}
  \put (-1,-2.5) {$0$}
 \put (19.5,-2.5) {$1$}
  \put (39.25,-2.5) {$2$}
  \put (59,-2.5) {$3$}
    \put (79,-2.5) {$4$}
     \put (99,-2.5) {$5$}
     \put (-2,10) {$\frac{1}{2}$}
     \put (-2,20) {$1$}
\end{overpic}
\end{center}
\vspace*{-0.4cm}\addtocounter{fig}{1}
\begin{center}
{\small Figure \arabic{fig}~: \hypertarget{autocommperpDelta}{Common
    perpendiculars between $\wt\Delta$ and its translates}.}
\end{center}

\dem Let $x=\frac{b}{d}$ and $y=\frac{a}{c}$, which are the two
endpoints at infinity of $\ga\cdot\wt\Delta$. Since $a,b,c,d>0$ and
$ad-bc=1$, we have $ad>bc$, hence $0<x<y$. The common perpendicular
between $\wt\Delta$ and $\ga\cdot\wt\Delta$ is a segment of the
Euclidean halfcircle centered at $0$ that intersects $\ga\cdot
\wt\Delta$ at a right angle, see Figure
\addtocounter{fig}{1}\hyperlink{computcommperp}{\arabic{fig}}\addtocounter{fig}{-1}.

\begin{center}
  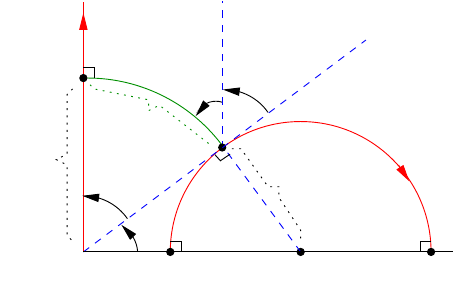
\end{center}
\vspace*{-0.8cm}\addtocounter{fig}{1}
\begin{center}
{\small Figure \arabic{fig}~: \hypertarget{computcommperp}{Computation of lengths of common
    perpendiculars}.}
\end{center}
The intersection point is the unique point on $\ga\cdot\wt\Delta$
where a Euclidean line $L_{x,y}$ through the origin is tangent to
$\ga\cdot\wt\Delta$. If $\phi\in\interval[open]0\pi$ is the angle that
$L_{x,y}$ makes with the positive real line at the origin, then
$\sin\phi=\frac{\frac{y-x}2}{\frac{y+x}2} =\frac{y-x}{y+x}$.  By the
angle of parallelism formula\footnote{See Equation
\eqref{eq:complexlengthbianchi} in Section \ref{sec:motohashi} for a
different computation using complex length.
} already used in Section \ref{sec:realhypgeom} and since
$ad-bc=1$, the length $\lambda$ of the common perpendicular from
$\wt\Delta$ to $\ga\cdot\wt\Delta$ is
\[
\lambda= \arcosh\big(\frac 1{\sin\phi}\big)=
\arcosh\big(\frac{y+x}{y-x}\big)=\arcosh
\Big(\frac{\frac{a}{c}+\frac{b}{d}}{\frac{a}{c}-\frac{b}{d}}\Big)
=\arcosh(1+2bc)\;.\quad\Box
\]

\bprop\label{prop:relatcommperpbinadddiv}
For every $s>0$, we have $\displaystyle \N_{\Delta,\Delta}(s)=
\sum_{k=1}^{\lfloor \frac 12(\cosh s-1)\rfloor}\nod(k)\,\nod(k+1)$.
\eprop

\dem By the comment just before Lemma \ref{lem:travaildivgeodmodular},
since the set of $\Ga_\ZZ$-translates of $\wt \Delta$ is in bijection
with the set of right cosets $\Ga_\ZZ/\Ga_{\wt \Delta}$, it follows
from Lemmas \ref{lem:travaildivgeodmodular} and
\ref{lem:travaildivgeodmodularbis} that $\N_{\Delta,\Delta}(s)$ is the
number of quadruples $(a,b,c,d)$ of positive integers such that
$ad-bc=1$ and $\arcosh(1+2bc)\leq s$, or equivalently $bc\leq n=
\big\lfloor\frac{1}{2}(\cosh s-1)\big\rfloor$ since $bc$ is a positive
integer. This number is $\sum_{k=1}^{n}\nod(k)\,\nod(k+1)$.
\cqfd

\medskip
Let us now relate Theorem \ref{theo:divergentperp2} (in the special
case $n=2$ and $\Ga_\ZZ=\PSL_2(\ZZ)$) with the known asymptotic
result on the binary additive divisor problem in $\ZZ$.  After the
work of Ingham \cite[p.~205]{Ingham27}, a major input by Estermann
\cite{Estermann31}, and various improvements on the error term by for
instance \cite[Thm.~2]{HeathBrown79} and \cite[Coro.~1]{Motohashi94},
we now know that there exist $a_1,a_2\in\RR$ such that (with a
simplified version of the best known error term)
\begin{equation}\label{eq:badpz}
  \sum_{k=1}^n\nod(k)\,\nod(k+1)=
  \frac 6{\pi^2}\, n(\ln n)^2+a_1\,n\ln n+a_2\,n+
\bigO(n^{\frac{5}{6}})\,.
\end{equation}
Using \cite[Eq.~(36)]{Estermann31}, we can compute the estimate
$a_1\simeq 1.574>0$.  By Proposition \ref{prop:relatcommperpbinadddiv}
and since $\lfloor \frac 12(\cosh s-1)\rfloor=\frac{1}{4} e^s+
\bigO(1)$, we thus have
\begin{equation}\label{eq:inghamgeom}
  \N_{\Delta,\Delta}(s)= 
  \sum_{k=1}^{\lfloor \frac 12(\cosh s-1)\rfloor}\nod(k)\,\nod(k+1)=
  \frac {3}{2\,\pi^2}\, s^2e^s+b_1\, s\, e^s +b_2\, e^s+
  \bigO(e^{\frac{5}{6}s})\,,
\end{equation}
where $b_1=\frac{a_1}{4}-\frac{6\ln 2}{\pi^2}\simeq -0,028\ne 0$.
Equation \eqref{eq:inghamgeom} agrees with the asymptotic
\[
\N_{\Delta,\Delta}(s)=
\frac{3}{2\pi^2}\,s^2e^{ s}+ \bigO(s\,e^{ s})
\]
given by Theorem \ref{theo:divergentperp2}, as seen for Equation
\eqref{eq:growthprimesansi}.  Furthermore, Equation
\eqref{eq:inghamgeom} gives an explicit nonzero term of the order $s\,
e^s$ and an error term of strictly smaller order.  This shows that the
size of the error term in Theorem \ref{theo:divergentperp2} is
optimal.

\section{The binary additive divisor problem for
  imaginary quadratic integers}
\label{sec:motohashi}

In this section, we use our asymptotic counting of common
perpendiculars between divergent geodesics proven in Theorem
\ref{theo:divergentperp2} in order to study the asymptotic binary
additive divisor problem for imaginary quadratic integers, confirming
a particular case of a conjecture of Motohashi
\cite[p.~277]{Motohashi01}.

Let $K, D_K,\OOO_K,\zeta_K$ and $\Naivenod_K$ be as in the
introduction.  Recall that the order $|\OOO_K^\times|$ of the group of
units $\OOO_K^\times$ of $\OOO_K$ is equal to $4$ when $D_K=-4$, to
$6$ when $D_K=-3$, and to $2$ otherwise.

\medskip\noindent {\bf Proof of Theorem \ref{theo:mainarithintro}. }
We start the proof by describing the relevant geometric framework.  As
in Section \ref{sec:realhypgeom}, 
let $\htr\subset
\CC\times\RR$ be the upper halfspace model of the real hyperbolic
$3$-space $\htr$. Its boundary at infinity is $\partial_\infty\htr=
(\CC\times\{0\})\cup \{\infty\}$, that we identify with
$\PP^1(\CC)=\CC\cup\{\infty\}$. The group $\PSLC$ acts isometrically
and faithfully on $\htr$, by the Poincaré extensions of the (complex)
homographies.  The {\em Bianchi group} $\Ga_{\!\OOO\!_K}=\PSLOK$ is a
nonuniform arithmetic lattice in $\PSLC$, with set of parabolic fixed
points $\operatorname{Par} _{\Ga_{\!\OOO\!_K}} =
\PP^1(K)=K\cup\{\infty\}$.  The {\em Bianchi orbifold}
$M_{\OOO\!_K}=\Ga_{\!\OOO\!_K}\bs\htr$ is a noncompact finite volume
complete connected real hyperbolic good orbifold. Its number of cusps
is equal to the class number of $K$, see for instance \cite[\S
  7.2]{ElsGruMen98}.

Let $\wt\ell:t\mapsto (0,e^t)$ be the geodesic line in $\htr$ through
$(0,1)\in\htr$ at time $t=0$, with endpoints at infinity $ 0$ and
$\infty$. Then $\ell=\Ga_{\!\OOO\!_K}\wt \ell$ is a divergent geodesic
in $M_{\OOO\!_K}$, converging at $\pm\infty$ to the cusp
$\Ga_{\!\OOO\!_K}\cdot\infty$ of $M_{\OOO\!_K}$. Note that $\ell$ is
strongly reversible since the image $\wt \Delta =\wt\ell(\RR)$ of $\wt \ell$ is
preserved by the involution $\iota= \begin{bsmallmatrix} 0&-1
    \\1&\ \; 0\end{bsmallmatrix} \in \Ga_{\!\OOO\!_K}$, whose fixed
point set is the geodesic line with points at infinity $-i$ and $i$,
that meets  $\wt \Delta$ perpendicularly at $(0,1)$.  The pointwise
stabilizer of $\wt \Delta$ in $\Ga_{\!\OOO\!_K}$ is the group
consisting of the diagonal elements $\begin{bsmallmatrix} u&0\;\;
\\0&\overline u\end{bsmallmatrix}$ for $u\in\OOO_K^\times$, and the
(global) stabilizer $\Ga_{\wt \Delta}$ of $\wt \Delta$ in
$\Ga_{\!\OOO\!_K}$ is the binary dihedral group generated by $\iota$
and the pointwise stabilizer.  Hence, with $\Delta= \ell(\RR)=
\Ga_{\!\OOO\!_K}\wt \Delta$ the image of $\ell$, we have
\begin{equation}\label{eq:stabfixDelta}
  m(\Delta)=\frac{|\OOO_K^\times|}{2}\quad\text{and}\quad
  |\Ga_{\wt \Delta}|=|\OOO_K^\times|\,.
\end{equation}
In particular, $m(\Delta)=1$ unless $D_K=-3$ or $D_K=-4$.

For every $k\in \OOO_K\ssm\{0,1\}$, the product $\nod_K(k)\,
\nod_K(k-1)$ is the number of representations of $1$ as the difference
$ad-bc$ for a quadruple $(a,b,c,d)$ of elements of $\OOO_K\ssm\{0\}$
such that $ad=k$ and $bc=k-1$.  Hence
\begin{equation}\label{eq:ddKtoSLO}
  \nod_K(k)\,\nod_K(k-1)=\card
  \big\{\ga=\begin{psmallmatrix}a&b\\c&d\end{psmallmatrix}
    \in\SLOK:a\,b\,c\,d\neq 0,\;ad=k\big\}\,.
\end{equation}

Let $\ga=\begin{bsmallmatrix} a&b\\c&d \end{bsmallmatrix}
\in\Ga_{\!\OOO\!_K}$. The geodesic lines $\wt \Delta$ and $\ga\wt
\Delta$ have no common endpoint at infinity if and only if $a\,b\,c\,d
\ne 0$, by the same argument as the one at the beginning of the proof
of Lemma \ref{lem:travaildivgeodmodular}.  As in that proof,
considering now the horoball $\H_\infty= \{(z,v)\in\htr:v\geq 1\}$ and
replacing $i\in\hdr$ by $(0,1)\in\htr$, we see that if $a\,b\,c\,d \ne
0$, then the geodesic lines $\wt \Delta$ and $\ga\wt \Delta$ have
empty intersection, since the $\Ga$-equivariant family
$(\ga\H_\infty)_{\ga\in\Ga/\Ga_\infty}$ is again precisely invariant,
the only horoballs in this family containing $(0,1)$ are $\H_\infty$
and $\iota\H_\infty$ and $\wt \Delta$ is contained in $\H_\infty\cup
\iota\H_\infty$. In particular, $\wt \Delta$ and $\ga\wt \Delta$ have
a common perpendicular if $a\,b\,c\,d \ne 0$.

Let $\lambda_\ga= d(\wt \Delta,\ga\wt \Delta)=d(p,q)>0$ be the length
of the common perpendicular $[p,q]$ between $\wt \Delta$ and $\ga\wt
\Delta$, with $p\in \wt\Delta$. Let $\theta_\ga$ be the angle at $p$
between the parallel transport of the oriented geodesic line $\ga
\wt\Delta$ along $[p,q]$ and the oriented geodesic line $\wt\Delta$.
By \cite[Lemma 2.2]{ParPau11MZ} and since $ad-bc=1$, we have
\begin{equation}\label{eq:complexlengthbianchi}
  \cosh\lambda_\ga+\cos\theta_\ga=
  2\,\frac{|\gamma\cdot\infty|}{|\gamma\cdot\infty-\gamma\cdot 0|}
  =2\,|ad|\,.
\end{equation}
For all $\lambda\geq 0$ and $\theta\in\RR$, we have
$|\frac{e^{-\lambda}}2 + \cos\theta|\le 2$, so that as $\lambda\ra
+\infty$, a first order approximation gives
\begin{equation}\label{eq:approx}
  \ln(\cosh\lambda+\cos\theta)=
  \ln\Big(\frac{e^\lambda}2+\frac{e^{-\lambda}}2+\cos\theta\Big)
  =\lambda-\ln 2+\bigO(e^{-\lambda})\,.
\end{equation}

The only unit normal vectors to $\wt \Delta$ that have a nontrivial
stabilizer in $\Ga_{\!\OOO\!_K}$ are the finitely many tangent vectors
$v_u$ at $(0,1)$ of the oriented geodesic lines with points at
infinity $u$ and $-u$ for $u\in\OOO_K^\times$, stabilized by the
element $\begin{bsmallmatrix} \ 0&u\\-\overline u&0
\end{bsmallmatrix}$. The number of common perpendiculars between
$\wt \Delta$ and its images by the elements of $\Ga_{\!\OOO\!_K}$,
whose initial tangent vectors are fixed to be $v_u$, grows at most
linearly in their length. In particular, all multiplicities
$m_{\wt\Delta,\ga\wt\Delta}$ of these common perpendiculars are equal
to $1$, except for a number of them that is linear in their length.  A
formula due to Humbert for the volume of
$M_{\OOO\!_K}=\Ga_{\!\OOO\!_K}\bs\htr$ gives
\[
\Vol(M_{\OOO\!_K})=
\frac{1}{4\pi^2}\;|D_K|^{3/2}\,\zeta_K(2)\;,
\]
see for instance Sections 8.8, 9.6 of \cite{ElsGruMen98}.  With $n=3$,
let us define
\begin{equation}\label{eq:cstcKdef}
  c_K=\frac{(n-1)\;\pi^{\frac{n}{2}-1}\;\Ga(\frac{n}{2})\,
    \iota_{\rm rev}(\Delta)^2}
{2^{n+1}\;\Ga(\frac{n+1}{2})^2\;m(\Delta)^2\;\Vol M_{\OOO\!_K}}
=\frac{\pi^3}{|\OOO_K^\times|^2\,|D_K|^{\frac 32}\, \zeta_K(2)}\,.
\end{equation}
By Theorem \ref{theo:divergentperp2} applied with $n=3$,
$M=M_{\OOO\!_K}$ and $D^-=D^+=\Delta$, we have
\begin{align}
\N_{\Delta,\Delta}(s)&=\card\big\{[\ga]\in\Ga_{\wt \Delta}\bs\Ga_{\!\OOO\!_K}/
\Ga_{\wt \Delta} : 0<d(\wt \Delta,\ga\,\wt \Delta)\leq s\big\}+\bigO(s)
\nonumber \\&=c_K\;s^2e^{2 s}+ \bigO(s\,e^{2 s})\,.
\label{eq:bianchicountdownstairs}
\end{align}

Using, in the following computations, respectively

\noindent$\bullet$~ Equations \eqref{eq:ddKtoSLO} and
\eqref{eq:complexlengthbianchi} for the first equality,

\noindent$\bullet$~ the facts that the kernel of the isometric action
of $\SL_2 (\OOO_K)$ on $\htr$ is the subgroup $\{\pm \id\}$ with order
$2$, that $\Ga_{\!\OOO\!_K}=\SL_2(\OOO_K) /\{\pm \id\}$, that the
assumptions on $\ga$ in the second line depend only on its double class
$[\ga]$ in $\Ga_{\wt \Delta}\bs\Ga_{\!\OOO\!_K}/\Ga_{\wt \Delta}$ and
that $|\Ga_{\wt \Delta}|=|\OOO_K^\times|$ by Equation
\eqref{eq:stabfixDelta} for the second equality,

\noindent$\bullet$~ Equation \eqref{eq:approx} for the third equality,

\noindent$\bullet$~ a partition of the set of $[\ga]\in\Ga_{\wt \Delta}
\bs\Ga_{\!\OOO\!_K}/\Ga_{\wt \Delta}$ with $\lambda_\ga+
\bigO(e^{-\lambda_\ga}) \le \ln(4N)$ into on the one hand the ones
with $\ln(4\, \sqrt{N}\,)<\lambda_\ga +\bigO(e^{-\lambda_\ga}) \le
\ln(4N)$, so that $\lambda_\ga \geq \ln(4\,\sqrt{N}\,)+\bigO(1)$ hence
$e^{-\lambda_\ga} =\bigO(N^{-\frac 12})$ thus by bootstrap we obtain  
$\lambda_\ga\leq \ln(4N)+ \bigO(N^{-\frac 12}) =\ln(4N+
\bigO(\sqrt{N}))$, and on the other hand the ones with
$\lambda_\ga+\bigO(e^{-\lambda_\ga})\leq \ln(4\,\sqrt{N}\,)$ so that
$\lambda_\ga\leq\ln(4\sqrt{N}\,)+\bigO(1)$, for the fourth equality,

\noindent$\bullet$~ Equation \eqref{eq:bianchicountdownstairs} for the
fifth equality and Equation \eqref{eq:cstcKdef} for the last one,

\noindent we have
\begin{align} \label{eq:Motohashigeneral}
  & \frac{1}{|\OOO_K^\times|^2} \sum_{\substack{k\in\OOO_K\sm\{0,1\}\\|k|\le N}}
  \nod_K(k)\,\nod_K(k-1)
  \nonumber \\
  =\;&\frac{1}{|\OOO_K^\times|^2}
\card\Big\{\ga\in\SLOK:\begin{array}{l}
\partial_\infty(\ga\wt \Delta)\cap\partial_\infty\wt \Delta=\emptyset
\\\cosh\lambda_\ga+\cos\theta_\ga\le 2N\end{array}\Big\}
\nonumber \\
=\;&2\card\Big\{[\ga]\in\Ga_{\wt \Delta}\bs\Ga_{\!\OOO\!_K}/\Ga_{\wt \Delta}:
\begin{array}{l}
\partial_\infty(\ga\wt \Delta)\cap\partial_\infty\wt \Delta=\emptyset
\\\cosh\lambda_\ga+\cos\theta_\ga\le 2N\end{array}\Big\}
\nonumber \\
=\;&2\card\Big\{[\ga]\in\Ga_{\wt \Delta}\bs\Ga_{\!\OOO\!_K}/\Ga_{\wt \Delta}:
\begin{array}{l}
\partial_\infty(\ga\wt \Delta)\cap\partial_\infty\wt \Delta=\emptyset
\\\lambda_\ga+\bigO(e^{-\lambda_\ga}) \le \ln(4N)\end{array}\Big\}
\nonumber \\
=\;&2\,\N_{\Delta,\Delta}\big(\ln(4N+\bigO(\sqrt{N}\,))\big)+
\bigO\big(\N_{\Delta,\Delta}(\ln(4\,\sqrt{N}\,)+\bigO(1))\big)
\nonumber \\
=\;&2\,c_K\,\ln^2(4N+\bigO(\sqrt{N}\,))(4N+\bigO(\sqrt{N}\,))^2
\nonumber \\
&\!\!\!\!\!\!+\!\bigO\!\big(\ln(4N\!+\!\bigO(\sqrt{N}\,))(4N\!+\!\bigO(\sqrt{N}\,))^2\big)\!
+\!\bigO\!\big((\ln(4\sqrt{N})+\bigO(1))^2(4\sqrt{N}\,)^2\big)
\nonumber \\
=\;&\frac{32\pi^3}{|\OOO_K^\times|^2\,|D_K|^{\frac 32}\,
  \zeta_K(2)}(\ln N)^2N^2+\bigO((\ln N)N^2)\;.
\end{align}
As the sums $\sum_{\substack{k\in\OOO_K\ssm\{0,-1\}\\|k|\le N}}
\nod_K(k) \,\nod_K(k+1)$ and
$\sum_{\substack{k\in\OOO_K\ssm\{0,1\}\\ |k|\le N}}
\nod_K(k)\,\nod_K(k-1)$ have the same asymptotic behaviour, Theorem
\ref{theo:mainarithintro} in the introduction follows by taking
$N=\lfloor \sqrt{X}\rfloor$ and canceling the first and last factors
$\frac{1}{|\OOO_K^\times|^2}$ from Equation
\eqref{eq:Motohashigeneral}.
\cqfd

\brema\label{rem:Moto}{\rm We take $K=\QQ(i)$ in these remarks. In
  this case, $D_K=-4$ and $|\OOO_K^\times|=4$, and Equation
  \eqref{eq:Motohashigeneral} becomes
$$
\sum_{\substack{k\in\OOO_K\sm\{0,1\}\\|k|\le N}}\nod_K(k)\,\nod_K(k-1) =
\frac{4\pi^3}{\zeta_K(2)}(\ln N)^2N^2+\bigO((\ln N)N^2)\;.
$$

\noindent
\hypertarget{Moto1}{(1)} Figure
\addtocounter{fig}{1}\hyperlink{Bianchicommperp}{\arabic{fig}}\addtocounter{fig}{-1}
shows $\Ga_{\!\OOO\!_K}$-translates of $\wt \Delta$ with points at
infinity in the sector of $\CC$ defined by the inequalities $\Re
z\ge0$ and $\Im z\ge0$, that are at distance at most $\arcosh(9)$ from
$\wt \Delta$. These translates correspond to $|ad|\le 5$ in the
notation of the above proof of Theorem \ref{theo:mainarithintro}. The
grey surface in Figure
\addtocounter{fig}{1}\hyperlink{Bianchicommperp}{\arabic{fig}}\addtocounter{fig}{-1}
is a truncated sector of the boundary of the
$\arcosh(9)$-neighbourhood of $\wt \Delta$.  Four of the six
$\PSLZ$-translates in $\hdr$ at distance $\arcosh(9)$ from $\wt
\Delta$, shown in Figure
\addtocounter{fig}{-1}\hyperlink{autocommperpDelta}{\arabic{fig}}\addtocounter{fig}{1}
after Lemma \ref{lem:travaildivgeodmodularbis}, are now visible as the
red arcs in the foreground of Figure
\addtocounter{fig}{1}\hyperlink{Bianchicommperp}{\arabic{fig}}\addtocounter{fig}{-1}
with the points of intersection with the grey surface marked with red
points.

\begin{center}
\begin{overpic}[trim=160 0 10 0, clip,width=0.95\textwidth,tics=10]
{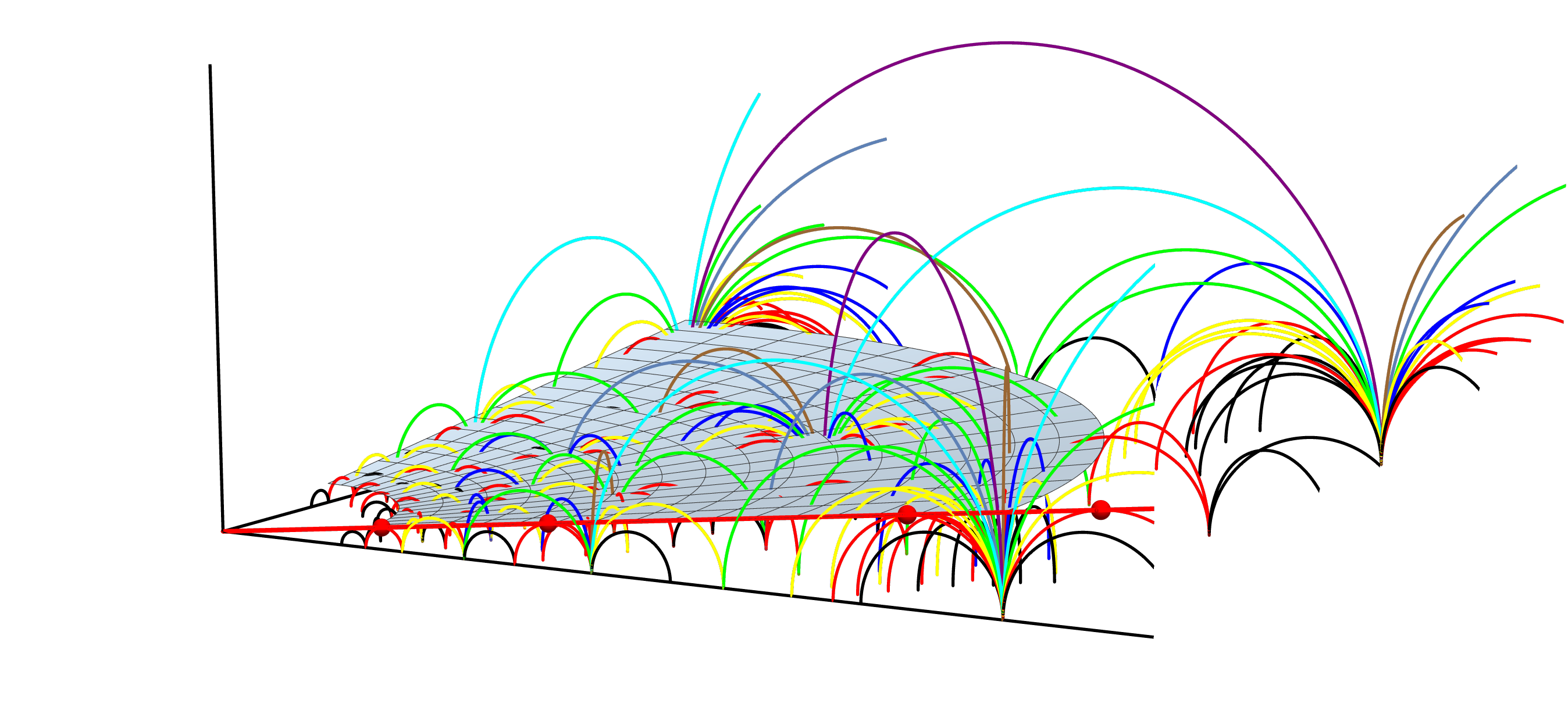}
  \put (1.5,11) {$0$}
 \put (19,8.5) {$\frac 13$}
  \put (28.3,7.5) {$\frac 12$}
   \put (38,6.5) {$\frac 23$}
  \put (58.5,4.5) {$1$}
\end{overpic}
\end{center}
\vspace*{-1.3cm}\addtocounter{fig}{1}
\begin{center}
{\small Figure \arabic{fig}~: \hypertarget{Bianchicommperp}{Translates
    of $\wt\Delta$ by the Gaussian Bianchi group $\PSL_2(\ZZ[i])$}.}
\end{center}

\renewcommand{\Motonod}{d}
\noindent \hypertarget{Moto2}{(2)} In this remark, we relate Theorem
\ref{theo:mainarithintro} with Motohashi's conjecture stated in the
third centered formula on page 277 of \cite{Motohashi01}, starting by
recalling the relevant definitions.\footnote{See the line after
Equation (9.7) in loc.~cit., that says that Motohashi's divisor
function $d$ is exactly our $d_{\QQ(i)}$ that we define below.} Let
$\I^+_K$ be the set of nonzero (integral) 
ideals of $\OOO_K$, and let 
$\Nr: \I^+_K\ra\NN\ssm\{0\}$ be the ideal norm, defined for every
$\aaa\in \I^+_K$ by $\Nr(\aaa)=[\OOO_K:\aaa]$. Let $\Motonod_K:
\I^+_K\to\NN$ be the {\em number of divisors function} of nonzero
ideals of $\OOO_K$, defined by $\Motonod_K(\aaa)=\card\{\bbb \in
\I^+_K:\bbb\mid\aaa\}$, so that $\zeta_K(s)^2=\sum_{\aaa\in \I^+_K}
\frac{\Motonod_K(\aaa)} {\Nr(\aaa)^s}$ for $\Re s>1$. For every
$x\in\OOO_K\ssm\{0\}$, let $\Nr(x)=\Nr(x\OOO_K)=|x|^2$ be the algebraic
norm and $\Motonod_K(x) =\Motonod_K(x\OOO_K)$.  Note that when
$\OOO_K$ is principal, and in particular when $K= \QQ(i)$, we have
$\Motonod_K(x)=\frac{\nod_K(x)} {|\OOO_K^\times|}$.  Theorem
\ref{theo:mainarithintro} when $K=\QQ(i)$ becomes
\begin{equation}\label{eq:motohashiconj}
\sum_{x\in\OOO_K\,\ssm\,\{0,-1\}\,:\; \Nr(x)\le X} \Motonod_K(x)\,\Motonod_K(x+1)
=\frac 1{16}\frac{\pi^3}{\zeta_K(2)}\; X(\ln X)^2+\bigO(X \ln X)\,.
\end{equation}

This confirms Motohashi's conjecture up to the usual multiplicative
factor $\frac 1{16}=\frac 1{|\OOO_K^\times|^2}$ in the particular case
when $f=1$ in his notation. See also \cite{SavVar03} for a similar
result on
\[
\sum_{k\in\OOO_K\ssm\{0,-f\}\,,\; \Nr(k)\le X}\Motonod_K(k)\,\Motonod_K(k+f)
\]
for all $f$, where the constant term in front of $X(\ln X)^2$ in
Equation \eqref{eq:motohashiconj} appears in a more complicated form
than above. See for instance \cite{GorNev15} for related counting
problems of integral points on homogeneous affine algebraic varieties,
as the one in $\M_2(\RR)$ defined by the equation $\det Y=f$ in the
variable $Y\in\M_2(\RR)$ for a fixed $f\in\ZZ$.}

\renewcommand{\Motonod}{d}
\smallskip\noindent (3) Numerical computations of the
ratio
\[
R(N)=\frac{4\,\zeta_K(2)}{\pi^3\;N^2(\ln N)^2}\;
\sum_{k\in\OOO_K\ssm\{0,-1\},\;|k|\le N} \Motonod_K(k)\,\Motonod_K(k+1)\,,
\]
for $K=\QQ(i)$ show that
\[
R(N)\simeq 1.213, 1.195, 1.18 \;\text{and}\; 1.167 \quad\text{when}\quad
N=2000, 4000, 8000\;\text{and}\;16000
\]
respectively.  This slow convergence of $N(R)$ to $1$ as $N\ra+\infty$
is similar to the case of integers: In Equation \eqref{eq:badpz}, the
ratio
\[
\frac{ \sum_{k=1}^n\nod(k)\,\nod(k+1)
  }{\frac 6{\pi^2}\, n(\ln n)^2}
\]
is approximately $1.18$ when $n=10^6$, and the ratio decreases closer
to $1$ very slowly, since the second term order in the development
differs from the main term order only by a logarithmic term.  On the
other hand, the ratio
\begin{equation}\label{eq:Z-ratio}
  \frac{ \sum_{k=1}^n\nod(k)\,\nod(k+1)}
       {\frac 6{\pi^2}\, n(\ln n)^2+a_1\,n\ln n}
\end{equation}
is approximately $0.997$ when $n=10^6$, giving already a much better
approximation.

The values of the ratio 
\[
\frac{\sum_{k\in\OOO_K\ssm\{0,-1\},\;|k|\le N} \Motonod_K(k)\,\Motonod_K(k+1)}
     {\frac{\pi^3}{4\,\zeta_K(2)}\;N^2(\ln N)^2+8.37\,N^2\ln N}
\]
analogous to the one in Equation \eqref{eq:Z-ratio} at $N=2000, 4000,
8000,16000$ are $1.00016$, $1.0001$, $1.00001$ and $0.99938$.  This
leads one to speculate a development similar to Equation
\eqref{eq:badpz}
\[
\sum_{k\in\OOO_K\ssm\{0,-1\},\;|k|\le N} \Motonod_K(k)\,\Motonod_K(k+1)=
\frac {\pi^3}{4\,\zeta_K(2)}\;N^2(\ln N)^2+ A_1N^2\ln N
+\smallo(N^2\ln N)\,,
\]
with $A_1\approx 8.4$.  
\erema

{\small \bibliography{../biblio} }

\bigskip
{\small
\noindent \begin{tabular}{l} 
Department of Mathematics and Statistics, P.O. Box 35\\ 
40014 University of Jyv\"askyl\"a, FINLAND.\\
{\it e-mail: jouni.t.parkkonen@jyu.fi}
\end{tabular}

\medskip\noindent\begin{tabular}{l}
Laboratoire de mathématique d'Orsay, UMR 8628 CNRS,\\
Universit\'e Paris-Saclay, 91405 ORSAY Cedex, FRANCE\\
{\it e-mail: frederic.paulin@universite-paris-saclay.fr}
\end{tabular}
}

\end{document}